\pdfoutput=1
\RequirePackage{ifpdf}
\ifpdf 
\documentclass[pdftex]{sigma}
\else
\documentclass{sigma}
\fi

\usepackage{amscd,euscript}
\usepackage{mathrsfs}
\usepackage{amsxtra}

\usepackage[all]{xy}

\usepackage{verbatim,epsf}
\usepackage{framed}
\usepackage{multirow}
\usepackage{array}
\usepackage{tabu}
\usepackage{mathtools}

\newcommand{\field}[1]{\ensuremath{\mathbb{#1}}}
\newcommand{\CC}{\field{C}}

\newcommand{\PP}{\field{P}}

\newcommand{\RR}{\field{R}}

\newcommand{\ZZ}{\field{Z}}
\newcommand{\calO}{\mathcal{O}}

\newcommand\sbul[1][.6]{\mathbin{\vcenter{\hbox{\scalebox{#1}{$\bullet$}}}}}

\newcommand{\curly}[1]{\mathscr{#1}}

\newcommand{\cC}{\curly{C}}

\newcommand{\cE}{\curly{E}}
\newcommand{\cF}{\curly{F}}

\newcommand{\cK}{\curly{K}}
\newcommand{\cL}{\curly{L}}
\newcommand{\cM}{\curly{M}}
\newcommand{\cN}{\curly{N}}
\newcommand{\cO}{\curly{O}}

\newcommand{\cS}{\curly{S}}

\newcommand{\cU}{\curly{U}}
\newcommand{\cV}{\curly{V}}
\newcommand{\cW}{\curly{W}}

\newcommand{\qB}{\mathcal{B}}
\newcommand{\qC}{\mathcal{C}}
\newcommand{\qH}{\mathcal{H}}
\newcommand{\qI}{\mathcal{I}}
\newcommand{\qR}{\mathcal{R}}
\newcommand{\qW}{\mathcal{W}}

\newcommand{\qp}{\boldsymbol{\cdot}}

\DeclareMathOperator{\Ad}{Ad}
\DeclareMathOperator{\Aut}{Aut}
\DeclareMathOperator{\Div}{Div}
\DeclareMathOperator{\End}{End}
\DeclareMathOperator{\pardeg}{\mathrm{par}\deg}
\DeclareMathOperator{\parmu}{\mathrm{par}\,\mu}
\DeclareMathOperator{\Res}{Res}
\DeclareMathOperator{\tr}{tr}
\DeclareMathOperator{\GL}{GL}

\numberwithin{equation}{section}

\newtheorem{Theorem}{Theorem}[section]
\newtheorem*{Theorem*}{Theorem}
\newtheorem{Corollary}[Theorem]{Corollary}
\newtheorem{Lemma}[Theorem]{Lemma}
\newtheorem{Proposition}[Theorem]{Proposition}
 { \theoremstyle{definition}
\newtheorem{Definition}[Theorem]{Definition}

\newtheorem{Remark}[Theorem]{Remark} }

\begin{document}
\allowdisplaybreaks

\renewcommand{\thefootnote}{}

\newcommand{\arXivNumber}{2012.13389}

\renewcommand{\PaperNumber}{062}

\FirstPageHeading

\ShortArticleName{Geometric Models and Variation of Weights on Moduli of Parabolic Higgs Bundles}

\ArticleName{Geometric Models and Variation of Weights on Moduli\\ of Parabolic Higgs Bundles over the Riemann Sphere:\\ a Case Study\footnote{This paper is a~contribution to the Special Issue on Mathematics of Integrable Systems: Classical and Quantum in honor of Leon Takhtajan.

~~\,The full collection is available at \href{https://www.emis.de/journals/SIGMA/Takhtajan.html}{https://www.emis.de/journals/SIGMA/Takhtajan.html}}}

\Author{Claudio MENESES}

\AuthorNameForHeading{C.~Meneses}

\Address{Mathematisches Seminar, Christian-Albrechts Universit\"at zu Kiel,\\ Heinrich-Hecht-Platz 6, 24118 Kiel, Germany}
\Email{\href{mailto:meneses@math.uni-kiel.de}{meneses@math.uni-kiel.de}}
\URLaddress{\url{http://www.math.uni-kiel.de/geometrie/de/claudio-meneses/}}

\ArticleDates{Received September 21, 2021, in final form July 28, 2022; Published online August 13, 2022}

\Abstract{We construct explicit geometric models for moduli spaces of semi-stable strongly parabolic Higgs bundles over the Riemann sphere, in the case of rank two, four marked points, arbitrary degree, and arbitrary weights. The mechanism of construction relies on elementary geometric and combinatorial techniques, based on a detailed study of orbit stability of (in general non-reductive) bundle automorphism groups on certain carefully crafted spaces. The aforementioned techniques are not exclusive to the case we examine, and this work elucidates a general approach to construct arbitrary moduli spaces of semi-stable parabolic Higgs bundles in genus~0, which is encoded into the combinatorics of weight polytopes. We also present a comprehensive analysis of the geometric models' behavior under variation of parabolic weights and wall-crossing, which is concentrated on their nilpotent cones.}


\Classification{14H60; 14D22; 32G13; 22E25}

\begin{flushright}
\begin{minipage}{88mm}
\it To Professor Leon A.~Takhtajan on his 70th birthday,\\ with profound respect and admiration.
\end{minipage}
\end{flushright}

\renewcommand{\thefootnote}{\arabic{footnote}}
\setcounter{footnote}{0}

\section{Introduction}

The moduli spaces of parabolic Higgs bundles on a Riemann surface of lowest possible dimension are elliptic surfaces in virtue of their Hitchin fibrations \cite{Hit87a}. These examples occur in low genus and encode at once, in its simplest possible form, the characteristic features of the rich and intricate geometry occuring in arbitrary dimensions. From the general theory of elliptic surfaces~\cite{Mir89}, they are necessarily biholomorphic to one of the fibrations in Kodaira's list. Following the seminal work of Gaiotto--Moore--Neitzke~\cite{GMN13}, their hyperk\"ahler geometry has recently attracted renewed interest in relation to the study of gravitational instantons of ALG type~\cite{FMSW22}.

The \textit{toy model} \cite{Hau98} for the Hitchin fibration and $\CC^{*}$-action of the moduli spaces of rank 2 parabolic Higgs bundles over $\CC\PP^{1}$ with four marked points, is the elliptic surface $\cM_{\mathrm{toy}}$ associated to a choice $(z_{1},z_{2},z_{3},z_{4})\in\mathrm{Conf}_{4}(\CC\PP^{1})$ as follows. Let $\Sigma_{D}\rightarrow \CC\PP^{1}$ be the elliptic curve of $D=z_{1}+z_{2}+z_{3}+z_{4}$, with involution $\tau$ and ramification points $\{w_{1},w_{2},w_{3},w_{4}\}$. The extension of $\tau$ to an action in $\Sigma_{D}\times \CC$ in terms of the character generated by $\rho(\tau) = -1$ determines a~2-dimensional orbifold $(\tau\times\rho) \backslash (\Sigma_{D}\times \CC)$ which is smooth away of the points $[(w_{i},0)]$. $\cM_{\mathrm{toy}}$~is defined as the desingularization of $(\tau\times\rho) \backslash (\Sigma_{D}\times \CC)$ by blow-up at those points. The elliptic fibration $\pi_{2}\colon\cM_{\mathrm{toy}}\rightarrow \CC$ is given by extending the map $[(w,z)]\mapsto z^{2}$ trivially along exceptional divisors, while its $\CC^{*}$-action is induced by the standard $\CC^{*}$-action on $\CC$. The loci where the (holomorphic) involution extending $\tau\times\rho$ acts freely is a Zariski open subset $\cU\cong T^{*}\CC\PP^{1}\backslash\CC\PP^{1}$, equipping $\cM_{\mathrm{toy}}$ with an additional projection $\pi_{1}\colon \cU\rightarrow\CC\PP^{1}$. Intuitively, a biholomorphism between $\cM_{\mathrm{toy}}$ and the aforementioned moduli spaces of semi-stable parabolic Higgs bundles could be constructed by relating $\pi_{1}$ to the forgetful map to parabolic structures, and $\pi_{2}$ to the determinant map on parabolic Higgs fields.

There are some caveats when trying to carry out the previous biholomorphisms, though. First of all, they are not expected to hold for other ranks and number of marked points.
A nontrivial feature of parabolic Higgs bundles is the dependence of their moduli problem on a choice of parabolic weights.
While several classical results on this dependence \cite{Nak96, Thad02} are linked to the study of the birational geometry of the moduli problem, not only they are trivial when the moduli space is a complex surface, but also they shed no light on other subtle phenomena. For instance, it is unclear how holomorphic invariants such as Harder--Narasimhan stratifications may depend on parabolic weights, or what happens to the moduli space when semi-stable parabolic bundles don't exist \cite{Bau91, Bis98}.\footnote{Incidentally, in \cite{Ray17} Rayan introduced the notion of ``quiver at the bottom of the nilpotent cone", in the context of twisted Higgs bundles on $\CC\PP^{1}$. The idea that for parabolic Higgs bundles, the subspace parametrized by these objects may not correspond to a moduli space of parabolic bundles, seems to be standard folklore \cite{Bla15}, and resurfaces in our work in an explicit way.} Moreover, the natural hyperk\"ahler structures on these moduli spaces depend nontrivially on parabolic weights \cite{KW18, MT21}, and the \textit{hyperk\"ahler moduli problem} motivates the search for suitable complex-analytic structures capturing this dependence. The will of attacking these problems is the main motivation to carry out this work. More specifically, this work sprung as a natural advance towards the parabolic Higgs bundle generalization of the results in \cite{MT21}, dealing with the study of K\"ahler forms on moduli spaces of semi-stable parabolic bundles, where the understanding of refined moduli space stratifications is crucial.

We propose a construction of moduli spaces of rank 2 stable parabolic Higgs bundles on $\CC\PP^{1}$ with four marked points, on which the dependence on parabolic weights and its wall-crossing phenomena becomes as transparent as possible. The construction is comprised of two steps. The first step is based on the observation that the peculiarities of genus 0 (namely, the infinitesimal rigidity of holomorphic bundles dictated by the Birkhoff--Grothendieck theorem) lead to the construction of Harder--Narasimhan strata for underlying holomorphic bundles $E\cong\calO(m_{1})\oplus\calO(m_{2})$ as semi-stable orbit spaces under the action of bundle automorphisms, since any isomorphism class of quasi-parabolic Higgs bundles on $E$ can be modeled as an orbit in an explicit complex manifold $\mathrm{QPH}(E)$ with respect to an action of the group $\PP(\Aut(E))$ of projective automorphisms of $E$.\footnote{We remark that bundle automorphism groups are in general non-reductive. It should be possible to reinterpret our constructions as explicit examples of variations of non-reductive GIT quotients \cite{BDHK18, BJK18} (cf.~\cite{Ham19}).} The mechanism of construction is entirely self-contained, and follows the general localization philosophy to study invariants of bundle automorphism actions proposed in \cite{Men18}. Once this is done, the remaining problem is to understand the way these strata can be subsequently glued into a smooth complex manifold modeling the moduli space.

An important property of our construction is that it admits a generalization to arbitrary rank and number of marked points.
This way, this is the first installment (or ``toy model'') of a~project focused on the study of complex-analytic structures of moduli spaces of parabolic Higgs bundles (which could be understood as an instance of the so-called \textit{Gelfand principle}), where we exhaustively explore the simplest possible nontrivial example, and we rigorously exhibit the independence of the biholomorphism type of these moduli spaces under variation of parabolic weights and wall-crossing. Since the latter is exclusive of the case of four marked points, this is yet another reason to present it separately.

\subsection{Structure of the paper}
Even
though the techniques employed in this work are all elementary, in the end the total number of steps involved in the construction of geometric models is substantial. This requires an adequate signposting strategy to guide the reader through it, which we now provide. The starting point are the model spaces $\mathrm{QPH}(E)$ of quasi-parabolic Higgs bundles on $E$, together with a pair of projections
\[
\xymatrix{
&\mathrm{QPH}(E) \ar[ld]_{\mathrm{par}} \ar[rd]^{\det} &\\
\displaystyle\mathrm{QP}(E) := \prod_{i=1}^{4}\PP(E\vert_{z_{i}}) && H^{0}\bigg(\CC\PP^{1},K^{2}_{\CC\PP^{1}}\bigg(\displaystyle\sum_{i=1}^{4}z_{i}\bigg)\bigg)
}
\]
$\mathrm{par}$ is the $\PP(\Aut(E))$-equivariant forgetful map to quasi-parabolic structures, while $\det$ is the $\PP(\Aut(E))$-invariant determinant map on parabolic Higgs fields.

The construction and classification of all different Harder--Narasimhan strata as spaces of semi-stable $\PP(\Aut(E))$-orbits in the complex manifolds $\mathrm{QPH}(E)$ follows from a systematic break down of the semi-stability condition for parabolic Higgs bundles into structural units, recollected in Sections~\ref{sec:ev-data-bundles}--\ref{sec:CS}. This process consists of four parts which are built progressively:

Section~\ref{sec:ev-data-bundles}:
$\mathrm{QP}(E)$ is stratified in terms of interpolation loci for line sub-bundles of $E$. The resulting stratification is then related to the stratification by pointwise $\PP(\Aut(E))$-stabilizer subgroups.

Section~\ref{sec:struct}:
$\det$ stratifies $\mathrm{QPH}(E)$ into a (possibly empty) bulk $\mathrm{QPH}_{\CC^{*}}(E)$ and its nilpotent locus $\mathrm{QPH}_{0}(E)$. Potentially unstable points can only occur in $\mathrm{QPH}_{0}(E)$ (Lemma~\ref{lemma:nilpotent}). This leads to a stratification of $\mathrm{QPH}_{0}(E)$ in terms of the interpolation properties of potentially destabilizing line sub-bundles, which arise as kernels of parabolic Higgs fields on $E$ when the latter are nonzero. The stratifications of $\mathrm{QPH}(E)$ and $\mathrm{QP}(E)$ are then related by the projection $\mathrm{par}$ (Proposition~\ref{prop:QPH-par}).

Section~\ref{sec:poly}:
A combinatorial classification of semi-stability walls $\qH$ and open chambers $\qC$ in the even and odd weight polytopes $\qW_{\sbul}$ is explained. Open chambers are defined as the connected components in the complement of the union of all semi-stability walls in $\qW_{\sbul}$, and parametrize the different stable loci $\mathrm{QPH}^{s}_{\qC}(E)$ in $\mathrm{QPH}(E)$. Open chambers are split into \textit{interior} and \textit{exterior} (denoted as $\qC_{\qI}$ and $\qC_{I}$ respectively; this notation is carefully explained), according to whether stable parabolic bundles exist or not.

Section~\ref{sec:CS}: In order to determine stable loci $\mathrm{QPH}^{s}_{\qC}(E)$, the stratifications of $\mathrm{QPH}_{0}(E)$ and the combinatorics of $\qW_{\sbul}$ are combined to classify conditionally stable loci (Definition~\ref{def:cond-stab}) as well as their (in general non-Hausdorff) orbit spaces. The latter are organized in terms of the combinatorics of some classical projective plane configurations, and then dissected into \textit{basic building blocks} (Definition~\ref{def:b-b-b}). Finally, the structure of basic building blocks is classified (they are biholomorphic to $\CC\PP^{1}$ or $\CC$), and a \textit{nilpotent cone assembly kit} (Definition~\ref{def:a-kit}) is conformed for every choice of open chamber $\qC\subset \qW_{\sbul}$.

Once all Harder--Narasimhan strata are constructed, their glueing is brought in, and the interpretation of semi-stability via the geometry of $\PP(\Aut(E))$-actions leads to an explicit classification of wall-crossing as nilpotent cone transformations, in terms of the combinatorics of~$\qW_{\sbul}$, as well as an additional characterization of parabolic $S$-equivalence along semi-stability walls. We list our main results, whose proof is presented in Section~\ref{sec:proof-m} together with a series of corollaries.

\begin{Theorem}[orbit space models for Harder--Narasimhan strata]\label{theo:main-1}
For any open chamber \mbox{$\qC\subset\qW_{\sbul}$} and vector bundle $E$, if $\mathrm{QPH}^{s}_{\qC}(E)\neq \varnothing $, then $\PP(\Aut(E))$ acts freely and properly on~it. There is an isomorphism
\[
\cM_{\qC}(E)\cong \mathrm{QPH}^{s}_{\qC}(E)/\PP(\Aut(E)),
\]
where $\cM_{\qC}(E)$ is the Harder--Narasimhan stratum associated to $E$. $\cM_{\qC}(E)$ is smooth if it is $2$-dimensional, and has at most two nodal singularities otherwise.
\end{Theorem}

\begin{figure}[!ht]
\centering
\includegraphics[width=2.8in]{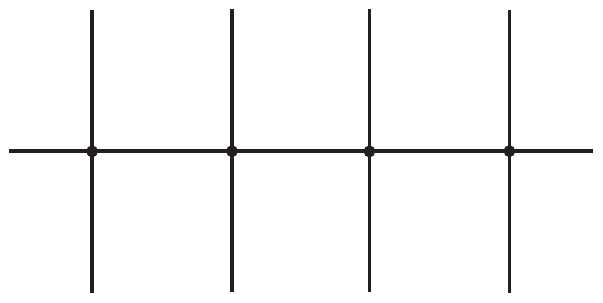}
\put(-169,54){\makebox(0,0)[lb]{\small$z_1$}}
\put(-122,54){\makebox(0,0)[lb]{\small$0$}}
\put(-76,54){\makebox(0,0)[lb]{\small$1$}}
\put(-29,54){\makebox(0,0)[lb]{\small$\infty$}}
\put(0,50){\makebox(0,0)[lb]{\small $\CC\PP^{1}$}}
\caption{The $D_{4}$-configuration.}\label{fig:D4-conf}
\end{figure}

\begin{Theorem}[glueing of Harder--Narasimhan strata]\label{theo:main-2}\quad
\begin{enumerate}\itemsep=0pt
\item[$(i)$] For $d$ even, the two orbit spaces
\[
\mathrm{QPH}_{\CC^{*}}(E)/\PP(\Aut(E)),\qquad m_{1}-m_{2} = 0,2,
\]
glue into a $2:1$-branched covering over $\CC\PP^{1}\times \CC^{*}$, ramified over $\{z_{1},0,1,\infty\}\times\CC^{*}$.

\item[$(ii)$] For every choice of open chamber $\qC\subset \qW_{\sbul}$, the components in the nilpotent cone assembly kit are glued into a $D_{4}$-configuration $($Figure~$\ref{fig:D4-conf})$. Moreover, there is an isomorphism
\[
\cM_{\qC}\cong \cM_{\mathrm{toy}},
\]
where $\cM_{\qC}$ is the smooth complex surface resulting from the glueing of geometric models of Harder--Narasimhan strata of corresponding degree.
\end{enumerate}
\end{Theorem}

\begin{Theorem}[wall-crossing classification; see Sections \ref{sec:poly} and~\ref{sec:CS} for notations]\label{theo:main-3}
For each boundary wall $\qH_{I,\vert I\vert -3}$ between the exterior and interior chambers $\qC_{I}$ and $\qC_{\qI(I)}$, there is an isomorphism $\cM_{\qC_{I}} \cong \cM_{\qC_{\qI(I)}}$ exchanging the nilpotent cone's central spheres $\cS_{I}$ and $\cS_{\qI(I)}$, and equal to the identity otherwise. For any pair of neighboring interior chambers $\qC_{\qI(I)}$ and $\qC_{\qI(I')}$, there is an isomorphism $\cM_{\qC_{\qI(I)}} \cong \cM_{\qC_{\qI(I')}}$ exchanging the compactified $I$- and $I'$-basic building blocks, and equal to the identity otherwise.
\end{Theorem}

After this, the Hitchin elliptic fibrations, their degeneration into their nilpotent cones, as well as a $\CC^{*}$-action and a collection of Hitchin sections, are explicitly identified as by-products of the method of construction. In Appendix~\ref{sec:Res-mod} we also introduce the alternative \textit{residue models} $\cM^{\mathrm{ev}}_{\qC}(E)$ for Harder--Narasimhan
strata.\footnote{Cf. Komyo--Saito \cite{KS19}, who provide explicit descriptions of Zariski open sets on moduli spaces of logarithmic connections and parabolic Higgs bundles in odd degree, by apparent singularities and their dual parameters as coordinates.} The generalization of this work to an arbitrary number of punctures requires a more careful analysis of automorphism group actions, including a suitable axiomatization of the intrinsic structures on display here, and will be treated elsewhere \cite{Mena}.

\subsection{Relation to other work}
The literature on moduli spaces emerging from parabolic structures on $\CC\PP^{1}$ is vast and comprises different perspectives, strategies, and objectives. We compiled works with ideas overlapping substantially with ours. Further details can be found in those works and references therein.

Related models were used by
Loray--Saito--Simpson \cite{LSS13}
in their proof of the non-abelian Hodge foliation conjecture for moduli of rank 2 logarithmic connections, in odd degree and $n=4$ marked points, in terms of transversal foliations arising as Lagrangian fibrations. There, a~bundle splitting type is necessarily generic, and their analysis of moduli dependence on parabolic weights leads to an analogous characterization of exterior chambers. These fibrations are studied for arbitrary $n$ in \cite{LS15} (cf.~\cite{KS19}), and in \cite{HL19} the study of moduli of holomorphic connections in genus~2 is reduced to the genus 0 case and $n=6$. The overall framework of these works is that of birational geometry. While similar in spirit, our proposal is considerably different in its methods and results.

Godinho--Mandini \cite{GM13,GM21} defined an isomorphism between hyperpolygon spaces and the Zariski open subset of rank 2 parabolic Higgs bundles of degree $0$ that are holomorphically trivial in terms of the residue data of the latter (a detailed description of the correspondence in the case $n=4$, in terms of the Nakajima quiver variety for the affine Dynkin diagram $\tilde{D_{4}}$, is discussed by Rayan--Schaposnik in \cite[Section 4.3]{RS21}). Similarly, Blaavand \cite{Bla15} provides a geometric construction of the nilpotent cone in degree $1$ and $n=4$ under a choice of parabolic weights constrained to lie inside a specific open chamber (see Remark~\ref{rem:Blaavand}). Heller--Heller \cite{HH16} applied the abelianization of logarithmic connections in degree $0$ when $n=4$ to compute symplectic volumes for moduli spaces of stable parabolic bundles.

\section{Conventions and index of notation}\label{sec:conv}

We will fix a point $z_{1}\in\CC\PP^{1}\backslash\{z_{2},z_{3},z_{4}\}$, where $z_{2}=[0:1]=0$, $z_{3}=[1:1]=1$, $z_{4}=[1:0]=\infty$. When necessary, we will consider the choice of cross-ratio such that $(z_{1},z_{2};z_{3},z_{4})=z_{1}$. Similarly, we will denote $D = \sum_{i=1}^{4}z_{i}$. A \emph{parabolic structure supported on $D$} on a rank 2 holomorphic vector bundle $E\rightarrow \CC\PP^{1}$ of degree $d$ consists of \emph{flags}
\[
F_{i} \subset E|_{z_{i}},\qquad i=1,2,3, 4
\]
weighted by real numbers $0\leq \beta_{i}<1$ such that $F_{i} = E\vert_{z_{i}}$ if and only if $\beta_{i}=0$. A \emph{quasi-parabolic structure} omits the choice of weights. Hereafter $E_{*}$ will denote a rank 2 \emph{parabolic bundle}, i.e., a~choice of parabolic structure over $D$ on $E$, and $E_{\qp}$ the underlying quasi-parabolic bundle. Given~$E_{*}$, to any subset $I\subset\{1,2,3,4\}$ we associate the weighted sum
\[
\beta_{I} = \sum_{i\in I}\beta_{i} - \sum_{j\in I^{\mathsf{c}}}\beta_{j}.
\]
For any $E_{\qp}$, a line sub-bundle $L\subset E$ determines a subset $I_{L,E_{\qp}}\subset \{1,2,3,4\}$ defined as
\[
I_{L,E_{\qp}} = \bigl\{i\in\{1,2,3,4\}\colon L|_{z_{i}} \subset F_{i}\bigr\}.
\]
For any $z\in\CC\PP^{1}$, $\mathfrak{n}(E\vert_{z})$ will denote the cone of nilpotent endomorphisms of $E\vert_{z}$, i.e., the singular quadric in $\mathfrak{sl}(E|_{z})$ defined as the zero locus of its determinant map, and for any $F\in\PP(E\vert_{z})$, $\mathfrak{n}(F)$ is the line of nilpotent endomorphisms $\phi$ such that $\ker(\phi) \supset F$.
Let $\End (E):=E^{\vee}\otimes E$. An~element $\Phi\in H^{0}(\CC\PP^{1},\End (E)\otimes K_{\CC\PP^{1}}(D))$ is a (\textit{strongly}) \textit{parabolic Higgs field} of $E_{*}$ if for each $i=1,2,3,4$, the residue $\Res_{z_{i}}\Phi$ satisfies
\[
\Res_{z_{i}}\Phi\in\mathfrak{n}(F_{i}).
\]
A \emph{parabolic Higgs bundle} is a pair of the form $(E_{*},\Phi)$. The pair $(E_{*},\Phi)$ is called \emph{stable} (resp.~\textit{semi-stable}) if for every $\Phi$-invariant line sub-bundle $L\subset E$ we have that
\begin{gather}\label{eq:effective}
\beta_{I_{L,E_{\qp}}} < d - 2\deg(L)\qquad (\text{resp.}~\leq).
\end{gather}
The (semi-)stability of a parabolic bundle $E_{*}$ is defined by letting $\Phi= 0$. $(E_{*},\Phi)$ is called \textit{strictly semi-stable} if it semi-stable but not stable.

Given a vector bundle $E$, a line bundle $L$, and their tensor product $E' = E\otimes L$, the bundle isomorphism $\End(E')\cong \End (E)$ induces a bijection between the sets of quasi-parabolic Higgs bundles on $E$ and $E'$. Since $L'$ is a line sub-bundle of $E'$ if and only if $L^{-1}\otimes L'$ is a line sub-bundle of $E$, it follows that the stability of a parabolic Higgs bundle is preserved under tensor product by a line sub-bundle.

\begin{Remark}\label{rem:effective}
In the case of rank 2 holomorphic vector bundles, the standard definition of parabolic stability \cite{MS80,Sim90} involves \emph{parabolic weights} $0\leq \alpha_{i1}\leq\alpha_{i2}<1$, $i=1,2,3,4$. In that case, the \emph{parabolic degree and slope} of $E_{*}$ are defined as
\[
\pardeg(E_{*}) := d + \sum_{i = 1}^{4}(\alpha_{i1}+ \alpha_{i2}),\qquad
\parmu (E_{*}) = \pardeg(E_{*})/2.
\]
In turn, a line sub-bundle $L\subset E$ acquires a set of parabolic weights $\alpha'_{i}$ defined as $\alpha_{i2}$ whenever $i\in I$ and $\alpha_{i1}$ otherwise. The induced (semi-)stability condition then reads
\[
\parmu (L_{*}):=\deg(L)+\sum_{i=1}^{4}\alpha'_{i} < \parmu(E_{*}) \qquad \text{(resp.~$\leq$)}.
\]
It can be verified that the last inequality only depends on the effective parameters
\[
\beta_{i} :=\alpha_{i2}-\alpha_{i1}\in[0,1),\qquad i=1,2,3,4,
\]
and reduces to \eqref{eq:effective}. A set of parabolic weights will be called \emph{admissible} whenever the equation $\pardeg(E_{*})=0$ is satisfied. The admissibility of parabolic weights, together with the stability condition for parabolic Higgs bundles, are necessary for the non-abelian Hodge correspondence to hold (see \cite{Sim90}), even though the moduli problem is meaningful without the former; e.g., in~\cite{Muk03} Mukai considers parabolic weights corresponding to the choices $\alpha_{i1}=0$ and $\alpha_{i2}=\beta_{i}$, for $n$ marked points on $\CC\PP^{1}$. On the other hand, Bauer \cite{Bau91} considers the $\mathrm{SU}(2)$-constraints
\[
\pardeg(E_{*})=0,\qquad 0<\alpha_{i1}=\alpha_{i}< 1/2,\qquad \alpha_{i2} = 1-\alpha_{i},\qquad i=1,\dots,n,
\]
fixing the value of $d$. When $d \equiv n\; (\bmod\; 2)$, up to the tensor product of $E$ by a suitable line bundle, every set $\{\beta_{i}\}$ can be lifted to a unique set of parabolic weights satisfying the $\mathrm{SU}(2)$-constraints. When $d \not\equiv n\; (\bmod\; 2)$, a lift to a set of admissible parabolic weights is only possible under further numerical constraints, but a weaker lift satisfying the $\mathrm{SU}(2)$-constraints can still be achieved in terms of an extra marked point $z_{0}$ with parabolic weights $\alpha_{01}=\alpha_{02}=1/2$ (so that $\beta_{0}=0$). This operation leaves the moduli problem unchanged, since it forces any flag over~$E|_{z_{0}}$ to be trivial and the corresponding residue $\Res_{z_{0}}\Phi$ to vanish.\footnote{An operation on parabolic bundles discussed in the literature is the \emph{parabolic tensor product}. We don't consider it here as it does not preserve the underlying associated bundle $\End (E)$.}
The alternative convention $-1/2<-\alpha'_{i}< \alpha'_{i}< 1/2$ for the $\mathrm{SU}(2)$-constraints is related to the former in terms of the transformations $\beta'_{i}= 1-\beta_{i}$, $i=1,2,3,4$.
\end{Remark}

\section[Geometry of P(Aut(E))-actions and stratifications in QP(E)]
{Geometry of $\boldsymbol{\PP(\Aut(E))}$-actions and stratifications in $\boldsymbol{\mathrm{QP}(E)}$}\label{sec:ev-data-bundles}

Let $E\rightarrow \CC\PP^{1}$ be a holomorphic vector bundle of rank 2 and degree $d$. A basic holomorphic invariant associated to $E$ is its \textit{Harder--Narasimhan filtration}. In genus 0, the Birkhoff--Grothendieck theorem postulates the existence of an isomorphism of vector bundles
\[
E\cong \calO(m_{1})\oplus\calO(m_{2}),\qquad m_{1}\geq m_{2},
\]
in such a way that $d=m_{1}+m_{2}$. Every rank 2 bundle $E$ is either semi-stable (i.e., it splitting coefficients satisfy $m_{1}=m_{2}$) or otherwise admits a unique filtration of the form
\[
E_{1}\subset E_{2} = E,\qquad \deg(E_{1}) >\deg(E_{2})/2.
\]
In other words, every isomorphism $E\cong \calO(m_{1})\oplus\calO(m_{2})$ is a refinement of the Harder--Nara\-sim\-han filtration of $E$. When $E$ is not semi-stable, we have that $E_{1}\cong \calO(m_{1})$, and a choice of Birkhoff--Grothendieck splitting for $E$ is equivalent to a choice of complementary sub-bundle $\calO(m_{2})\subset E$. The integers $m_{1}\leq m_{2}$ are holomorphic invariants for $E$. $E$ is called \emph{evenly-spit} if $m_{1}-m_{2}\leq 1$, or equivalently if $H^{1}(\CC\PP^{1},\End (E))=0$, where $\End(E):= E\otimes E^{\vee}$. Consequently, on any holomorphic family $\cF\rightarrow \mathrm{B}\times\CC\PP^{1}$ of rank 2 and degree $d$ holomorphic vector bundles $E\rightarrow\CC\PP^{1}$, the existence of evenly-split bundles is an open condition.

Pointwise multiplication endows the space $H^{0}(\CC\PP^{1},\End (E))$ with an associative algebra structure. The group of global bundle automorphisms of $E$ is the Lie group of invertible ele\-ments
\[
\Aut(E)\subset H^{0}\bigl(\CC\PP^{1},\End (E)\bigr).
\]
Its projectivization is defined as $\PP(\Aut (E)):=\Aut(E)/\mathrm{Z}(\Aut(E))$, where its center $\mathrm{Z}(\Aut(E))$ consists of all nonzero multiples of the identity. $\Aut(E)\cong\mathrm{GL}(2,\CC)$ in the special case when $E$ is semi-stable. Otherwise, when $m_{1}>m_{2}$, $\Aut(E)$ is $(m_{1}-m_{2}+3)$-dimensional and preserves the Harder--Narasimhan filtration of $E$. Its \emph{unipotent radical} is the normal subgroup of unipotent automorphisms of $E$, and will be denoted by $\mathrm{R}(\Aut(E))$. $\mathrm{R}(\Aut(E))$ is maximally abelian, has codimension 2 in $\Aut(E)$, and acts freely and transitively on the space $\mathrm{L}_{m_{2}}(E)$. In particular, the subgroup $\mathrm{R}(\Aut(E))\cap \mathrm{Z}(\Aut(E))$ is trivial.

Consider a point $z\in\CC\PP^{1}$. For any $F\in\PP(E\vert_{z})$, let $\mathrm{P}(F)\subset\mathrm{GL}(E\vert_{z})$ be the parabolic subgroup stabilizing the flag $F\subset E\vert_{z}$, and $\mathrm{R}(F)$ its unipotent radical. The evaluation maps
\[
\Aut(E)\mapsto \Aut(E)\vert_{z}\subset \mathfrak{gl}(E\vert_{z})
\]
determine a subgroup of $\GL(E|_{z})$. This subgroup coincides with $\GL(E\vert_{z})$ when $m_{1} = m_{2}$, in which case the evaluation map is an isomorphism. Otherwise, when $m_{1} > m_{2}$, $\Aut(E)\vert_{z}$ coincides with the parabolic subgroup $\mathrm{P}(E_{1}\vert_{z})$ that is induced by the Harder--Narasimhan filtration of $E$.
In this case we define the affine line
\[
\cV_{z}(E) := \PP(E\vert_{z})\backslash\{E_{1}\vert_{z}\}.
\]
Since there is an isomorphism
$\Aut(E)|_{z} \cong \mathrm{P}(E_{1}\vert_{z})$,
and $\PP(\mathrm{R}(E_{1}\vert_{z}))\cong \CC$ acts freely and transitively on $\cV_{z}(E)$, the cell decomposition
\[
\PP(E\vert_{z}) = \{E_{1}\vert_{z}\}\sqcup\cV_{z}(E)
\]
is the partition into two orbits determined by the action of $\PP(\Aut(E)|_{z})$ on $\PP(E|_{z})$.

\subsection{Line sub-bundles and interpolation of quasi-parabolic structures}

For any $E$, we will denote the space of all line sub-bundles $L\subset E$ of degree $j$ by $\mathrm{L}_{j}(E)$. Given an isomorphism $E\cong \calO(m_{1})\oplus\calO(m_{2})$, $\mathrm{L}_{j}(E)$ is modeled by the Zariski open set
\[
\mathrm{L}_{j}(\calO(m_{1})\oplus\calO(m_{2})) \subset \PP\bigl(H^{0}\big(\CC\PP^{1},\calO(m_{1}-j)\big)\oplus H^{0}\big(\CC\PP^{1},\calO(m_{2}-j)\big)\bigr)
\]
generated by pairs of holomorphic sections $(s_{1},s_{2})$ that are not simultaneously zero and have disjoint zero sets. In the particular case when $m_{1}=m_{2}=m$, we have the isomorphisms
\[
\mathrm{L}_{m}(E)\cong\CC\PP^{1},\qquad \mathrm{L}_{m-1}(E)\cong\CC\PP^{3}\backslash \mathrm{Segre}\big(\CC\PP^{1}\times\CC\PP^{1}\big).
\]
On the other hand, when $m_{1} > m_{2}$ there are no line sub-bundles $L\subset E$ of degree $j$ for $m_{1} > j > m_{2}$, and
\[
\mathrm{L}_{m_{2}}(E)\cong \mathbb{A}^{m_{1}-m_{2}+1}\subset \CC\PP^{m_{1}-m_{2}+1}.
\]
Given $I\subset\{1,2,3,4\}$ and $j\leq m_{2}$, we define the $\calO(j)$-interpolation locus on $\mathrm{QP}(E)$, relative to~$I$,~as
\begin{gather}\label{eq:B_{I,j}(E)}
\mathrm{B}_{I,j}(E) := \bigl\{(F_{1},F_{2},F_{3},F_{4})\in\mathrm{QP}(E)\colon \exists L\in\mathrm{L}_{j}(E),\; F_{k} = L\vert_{z_{k}}\; \forall k\in I \bigr\},
\end{gather}
which satisfies
\begin{gather}\label{eq:I<I'}
\mathrm{B}_{I',j}(E)\subset \mathrm{B}_{I,j}(E)\qquad \text{whenever}\quad I\subset I'.
\end{gather}
For any vector bundle $E$, there is a natural action of $\PP(\Aut(E))$ on its spaces of line sub-bund\-les~$\mathrm{L}_{j}(E)$. In the particular case when $m_{1}=m_{2}=m$, this reduces to the standard action on the projective line $\mathrm{L}_{m}(E)$ by M\"obius transformations, and moreover, $\mathrm{L}_{m-1}(E)$ is a principal homogeneous space for $\PP(\Aut(E))$. Otherwise, in the case when $m_{1}>m_{2}$, $\mathrm{L}_{m_{2}}(E)$ is a principal homogeneous space for $\mathrm{R}(\Aut(E))$.
Since the action of $\Aut (E)$ on the spaces $\mathrm{L}_{j}(E)$ commutes with the operations of evaluation, every $\mathrm{B}_{I,j}(E)$ is invariant under the induced $\PP(\Aut(E))$-action. In the particular case when $I=\{1,2,3,4\}$, it follows from the general properties of polynomial interpolation under an arbitrary choice of bundle isomorphism $E\cong \calO(m_{1})\oplus\calO(m_{2})$ that if $\dim \mathrm{L}_{j}(E)<4$, then the induced line bundle evaluation map
\[
\mathrm{ev}_{\mathrm{L}_{j}(E)}\colon\quad \mathrm{L}_{j}(E)\rightarrow \mathrm{QP}(E),\qquad L\subset E\mapsto \bigl(L\vert_{z_{1}},L\vert_{z_{2}},L\vert_{z_{3}},L\vert_{z_{4}}\bigr)
\]
is a holomorphic embedding, and
\[
\mathrm{ev}_{\mathrm{L}_{j}(E)}(\mathrm{L}_{j}(E)) = \mathrm{B}_{\{1,2,3,4\},j}(E).
\]
In addition, when $m_{1} > m_{2}$, $\calO(m_{1})$-interpolation relative to a subset $I$ endows $\mathrm{QP}(E)$ with the following $\PP(\Aut(E))$-invariant stratification
\[
\mathrm{QP}(E) = \bigsqcup_{I\subset\{1,2,3,4\}}\mathrm{B}_{I,m_{1}}(E),\qquad
\mathrm{B}_{I,m_{1}}(E) := \bigg\{\prod_{i\in I}\{E_{1}\vert_{z_{i}}\} \bigg\} \times \bigg\{\prod_{j\in I^{\mathsf{c}}}\cV_{z_{i}}(E) \bigg\}.
\]

\subsection[Combinatorial stratification of orbit spaces in QP(E)]
{Combinatorial stratification of orbit spaces in $\boldsymbol{\mathrm{QP}(E)}$}
We will now describe the action of $\PP(\Aut(E))$ on the spaces $\mathrm{QP}(E)$. Due to the structural differences of the groups $\PP(\Aut(E))$ depending on whether $E$ is semi-simple or not, the treatment is split into the two cases $m_{1} = m_{2}$ and $m_{1} > m_{2}$. In the first case, we define the \emph{configuration space locus} to be the following Zariski open subset of $\mathrm{QP}(E)$
\begin{gather*}
\mathrm{C}(E):=\bigcap_{\vert I\vert =2} \mathrm{B}^{\mathsf{c}}_{I,m}(E)\cong \mathrm{Conf}_{4}\big(\CC\PP^{1}\big).
\end{gather*}

\begin{Proposition}\label{prop:U}
If $m_{1}-m_{2} \leq 2$, then $\PP(\Aut(E))$ acts freely on the Zariski open subset
\[
\mathrm{U}(E) :=
\displaystyle\mathrm{QP}(E)\Big\backslash \bigcup_{I\sqcup I' = \{1,2,3,4\}} \bigl\{\mathrm{B}_{I,m_{1}}(E)\cap\mathrm{B}_{I',m_{2}}(E)\bigr\},
\]
which is a $\PP(\Aut(E))$-principal homogeneous space if $m_{1}-m_{2}=2$, and empty if $m_{1}-m_{2}>2$. There is a $\PP(\Aut(E))$-invariant ``cross-ratio" map $\operatorname{cr}\colon \mathrm{U}(E)\rightarrow \CC\PP^{1}$ identifying the $\mathrm{U}(E)$-orbit space with a projective line with three $($resp.~four$)$ double points if $m_{1} = m_{2}$ $($resp.~$m_{1} - m_{2} = 1)$, and such that $\operatorname{cr}(\mathrm{U}'(E)) = \CC\PP^{1}\backslash\{z_{1},0,1,\infty\}$ for
\[
\mathrm{U}'(E):= \begin{cases}
\mathrm{C}(E)\big\backslash\mathrm{B}_{\{1,2,3,4\},m-1}(E), & m_{1}=m_{2}=m,
\\[1ex]
\displaystyle\mathrm{B}_{\varnothing ,m_{1}}(E)\Big\backslash\bigsqcup_{\vert I'\vert \geq 3} \mathrm{B}_{I',m_{2}}(E), & m_{1}-m_{2}=1.
\end{cases}
\]
\end{Proposition}

\begin{proof}
\noindent \textit{Case} $m_{1} = m_{2} = m$. All nonempty loci $\mathrm{B}_{I,m}(E)$ provide an invariant characterization of the trivialization-dependent identities
\[
F_{i}=F_{j},\qquad \forall\; i,j\in I,
\]
since for any $\vert I \vert \geq 2$
\[
\mathrm{B}_{I,m}(E) = \bigcap_{I'\subset I,\; \vert I'\vert =2} \mathrm{B}_{I',m}(E).
\]
Therefore, $\mathrm{C}(E)$ is the locus in $\mathrm{U}(E)$, where $\PP(\Aut(E))$ acts properly, in such a way that $\mathrm{C}(E)/\PP(\Aut(E))\cong \CC\PP^{1}\backslash\{0,1,\infty\}$. Moreover, the six connected components in
\[
\mathrm{U}(E)\backslash\mathrm{C}(E) = \bigsqcup_{\vert I\vert =2} \bigl\{\mathrm{B}_{I,m}(E)\cap\mathrm{U}(E)\bigr\}
\]
are principal homogeneous spaces for $\PP(\Aut(E))$, and it follows that there is a unique $\PP(\Aut(E))$-invariant map $\operatorname{cr}\colon \mathrm{U}(E)\rightarrow \CC\PP^{1}$ normalized in the following manner
\begin{gather*}
\operatorname{cr}\bigl(\mathrm{B}_{\{1,2\},m}(E)\cap\mathrm{U}(E)\bigr) = \operatorname{cr}\bigl(\mathrm{B}_{\{3,4\},m}(E)\cap\mathrm{U}(E)\bigr) = 0,
\\
\operatorname{cr}\bigl(\mathrm{B}_{\{1,3\},m}(E)\cap\mathrm{U}(E)\bigr) = \operatorname{cr}\bigl(\mathrm{B}_{\{2,4\},m}(E)\cap\mathrm{U}(E)\bigr) = 1,
\\
\operatorname{cr}\bigl(\mathrm{B}_{\{1,4\},m}(E)\cap\mathrm{U}(E)\bigr) = \operatorname{cr}\bigl(\mathrm{B}_{\{2,3\},m}(E)\cap\mathrm{U}(E)\bigr) = \infty.
\end{gather*}
In turn, since $ \mathrm{B}_{\{1,2,3,4\},m-1}(E)$ is also a principal homogeneous space for $\PP(\Aut(E))$, we have that $\operatorname{cr}(\mathrm{U}'(E)) = \CC\PP^{1}\backslash\{z_{1},0,1,\infty\}$, since under the previous
normalization
\[
\operatorname{cr}\bigl(\mathrm{B}_{\{1,2,3,4\},m-1}(E)\bigr) = (z_{1},z_{2};z_{3},z_{4}) =z_{1}.
\]

\medskip\noindent\textit{Case $m_{1}>m_{2}$.}
In general, we have that for any partition $I\sqcup I' = \{1,2,3,4\}$ and $\vert I\vert > 2 - m_{1} + m_{2}$,
\[
\mathrm{B}_{I,m_{1}}(E)\cap\mathrm{B}_{I',m_{2}}(E) = \mathrm{B}_{I,m_{1}}(E).
\]
Moreover, when $\vert I\vert \leq 2 - m_{1} + m_{2}$, the $\PP(\Aut(E))$-action on $\mathrm{B}_{I,m_{1}}(E)\backslash \mathrm{B}_{I^{\mathsf{c}},m_{2}}(E)$ is free and proper, as it can be factored through the following intermediate affine quotients
\[
\bigl\{\mathrm{B}_{I,m_{1}}(E)\backslash\mathrm{B}_{I^{\mathsf{c}},m_{2}}(E)\bigr\}/\mathrm{R}(\Aut(E))\cong \CC^{3-m_{1}+m_{2}-\vert I\vert }\backslash\{0\},
\]
on which $\PP(\Aut (E))/\mathrm{R}(\Aut(E))\cong\CC^{*}$ acts in the standard way. Since we can re-express
\[
\mathrm{U}(E) = \bigsqcup_{\vert I\vert \leq 2-m_{1}+m_{2}}\bigl\{\mathrm{B}_{I,m_{1}}(E)\backslash\mathrm{B}_{I^{\mathsf{c}},m_{2}}(E)\bigr\},
\]
it follows that $\mathrm{U}(E)$ is empty if $m_{1}-m_{2}>2$, a principal homogeneous space for $\PP(\Aut(E))$ if $m_{1}-m_{2}=2$, and moreover, that there exists a $\PP(\Aut(E))$-invariant map $\mathrm{U}(E)\rightarrow \CC\PP^{1}$ when $m_{1}-m_{2}=1$, which determines an isomorphism
\[
\mathrm{U}'(E)/\PP(\Aut(E))\cong\CC\PP^{1}\backslash\{z_{1},0,1,\infty\},
\]
and the complement $\mathrm{U}(E)\backslash\mathrm{U}'(E)$ consists of eight connected components, namely
\[
\mathrm{B}_{\varnothing ,m_{1}}(E)\cap\mathrm{B}_{\{j,k,l\},m_{2}}(E)\qquad \text{and}\qquad \mathrm{B}_{\{i\},m_{1}}(E)\backslash\mathrm{B}_{\{j,k,l\},m_{2}}(E),
\]
each of which is a principal homogeneous space for $\PP(\Aut(E))$. These eight orbits determine four different pairs of compactification points for $\mathrm{U}'(E)/\PP(\Aut(E))$.
\end{proof}

\begin{Proposition}\label{prop:X-Y}\quad
\begin{enumerate}\itemsep=0pt
\item[$(i)$] The elements of the sets
\begin{gather*}
\mathrm{X}(E) :=
\displaystyle\bigcup_{\vert I\vert = 2 - m_{1}+m_{2}}\bigl\{\mathrm{B}_{I,m_{1}}(E)\cap\mathrm{B}_{I^{\mathsf{c}},m_{2}}(E)\bigr\} \big\backslash\mathrm{B}_{\{1,2,3,4\},m_{1}}(E),
\\
\mathrm{Y}(E) :=
\displaystyle\bigcup_{\substack{\vert I\vert = 3 - m_{1}+m_{2},\\\vert I\vert = 1 - m_{1}+m_{2}}}\bigl\{\mathrm{B}_{I,m_{1}}(E)\cap\mathrm{B}_{I^{\mathsf{c}},m_{2}}(E)\bigr\} \big\backslash\mathrm{B}_{\{1,2,3,4\},m_{1}}(E)
\end{gather*}
have pointwise stabilizer subgroups $\PP(\Aut(E))_{(F_{1},F_{2},F_{3},F_{4})}\cong \CC^{*}$.

\item[$(ii)$] If $m_{1}-m_{2} \leq 2$ then $\mathrm{X}(E)$ and $\mathrm{Y}(E)$ are degeneration loci for the orbits in $\mathrm{U}(E)\backslash \mathrm{U}'(E)$, where $\mathrm{U}'(E)=\varnothing $ if $m_{1}-m_{2}=2$, i.e.,
\[
\mathrm{X}(E)\cup \mathrm{Y}(E)\subset \overline{\mathrm{U}(E)\backslash\mathrm{U}'(E)}.
\]
\end{enumerate}
\end{Proposition}

\begin{proof}
The first claim is verified in analogy to the proof of Proposition~\ref{prop:U}. Notice that when $m_{1}>m_{2}$, $\{\mathrm{B}_{I,m_{1}}(E)\cap\mathrm{B}_{I^{\mathsf{c}},m_{2}}(E)\} \cap\mathrm{B}_{\{1,2,3,4\},m_{1}}(E)=\varnothing $, $\mathrm{X}(E)=\varnothing $ if $m_{1}-m_{2}> 2$, and the stabilizer of any point in $\mathrm{X}(E)$ or $\mathrm{Y}(E)$ is biholomorphic to $\PP(\Aut(E))/\mathrm{R}(\Aut(E))\cong\CC^{*}$. The~second claim is trivial when $m_{1}-m_{2}=2$. When $m_{1}=m_{2}=m$, the result follows since for any partition $I\sqcup I' =\{1,2,3,4\}$, $\vert I\vert =2$
\[
\bigl\{\mathrm{B}_{I,m}(E)\cap\mathrm{B}_{I',m}(E)\bigr\} \subset \overline{\mathrm{B}_{I,m}(E)\cap\mathrm{U}(E)}\cap\overline{\mathrm{B}_{I',m}(E)\cap\mathrm{U}(E)},
\]
while when $\vert I\vert = 3$
\[
\mathrm{B}_{I,m}(E) \subset \bigcap_{I'\subset I,\, \vert I'\vert = 2}\overline{\mathrm{B}_{I',m}(E)\cap\mathrm{U}(E)}.
\]
The case $m_{1}-m_{2} =1$, $\vert I\vert =1$, follows since for each partition $\{i\}\sqcup\{j,k,l\} = \{1,2,3,4\}$
\[
\mathrm{B}_{\{i\},m_{1}}(E)\cap\mathrm{B}_{\{j,k,l\},m_{2}}(E) \subset \overline{\mathrm{B}_{\varnothing ,m_{1}}(E)\cap\mathrm{B}_{\{j,k,l\},m_{2}}(E)} \cap \overline{\mathrm{B}_{\{i\},m_{1}}(E)\backslash\mathrm{B}_{\{j,k,l\},m_{2}}(E)},
\]
and is trivial when $\vert I\vert =2$ or $\vert I \vert =0$ from the definition of $\mathrm{U}'(E)$.
\end{proof}

The stratification of $\mathrm{QP}(E)$ can be completed by also grouping quasi-parabolic structures with stabilizer subgroups of higher dimension. We don't describe these additional strata explicitly since they won't appear in subsequent constructions.

\section[Structural results and stratification of QPH(E)]
{Structural results and stratification of $\boldsymbol{\mathrm{QPH}(E)}$}\label{sec:struct}

The construction of geometric models for Harder--Narasimhan strata will be based on the following series of results on the structure of spaces of quasi-parabolic Higgs bundles. For any fixed holomorphic rank 2 bundle $E\rightarrow \CC\PP^{1}$, let $\mathrm{QPH}(E)$ denote the space of all quasi-parabolic Higgs bundles $(E_{\qp},\Phi)$ with quasi-parabolic structure supported on $D$. We can identify $\mathrm{QPH}(E)$ with an algebraic subvariety of the product
\[
\mathrm{QP}(E)\times H^{0}\bigl(\CC\PP^{1},\End(E)\otimes K_{\CC\PP^{1}}(D)\bigr)
\]
characterized as the locus of all points $(F_{1},F_{2},F_{3},F_{4},\Phi)$ satisfying the incidence constraints
\begin{gather*}
\Res_{z_{i}}\Phi\in\mathfrak{n}(E\vert_{z_{i}}),\qquad
F_{i}\subset\ker (\Res_{z_{i}}\Phi),\qquad
i=1,2,3,4.
\end{gather*}
Equivalently, if we denote by $\mathrm{Int}(E)$ the intersection of the four quadrics in the space
\[
\bigl\{ \Phi \in H^{0}\big(\CC\PP^{1},\End(E)\otimes K_{\CC\PP^{1}}(D)\big)\colon \tr(\Phi)\vert_{z_{i}}=0,\; i=1,2,3,4\bigr\}
\]
defined by the equations
\[
\det(\Phi)\vert_{z_{i}} =0,\qquad i=1,2,3,4,
\]
then $\mathrm{QPH}(E)$ is modeled by the blow-up of $\mathrm{Int}(E)$ along each of the codimension 2 loci
\begin{gather}\label{eq:blow-up-loci}
\Res_{z_{i}}\Phi =0,\qquad i=1,2,3,4.
\end{gather}
There is a natural action of $\Aut(E)$ on $\mathrm{QPH}(E)$, whose point stabilizers contain $\mathrm{Z}(\Aut(E))$, descending to an action of the quotient group $\PP(\Aut(E))$. The proof of the following lemma is obvious.

\begin{Lemma}
There is a bijective correspondence between isomorphism classes of quasi-parabolic Higgs bundles $\{(E_{\qp},\Phi)\}$ and $\PP(\Aut(E))$-orbits in $\mathrm{QPH}(E)$.
\end{Lemma}

Since in the present case $\dim H^{0}(\CC\PP^{1},K^{2}_{\CC\PP^{1}}(D)) =1$, we will define the \textit{bulk} as
\[
\mathrm{QPH}_{\CC^{*}}(E) := \det{}^{-1}\bigl(H^{0}\big(\CC\PP^{1},K^{2}_{\CC\PP^{1}}(D)\big)\backslash\{0\}\bigr),
\]
and the \emph{nilpotent locus} as
\[
\mathrm{QPH}_{0}(E):=\det{}^{-1}(0).
\]
For the embedding $\iota \colon \mathrm{QP}(E)\hookrightarrow \mathrm{QPH}_{0}(E)$ defined as $\iota:= \{0\}\times\textrm{id}$, we will denote
\[
\mathrm{Q}(E) := \iota(\mathrm{QP}(E)),
\]
which is also characterized as the fixed-point locus of the holomorphic involution defined as
\begin{gather}\label{eq:def-tau}
\tau(F_{1},F_{2},F_{3},F_{4},\Phi):= (F_{1},F_{2},F_{3},F_{4},-\Phi).
\end{gather}
More generally, $\mathrm{QPH}(E)$ is equipped with a $\CC^{*}$-action,
defined for any $c\in\CC^{*}$ as
\begin{gather}\label{eq:def-C*}
c\cdot (F_{1},F_{2},F_{3},F_{4},\Phi) = (F_{1},F_{2},F_{3},F_{4},c\Phi),\qquad c\in\CC^{*}.
\end{gather}
Since both $\tau$ and this $\CC^{*}$-action are $\PP(\Aut(E))$-equivariant, they descend to an involution and a $\CC^{*}$-action on the space of $\PP(\Aut(E))$-orbits of $\mathrm{QPH}(E)$. While they are only trivial on $\mathrm{Q}(E)$, we will verify that their descents
possess additional fixed loci.
The next results characterize the structure of $\mathrm{QPH}(E)$ in terms of the stratification
\[
\mathrm{QPH}(E) = \mathrm{QPH}_{0}(E)\sqcup \mathrm{QPH}_{\CC^{*}}(E).
\]

\begin{Lemma}\label{lemma:nilpotent}
A nonzero parabolic Higgs field $\Phi$ preserves a line sub-bundle $L(\Phi)\subset E$ if and only if it is nilpotent, in which case $L(\Phi)$ is unique and determined by $\ker(\Phi)$.
\end{Lemma}

\begin{proof}
Given $\Phi\neq 0$, the equation $\Phi\cdot L = s L$ with $s\in H^{0}(\CC\PP^{1},K_{\CC\PP^{1}}(D))$ requires that $s^{2} = -\det \Phi\in H^{0}(\CC\PP^{1},K^{2}_{\CC\PP^{1}}(D))$, from which $s\equiv 0$, i.e., $\Phi$ is nilpotent. Conversely, if $\Phi$ is nilpotent, we can reconstruct $L(\Phi)$ in a unique way. Since the zero set of $\Phi$ as a holomorphic section of $\End(E)\otimes K_{\CC\PP^{1}}(D)$ is finite, we can associate an effective divisor $(\Phi)$ to $\Phi$. The complement of such zero set is the Zariski open set $\cU$ of points $z\in\CC\PP^{1}$, where $\ker(\Phi)(z)\subset E\vert_{z}$ (or $\ker(\Res_{z_{i}}\Phi)$ if $z=z_{i}$, $i=1,2,3,4$) is 1-dimensional. Hence $\ker(\Phi)\vert_{\cU}$ is a line sub-bundle of $E\vert_{\cU}$. Let $\sigma$ be a holomorphic section of the line bundle $[(\Phi)]$ such that $(\sigma)=(\Phi)$. Then by construction $\sigma^{-1}\Phi$ is a holomorphic section of $\End E\otimes K_{\CC\PP^{1}}(D-(\Phi))$ such that $\ker(\sigma^{-1}\Phi)(z)$ is 1-dimensional $\forall\; z\in\CC\PP^{1}$. Therefore $L(\Phi):=\ker(\sigma^{-1}\Phi)$ is a line sub-bundle of $E$ such that $L\vert_{\cU} =\ker(\Phi)\vert_{\cU}$ and invariant under $\Phi$ over all $\CC\PP^{1}$.
\end{proof}

Two important conclusions follow from Lemma~\ref{lemma:nilpotent}: when it exists, an element of $\mathrm{QPH}_{\CC^{*}}(E)$ can never be unstable, and moreover, the elements in $\mathrm{QPH}_{0}(E)\backslash \mathrm{Q}(E)$ can be classified according to their unique invariant line sub-bundles, inducing a finer stratification of $\mathrm{QPH}_{0}(E)$ that we will describe in detail.

\begin{Lemma}\label{lemma:cases}
$\mathrm{QPH}_{\CC^{*}}(E)$ is nonempty if and only if $m_{1}-m_{2}\leq 2$. Given a line sub-bundle $L\subset E$, there is a nonzero nilpotent parabolic Higgs field $\Phi$ such that $L(\Phi)=L$ if and only if
\begin{gather}\label{eq:bound-deg(L)}
2(\deg(L)+1) \geq \deg (E).
\end{gather}
\end{Lemma}

\begin{proof}
Any given choice of bundle isomorphism $E\cong \calO(m_{1})\oplus\calO(m_{2})$ induces an isomorphism
\begin{gather*}
\End(E)\otimes K_{\CC\PP^{1}}(D)\cong
\begin{pmatrix}
\calO(2) & \calO(m_{1}-m_{2}+2)\\
\calO(m_{2}-m_{1}+2) & \calO(2)
\end{pmatrix}\!,
\end{gather*}
from which it follows that a parabolic Higgs field on $E$ gets characterized as a holomorphic section
\begin{gather}\label{eq:Phi-1}
\Phi =
\begin{pmatrix}
u_{2} & -v_{m_{1}-m_{2}+2}\\
w_{m_{2}-m_{1}+2} & -u_{2}
\end{pmatrix}
\end{gather}
whose pointwise evaluation on $D$ moreover satisfies the set of nilpotency constraints
\[
\det(\Phi)\vert_{z_{i}} = \bigl({-}u_{2}^{2} + v_{m_{1}-m_{2}+2}w_{m_{2}-m_{1}+2}\bigr)\big\vert_{z_{i}}=0,\qquad i=1,2,3,4.
\]
It follows that the map $\det$ is surjective if and only if $m_{1}-m_{2}\leq 2$, from which the first claim follows. To prove the second claim, first assume that $m_{1} > m_{2}$ and $L=E_{1}$. A nonzero nilpotent parabolic Higgs field $\Phi$ preserving $E_{1}$ is characterized in terms of the identity $w_{m_{2}-m_{1}+2}\equiv 0$ (which is independent of the choice of isomorphism $E \cong \calO(m_{1})\oplus\calO(m_{2})$) since the nilpotency constraints imply that $u_{2}\equiv 0$ necessarily. Since $\Phi$ is identified with a nonzero holomorphic section of $\calO(m_{1}-m_{2}+2)$, it follows that there is always some $\Phi$ such that $L(\Phi)=E_{1}$.

Otherwise, assume that $L\neq E_{1}$ if $m_{1}>m_{2}$, or equivalently, that $\deg(L)\leq m_{2}$. Given a nonzero nilpotent parabolic Higgs field $\Phi$ such that $L(\Phi)=L$, we can choose $\sigma_{k}\neq 0\in H^{0}(\CC\PP^{1},[(\Phi)])$ for some $0\leq k\leq 2-m_{1}+m_{2}$, such that $(\sigma_{k})=(\Phi)$ and express
\begin{gather}\label{eq:Phi-2}
\Phi =
\sigma_{k}\begin{pmatrix}
u_{2-k} & -v_{2-m_{2}+m_{1}-k}\\
w_{2-m_{1}+m_{2}-k} & -u_{2-k}
\end{pmatrix}\!,
\end{gather}
such that $w_{2-m_{1}+m_{2}-k}\neq 0$ if $m_{1} > m_{2}$, and the following equation holds in $H^{0}(\CC\PP^{1},\calO(4-2k))$
\begin{gather}\label{eq:det=0}
u_{2-k}^{2}=v_{2-m_{2}+m_{1}-k}w_{2-m_{1}+m_{2}-k}.
\end{gather}
The holomorphic sections $u_{2-k}$ and $w_{2-m_{1}+m_{2}-k}$ could have a (simple) common zero only in the exceptional case when $m_{1}=m_{2}=m$ and $k=0$, when a parabolic Higgs field takes the form
\begin{gather}\label{eq:Phi(m-1)}
\Phi = \begin{pmatrix}
u'_{1}v'_{1} & -u'^{2}_{1}\\
v'^{2}_{1} & -u'_{1}v'_{1}
\end{pmatrix}\!.
\end{gather}
Since
$L(\Phi)$ is generated by $(u'_{1},v'_{1})$, it follows that $L(\Phi)\cong \calO(m-1)$.
Otherwise, we can assume that $u_{2-k}$ and $w_{2-m_{1}+m_{2}-k}$ don't have a common zero. Then $L(\Phi)$ is generated by
\[
(u_{2-k},w_{2-m_{1}+m_{2}-k})\in H^{0}\bigl(\CC\PP^{1},\calO(2-k)\oplus\calO(2-m_{1}+m_{2}-k)\bigr)
\]
and $\deg(L(\Phi)) = m_{1}-2 + k\geq m_{1}-2$. If $m_{1}-m_{2}=2$, then $k=0$ and $\deg(L)=m_{2}$. If~$m_{1}-m_{2}=0$ (resp.~$1$) and $k=0,1$ (resp.~$0$), then $w_{2-m_{1}+m_{2}-k}$ would have at least one zero, and it would follow from \eqref{eq:det=0} that $u_{2-k}$ and $w_{2-m_{1}+m_{2}-k}$ have a common zero, a contradiction. Therefore $k=2$ (resp.~1), or equivalently, $\deg(L(\Phi)) = m_{2}$. In all of the four possible cases, we conclude that
\[
k= 2(\deg(L(\Phi))+1) -\deg(E).
\]
Therefore, $k \geq 0$ is equivalent to the lower bound \eqref{eq:bound-deg(L)}, and the claim follows.

Conversely, consider first any line sub-bundle $L\subset E$ such that $\deg(L)=m_{2}$, generated by a holomorphic section $(u_{m_{1}-m_{2}},w_{0})$ with $w_{0}\neq 0$, and assume that \eqref{eq:bound-deg(L)} is satisfied, i.e., $m_{1}-m_{2}\leq 2$. It readily follows that there exists a holomorphic section
\[
v_{2m_{1}-2m_{2}} \in H^{0}\bigl(\CC\PP^{1},\calO(2m_{1}-2m_{2})\bigr),\qquad
(v_{2m_{1}-2m_{2}}) = 2(u_{m_{1}-m_{2}}),
\]
solving \eqref{eq:det=0}, for $k =2-m_{1}+m_{2}$. Any $\Phi$ constructed according to \eqref{eq:Phi-2}, for an arbitrary choice of $\sigma_{2-m_{1}+m_{2}}\in H^{0}(\CC\PP^{1},\calO(2-m_{1}+m_{2}))\backslash\{0\}$, would be such that $L(\Phi)=L$. The remaining case $m_{1}=m_{2}=m$, $\deg(L)=m-1$ follows analogously.
\end{proof}

The next result, on the stratifications of $\mathrm{QPH}_{0}(E)$ for any bundle splitting type, is immediate from Lemma~\ref{lemma:cases}. Notice that the $\CC^{*}$-action on $\mathrm{QPH}(E)$ preserves its stratification, and that of $\mathrm{QPH}_{0}(E)$ in particular.

\begin{Corollary}\label{cor:strata}
Let $\mathrm{R}(E)$ $($resp.~$\mathrm{S}_{j}(E))$ denote the locus of quasi-parabolic Higgs bundles such that $\Phi\neq 0$, for which $L(\Phi)=E_{1}$ $($resp.~$L(\Phi)\in \mathrm{L}_{j}(E))$. $\mathrm{QPH}_{0}(E)$ is stratified as
\[
\mathrm{QPH}_{0}(E) =
\begin{cases}
\mathrm{Q}(E) \sqcup \mathrm{S}_{m}(E)\sqcup \mathrm{S}_{m-1}(E) & \text{if} \quad E\cong\calO(m)\oplus\calO(m),
\\
\mathrm{Q}(E) \sqcup \mathrm{R}(E) \sqcup \mathrm{S}_{m}(E) & \text{if} \quad E\cong\calO(m+1)\oplus\calO(m),
\\
\mathrm{Q}(E) \sqcup \mathrm{R}(E) \sqcup \mathrm{S}_{m-1}(E) & \text{if} \quad E\cong\calO(m+1)\oplus\calO(m-1),
\\
\mathrm{Q}(E)\sqcup \mathrm{R}(E) & \text{if} \quad m_{1}-m_{2}\geq 3.
\end{cases}
\]
\end{Corollary}

Let $\mathrm{P}_{j}(E)\rightarrow \mathrm{L}_{j}(E)$ be the tautological principal $\CC^{*}$-bundle associated to the canonical projective embeddings of $\mathrm{L}_{j}(E)$. For any $I\subset\{1,2,3,4\}$, $\mathrm{Bl}_{I}(\CC\PP^{1})$ will denote the rational nodal curve resulting from blowing-up $\CC\PP^{1}\subset \CC\PP^{2}$ at $z_{i}$ $\forall i\in I$. In particular, $\mathrm{Bl}_{\{1,2,3,4\}}(\CC\PP^{1})$ is the $D_{4}$-configuration. Consider the $\PP(\Aut(E))$-equivariant maps
\[
L_{E,j}\colon\quad \mathrm{S}_{j}(E)\rightarrow \mathrm{L}_{j}(E),\qquad (E_{\qp},\Phi)\mapsto L(\Phi),
\]
and in the case $k= 2(j+1) -\deg(E) > 0$, also the $\PP(\Aut (E))$-invariant maps
\[
\Div_{E,j}\colon\quad \mathrm{S}_{j}(E)\rightarrow S^{k}(\CC\PP^{1}),\qquad (E_{\qp},\Phi)\mapsto (\Phi).
\]

\begin{Proposition}\label{prop:S_{j}(E)}
The restriction of the $\CC^{*}$-action to each $\mathrm{S}_{j}(E)$ turns it into a principal $\CC^{*}$-bundle, whose base is modeled by the $3$-dimensional subvariety of the product
\[
\mathrm{L}_{j}(E) \times S^{2(j+1) -\deg(E)}\big(\CC\PP^{1}\big)\times \mathrm{QP}(E),
\]
defined by the incidence constraints $F_{i} = L(\Phi)\vert_{z_{i}}$ whenever $z_{i}\not\in(\Phi)$, and projection given by the map $L_{E,j} \times \Div_{E,j} \times \mathrm{par}$.
 In particular, there is an isomorphism
\[
\mathrm{S}_{j}(E) \cong
\begin{cases}
\mathrm{P}_{j}(E), & 2(j+1) -\deg(E) = 0,\\
\mathrm{Bl}_{\{1,2,3,4\}}(\CC\PP^{1})\times \mathrm{P}_{j}(E), & 2(j+1) -\deg(E) = 1.
\end{cases}
\]
In general,
\[
\mathrm{par}(\mathrm{S}_{j}(E)) = \bigcup_{\vert I\vert = 2(j+1) -\deg(E)}\mathrm{B}_{I,j}(E).
\]
\end{Proposition}

\begin{proof}
The map $L_{E,j}\times\Div_{E,j}\times\mathrm{par}$ is invariant under the $\CC^{*}$-action on $\mathrm{QPH}(E)$ by definition. The first claim follows from the definition of $\mathrm{S}_{j}(E)$ as a stratum of $\mathrm{QPH}(E)$, and the expression of the parabolic Higgs field $\Phi$ of any element $(E_{\qp},\Phi)\in\mathrm{S}_{j}(E)$ in one of the two canonical forms~\eqref{eq:Phi-2} and~\eqref{eq:Phi(m-1)} under a choice of isomorphism $E\cong \calO(m_{1})\oplus\calO(m_{2})$ as in proof of Lemma~\ref{lemma:cases}, with $k=2(j+1) -\deg(E)$, and in such a way that $(\Phi) = (\sigma_{k})$.
\end{proof}

\begin{Proposition}\label{prop:QPH-par}\quad
\begin{enumerate}\itemsep=0pt
\item[$(i)$] If $E$ is evenly-split, then
\[
\mathrm{par}\bigl(\mathrm{QPH}_{\CC^{*}}(E)\bigr) = \mathrm{U}'(E)\sqcup \mathrm{X}(E),
\]
and
\[
\mathrm{par}\bigl(\mathrm{QPH}_{0}(E)\backslash\mathrm{Q}(E)\bigr) = \mathrm{X}(E)\sqcup\mathrm{Y}(E).
\]
\item[$(ii)$] If $m_{1}-m_{2}=2$, then
\[
\mathrm{par}\bigl(\mathrm{QPH}_{\CC^{*}}(E)\bigr) = \mathrm{par}(\mathrm{S}_{m-1}(E)) = \mathrm{X}(E).
\]
\item[$(iii)$] If $m_{1}-m_{2}\geq 2$, then $\mathrm{par}\vert_{\mathrm{R}(E)}$ is surjective.
\end{enumerate}
\end{Proposition}

\begin{proof}
{\sloppy $(i)$ Assume that $E$ is evenly-split. We will treat both cases independently using the parametrization \eqref{eq:Phi-1} of parabolic Higgs fields on a~quasi-parabolic bundle $E_{\qp}$ depen\-ding on a choice of isomorphism $E\cong\calO(m_{1})\oplus\calO(m_{2})$. When $m_{1}=m_{2} = m$, for a partition $\{1,2,3,4\} = \{i\}\sqcup\{j,k,l\}$, consider the triple of Lagrange interpolating sections $\{s_{j},s_{k},s_{l}\}$ spanning $H^{0}(\CC\PP^{1},\calO(2))$ relative to the triple $\{z_{j},z_{k},z_{l}\}$. Since the ``residue evaluation" map
\[
\mathfrak{n}(F_{j})\oplus\mathfrak{n}(F_{k})\oplus\mathfrak{n}(F_{l}) \rightarrow\mathfrak{sl}(E\vert_{z_{i}}),\qquad
(\phi_{j},\phi_{k},\phi_{l})\mapsto (\phi_{j}s_{j}+\phi_{k}s_{k}+\phi_{l}s_{l})\vert_{z_{i}}
\]}\noindent
is an isomorphism, there is a line in $\mathfrak{n}(F_{j})\oplus\mathfrak{n}(F_{k})\oplus\mathfrak{n}(F_{l})$ mapping to $\mathfrak{n}(F_{i})$. This implies that any quasi-parabolic Higgs field $\Phi$ is expressed as $\phi_{j}s_{j}+\phi_{k}s_{k}+\phi_{l}s_{l}$, and consequently $\mathrm{par}(\mathrm{QPH}(E)\backslash\mathrm{Q}(E)) = \mathrm{QP}(E)$.

If $\mathrm{par}(E_{\qp},\Phi)\in\mathrm{B}_{I,m}(E)$ for some $\vert I\vert \geq 3$ and $\Phi \neq 0$, then after representing $\Phi = \phi_{j}s_{j}+\phi_{k}s_{k}+\phi_{l}s_{l}$ as before for some $\{j,k,l\}\subset I$, we see that $\dim (\mathrm{Span}\{\phi_{j},\phi_{k},\phi_{l}\}) = 1$, and consequently $(E_{\qp},\Phi)\in\mathrm{S}_{m}(E)$. An analogous argument for $\mathrm{par}(E_{\qp},\Phi)\in\mathrm{B}_{\{1,2,3,4\},m-1}(E)$ implies that $(E_{\qp},\Phi)\in\mathrm{S}_{m-1}(E)$. Finally, if $\mathrm{par}(E_{\qp},\Phi)\in\mathrm{B}_{I,m}(E)$ for some $(E_{\qp},\Phi)\in\mathrm{QPH}_{\CC^{*}}(E)$ and $\vert I\vert = 2$, then the same argument implies that $\mathrm{par}(E_{\qp},\Phi)\in\{\mathrm{B}_{I,m}(E)\cap \mathrm{B}_{I^{\mathsf{c}},m}(E)\}\backslash\mathrm{B}_{\{1,2,3,4\},m}(E)\subset \mathrm{X}(E)$. This last possibility occurs when the orbit of $(E_{\qp},\Phi)$ contains an element with parabolic Higgs field
\begin{gather}\label{eq:HX-0}
\Phi = \begin{pmatrix}
0 & -v_{2}\\
w_{2} & 0
\end{pmatrix}\!,
\end{gather}
with $(v_{2}) = z_{i}+z_{j}$ and $(w_{2}) = z_{k}+z_{l}$ for a partition $\{i,j\}\sqcup\{k,l\} = \{1,2,3,4\}$. The $\PP(\Aut(E))$-stabilizers of these points are trivial, and their projection exhaust $\mathrm{X}(E)$. The surjectivity of $\mathrm{par}\vert_{\mathrm{QPH}(E)\backslash\mathrm{Q}(E)}$ and the definition of $\mathrm{U}'(E)$ (Proposition~\ref{prop:U}) imply that $\mathrm{QPH}_{\CC^{*}}(E)$ necessa\-rily projects onto $\mathrm{U}'(E)\sqcup \mathrm{X}(E)$. The characterization of $\mathrm{par}(\mathrm{QPH}_{0}(E)\backslash\mathrm{Q}(E))$ follows from Corollary~\ref{cor:strata} and Proposition~\ref{prop:S_{j}(E)}.

When $E\cong \calO(m+1)\oplus\calO(m)$, we recall that the divisor $(w_{1})$ associated to the expression~\eqref{eq:Phi-1} is independent of the choice of Birkhoff--Grothendieck splitting and an invariant of any point $\mathrm{QPH}(E)\backslash\mathrm{Q}(E)$. The projection of the Zariski open subset of $\mathrm{QPH}_{\CC^{*}}(E)$ for which $(w_{1})\in\CC\PP^{1}\backslash\{z_{1},0,1,\infty\}$ equals $\mathrm{U}'(E)$ (defined in Proposition~\ref{prop:U}), since for any $(F_{1},F_{2},F_{3},F_{4})\in\mathrm{U}'(E)$ and any $u_{2}$ such that $z_{i}\not\in(u_{2})$ there exists a unique $v_{3}$ such that $\ker(\Res_{z_{i}}\Phi) = F_{i}$, $i=1,2,3,4$, and on the other hand, an element $\mathrm{QPH}(E)\backslash\mathrm{Q}(E)$ projecting to $\mathrm{B}_{I,m}(E)$ for any $\vert I\vert \geq 3$ necessa\-rily belongs to $\mathrm{S}_{m}(E)$. In turn, the projection of the complementary loci of points in $\mathrm{QPH}_{\CC^{*}}(E)$ for which $(w_{1})=z_{i}$, $i=1,2,3,4$,
exhaust the exceptional locus $\mathrm{X}(E)$ in Proposition~\ref{prop:X-Y}, via the $\PP(\Aut(E))$-orbits of elements $(E_{\qp},\Phi)\in\mathrm{QPH}_{\CC^{*}}(E)$ of the form
\begin{gather}\label{eq:HX-1}
\Phi = \begin{pmatrix}
0 & -v_{3}\\
w_{1} & 0
\end{pmatrix}\!,
\end{gather}
such that $(v_{3}) = z_{j} + z_{k} + z_{l}$ for the partition $\{i\}\sqcup\{j,k,l\}=\{1,2,3,4\}$. Moreover, it follows that the $\PP(\Aut(E))$-stabilizers of these points are trivial as well. Finally, it also follows from~\eqref{eq:Phi-1} and Proposition~\ref{prop:S_{j}(E)} that any quasi-parabolic structure in the complement of $\mathrm{U}'(E)$ can be always represented as the projection of an element in either $\mathrm{R}(E)$ or $\mathrm{S}_{0}(E)$, from which the full claim follows.

$(ii)$ When $m_{1}-m_{2}=2$, a point in $\mathrm{QPH}_{\CC^{*}}(E)\sqcup \mathrm{S}_{0}(E)$, satisfies $w_{0}\neq 0$, and is fully determined by $\Phi$. Under a choice of isomorphism $E\cong\calO(m+1)\oplus\calO(m-1)$, the $\PP(\Aut(E))$-orbit of $\Phi$ always contains an element of the form
\[
\Phi = \begin{pmatrix}
0 & -v_{4}\\
w_{0} & 0
\end{pmatrix}\!,
\]
from which the claim follows.

$(iii)$ The result follows from \eqref{eq:Phi-1} and the constraints $w_{m_{2}-m_{1}+2}=0$ and $u_{2}=0$ under a~choice of isomorphism $E\cong \calO(m_{1})\oplus\calO(m_{2})$.
\end{proof}

\section{Combinatorial description of parabolic weight polytopes}\label{sec:poly}

We will present a thorough account of the combinatorial structure on $[0,1]^{4}$ that encodes wall-crossing phenomena for the toy model, and refines an explicit embedding of the convex polytope known as the 4-\textit{demicube} or \textit{demitesseract}.\footnote{The latter was described by Bauer \cite{Bau91} under the $\mathrm{SU}(2)$-constraints (requiring $d$ to be even in the case of 4 marked points), and later by Biswas \cite{Bis98} for parabolic degree 0, in relation to the moduli problem for parabolic bundles of rank 2.} In order to emphasize the features resulting from the inclusion of parabolic Higgs fields into the moduli problem, we will construct this structure from first principles. Consider the functions
\[
\beta_{I}\colon\quad [0,1]^{4}\rightarrow \RR,\qquad I\subset\{1,2,3,4\},
\]
which by definition satisfy the relations $\beta_{I}+\beta_{I^{\mathsf{c}}}=0$ and $\beta_{I'}\leq \beta_{I}$ whenever $I'\subset I$, as well as the bounds
\begin{gather}\label{eq:bounds}
-\vert I^{\mathsf{c}}\vert \leq \beta_{I}\leq \vert I\vert.
\end{gather}
Without any loss of generality (see Remark~\ref{rem:effective}), we will reformulate the necessary and sufficient conditions on parabolic weights for the existence of stable parabolic bundles in terms of $\beta$-weights and the inequalities \eqref{eq:effective}.

\begin{Proposition}[\cite{Bau91, Bis98}]\label{prop:B-B}
There exists a rank $2$ semi-stable parabolic bundle $E_{*}$ of degree $d$ with respect to $\pmb{\beta}\in[0,1)^{4}$ if and only if for any $I\subset\{1,2,3,4\}$ such that $d \equiv \vert I\vert -1\; (\bmod\; 2)$,
\begin{gather}\label{eq:Biswas-ineq}
 \beta_{I}(\pmb{\beta}) \leq \vert I\vert-1.
\end{gather}
$E_{*}$ is necessarily stable if all inequalities are strict.
\end{Proposition}

Proposition \ref{prop:B-B} states necessary and sufficient conditions on $\pmb{\beta}\in[0,1)^{4}$ to grant the semi-stability of generic quasi-parabolic structures on evenly-split bundles, since then
\[
\vert I\vert = d-2\deg(L)+1 = \dim \PP\bigl(H^{0}\big(\CC\PP^{1},E\otimes L^{-1}\big)\bigr)
\]
is the minimum number of flags in any $E_{\qp}$ that can be interpolated by some line sub-bundle $L\subset E$, for which we have that $d \equiv \vert I\vert -1\; (\bmod\; 2)$.
The inequalities \eqref{eq:Biswas-ineq} only depend on the parity of $d$, and each possibility will be treated independently. In the case when $d$ is even, these inequalities are equivalent to the following
interval bounds
\begin{gather}\label{eq:bound1}
-2\leq\beta_{I} \leq 0, \qquad \vert I\vert=1\quad \Leftrightarrow\quad 0\leq\beta_{I'} \leq 2, \qquad \vert I'\vert=3,
\end{gather}
while when $d$ is odd, we obtain the interval bounds
\begin{gather}\label{eq:bound2}
-3 \leq\beta_{\varnothing }\leq -1\Leftrightarrow 1 \leq \beta_{\{1,2,3,4\}} \leq 3
\\ \label{eq:bound3}
-1 \leq \beta_{I} \leq 1,\qquad \vert I\vert=2.
\end{gather}
For both degree parities, the total number of independent interval bounds is equal to four.
Given $I\subset \{1,2,3,4\}$, let $\chi_{I}\colon \{1,2,3,4\}\rightarrow\{0,1\}$ be its characteristic function. Since
\[
\beta_{I}^{-1}(-\vert I^{\mathsf{c}}\vert) = \sum_{i=1}^{4}\chi_{I^{\mathsf{c}}}(i)\,\mathbf{e}_{i},
\]
the map $I\mapsto v_{I}:=\beta^{-1}_{I}(-\vert I^{\mathsf{c}}\vert)$ determines a bijective correspondence between subsets $I\subset \{1,2,3,4\}$ and vertices in $\partial [0,1]^{4}$. A vertex $v_{I}\in
\partial[0,1]^{4}$ will be called \emph{even} (resp.~\emph{odd}) if $\vert I\vert$ is even (resp.~odd). The parity of a vertex is also equal to the parity of its number of nonzero entries. More generally, we have the following convexity result, whose proof is straightforward.

\begin{Lemma}\label{lemma:level-set}
Any $I\subset \{1,2,3,4\}$ induces a partition of the set of vertices in $\partial [0,1]^{4}$ by the level sets of $\beta_{I}$. For every $j=0,\dots,4$, $\beta_{I}^{-1}(j - \vert I^{\mathsf{c}}\vert)$ contains exactly ${4\choose j}$ vertices, namely, those whose Hamming distance to $v_{I}$ is $j$. Each level set is a convex hull of its set of vertices. The partitions induced by $I$ and $I^{\mathsf{c}}$ coincide.
\end{Lemma}

The subset of the power set $P(\{1,2,3,4\})$ containing sets of even (resp.~odd) cardinality will be denoted by $P_{0}(\{1,2,3,4\})$ (resp.~$P_{1}(\{1,2,3,4\})$).

{\sloppy\begin{Definition}
An even (resp.~odd) \textit{partition set} is any subset of $P_{0}(\{1,2,3,4\})$ (resp. $P_{1}(\{1,2,3,4\})$) containing exactly one element of the pair $\{I,I'\}$ for every partition $\{1,2,3,4\}=I\sqcup I'$ by even (resp.~odd) subsets.
\end{Definition}}

\begin{Proposition} \label{prop:W-hull}
Let $\qB_{0}$ $($resp.~$\qB_{1})$ be the convex $4$-polytope in $[0,1]^{4}$ determined by the interval bounds \eqref{eq:bound1} $($resp.~\eqref{eq:bound2}--\eqref{eq:bound3}$)$.
$\qB_{0}$ is the convex hull of all even vertices in $\partial [0,1]^{4}$. Similarly, $\qB_{1}$ is the convex hull of all odd vertices in $\partial [0,1]^{4}$. Consequently, both are isomorphic to the $4$-demicube. $\partial\qB_{\sbul}$ is a triangulation with $16$ tetrahedral cells, parametrized by partition sets $\qI\subset P_{\sbul}(\{1,2,3,4\})$ of given parity, as the convex hulls of the sets $\{v_{I}\in\qB_{\sbul}\colon I\in\qI\}$.
\end{Proposition}

\begin{proof}
The first two statements are straightforward from the explicit form of the interval bounds \eqref{eq:bound1}--\eqref{eq:bound3}.
Moreover, 3-cells in each boundary $\partial\qB_{\sbul}$ correspond to convex hulls of sets of vertices in $\qB_{\sbul}$ whose Hamming distance is exactly 2. The latter are in bijective correspondence with partitions sets $\qI\subset P_{\sbul}(\{1,2,3,4\})$ of the given parity, whose cardinalities are equal to 4, as~$\qI \leftrightarrow \{v_{I}\in\qB_{\sbul}\colon I\in\qI\}$.
\end{proof}

\begin{Definition}
A \emph{semi-stability wall} is any hyperplanar region in $[0,1]^{4}$ of the form
\[
\qH_{I,k}= \beta_{I}^{-1}(k),\qquad k\in\{\vert I\vert -1, \vert I\vert -2, \vert I\vert -3\}.
\]
A wall $\qH_{I,k}$ is \emph{even} (resp.~\emph{odd}) if $k$ is even (resp.~odd). The walls $\qH_{I,\vert I\vert -2}$ will be called \emph{interior}; otherwise, they will be called \emph{boundary}.
\end{Definition}

{\sloppy
The definition of semi-stability walls is based on the condition \eqref{eq:effective} and the interval bounds~\eqref{eq:bounds}. Depending on their parity, semi-stability walls are subsets of either $\qB_{0}$ or $\qB_{1}$. Since
\begin{gather}\label{eq:wall-id-1}
\qH_{I,\vert I\vert -3} = \qH_{I^{\mathsf{c}},\vert I^{\mathsf{c}}\vert -1},
\end{gather}}\noindent
even (resp.~odd) boundary semi-stability walls are in bijection with subsets $I$ of odd (resp.~even) cardinality under the correspondence $I\mapsto\qH_{I,\vert I\vert - 3}$. In addition, there are four interior semi-stability walls for each parity, corresponding to the four partitions $I \sqcup I' = \{1,2,3,4\}$ into subsets of the same cardinality parity, in the form
\begin{gather}\label{eq:wall-id-2}
\qH_{I,\vert I\vert -2} = \qH_{I',\vert I'\vert -2}.
\end{gather}
For both parities, $v_{1/2} := (1/2,1/2,1/2,1/2)$ is the intersection of all four interior semi-stability walls. Tables \ref{tab:wall-even} and \ref{tab:wall-odd} list all semi-stability walls for both parities.
\begin{table}[h!]
\centering
\caption{List of even semi-stability walls.}\label{tab:wall-even}\vspace{1mm}
\tabulinesep=1.2mm
{\small\begin{tabu}{c|c}
 \hline
\multirow{2}{*}{Boundary}
& $\qH_{\{i\},-2} = \qH_{\{j,k,l\},2},\quad \{i\}\sqcup \{j,k,l\}=\{1,2,3,4\}$
\\
\cline{2-2}
& $\qH_{\{j,k,l\},0} = \qH_{\{i\},0},\quad \{i\}\sqcup \{j,k,l\}=\{1,2,3,4\}$
\\
\hline
\multirow{2}{*}{Interior} & $\qH_{\varnothing ,-2} = \qH_{\{1,2,3,4\},2}$
\\
\cline{2-2}
& $\qH_{\{i,j\},0} = \qH_{\{k,l\},0},\quad \{i,j\}\sqcup \{k,l\}=\{1,2,3,4\}$
\\
\hline
\end{tabu}}
\end{table}

\begin{table}[h!]
\centering
\caption{List of odd semi-stability walls.}\label{tab:wall-odd} \vspace{1mm}
\tabulinesep=1mm
{\small\begin{tabu}{c|c}
 \hline
\multirow{3}{*}{Boundary} & $\qH_{\varnothing ,-3} = \qH_{\{1,2,3,4\},3}$
\\
\cline{2-2}
& $\qH_{\{i,j\}, -1} = \qH_{\{k,l\}, 1},\quad \{i,j\}\sqcup \{k,l\}=\{1,2,3,4\}$
\\
\cline{2-2}
& $\qH_{\{1,2,3,4\},1} = \qH_{\varnothing ,-1}$
\\
\hline
\multirow{1}{*}{Interior} & $ \qH_{\{i\},-1} = \qH_{\{j,k,l\},1},\quad \{i\}\sqcup \{j,k,l\}=\{1,2,3,4\}$
\\
\hline
\end{tabu}}
\end{table}

\begin{Definition}
The even (resp.~odd) \textit{parabolic weight polytope} is the refinement $\qW_{0}$ (resp.~$\qW_{1}$) of $[0,1]^{4}$ resulting after the inclusion of all semi-stability walls into $\qB_{0}$ (resp.~$\qB_{1}$). An even (resp.~odd) {open chamber} is the interior of a 4-cell in $\qW_{0}$ (resp.~$\qW_{1}$), or equivalently, any connected component in
\[
(0,1)^{4}\backslash\cup \qH_{I,k},\qquad k\in 2\ZZ\qquad (\text{resp.~} k\in\ZZ\backslash 2\ZZ).
\]
An open chamber inside $\qB_{\sbul}$ is called \emph{interior}, and \emph{exterior} otherwise. A tetrahedral cell in $\partial\qB_{\sbul}$ is of \textit{type A} if it is also contained in $\partial[0,1]^{4}$, and of \textit{type B} otherwise.
\end{Definition}

\begin{Proposition}
There are exactly 8 tetrahedral cells of type $A$ in $\partial\qB_{\sbul}$.
Tetrahedral cells of type $B$ correspond to the $8$ boundary semi-stability walls of a given parity.
\end{Proposition}

\begin{proof}
A tetrahedral cell belongs to $\partial\qB_{\sbul}\cap\partial[0,1]^{4}$ if and only if it is the convex hull of a set of 4 vertices of a given parity with a fixed coordinate component value $\beta_{i}=0,1$, $i=1,2,3,4$, yielding 8 possibilities for each parity. It follows from Lemma~\ref{lemma:level-set} and Proposition~\ref{prop:W-hull} that the 8 remaining tetrahedral cells in $\partial\qB_{\sbul}$ are exhausted by the boundary semi-stability walls $\qH_{I,\vert I\vert -3}$ for all $I\subset\{1,2,3,4\}$ of a given parity.
\end{proof}

\begin{Corollary}[classification of exterior chambers]\label{cor:ex-chambers}
At most one inequality in either \eqref{eq:bound1} or~\eqref{eq:bound2}--\eqref{eq:bound3} can fail to hold on $(0,1)^{4}$, and each boundary semi-stability wall $\qH_{I,\vert I\vert - 3}$ determines a unique exterior chamber whose closure contains it. Every exterior chamber to~$\qB_{0}$ $($resp.~$\qB_{1})$ is the interior of the convex hull of $\qH_{I,\vert I\vert - 3}$ and the vertex $v_{I}$ for a unique $I\subset \{1,2,3,4\}$ of odd $($resp.~even$)$ cardinality.
\end{Corollary}

\begin{proof}
The first claim follows from the convexity of the regions $\beta_{I} < \vert I\vert - 3$. Similarly, it follows from the identity \eqref{eq:wall-id-1} that there is a bijective correspondence between exterior chambers and regions $\beta^{-1}_{I}((-\vert I^{\mathsf{c}}\vert ,\vert I\vert - 3))$ for each $I\subset\{1,2,3,4\}$, in such a way that their chamber parity is given by the parity of $\vert I\vert -3$.
\end{proof}

On the other hand, interior chambers on each parabolic weight polytope are effectively classified in terms of the convex geometry of interior semi-stability walls.
This is formalized in the next two results.

 \begin{Proposition}\label{prop:int-w}
Every $2$-cell in the boundary of a given tetrahedral cell in $\partial\qB_{\sbul}$ is the intersection of the latter and an interior semi-stability wall of the given parity. Consequently, the restriction $\qW_{\sbul}\vert_{\qB_{\sbul}}$ is equal to the pyramid of $\qB_{\sbul}$ over the apex $v_{1/2}$.
\end{Proposition}

\begin{proof}
By Lemma~\ref{lemma:level-set}, the interior semi-stability wall corresponding to a partition $I\sqcup I' = \{1,2,3,4\}$ as in \eqref{eq:wall-id-2} is the convex hull of 6 vertices, obtained by removing $\{v_{I}, v_{I'}\}$ (of pairwise Hamming distance 4) from the set of 8 vertices of given parity. In this set of 6 vertices, there exist exactly 8 triples of vertices whose pairwise Hamming distance is equal to 2, so that each triple spans a 2-cell in $\partial \qB_{\sbul}$. It follows from Proposition~\ref{prop:W-hull} that every tetrahedral cell in $\partial \qB_{\sbul}$ is uniquely expressed as the convex hull of one of these triples and either $v_{I}$ or $v_{I'}$, in such a way that its intersection with each of the 4 interior semi-stability walls in $\qB_{\sbul}$ retrieves the 2-cells in its boundary. The second claim is immediate from the convexity of $\qB_{\sbul}$, since $v_{1/2}$ is the common intersection of all interior semi-stability walls of given parity.
\end{proof}

\begin{Corollary}[classification of interior chambers]\label{cor:in-chambers}
Every interior chamber is the intersection of a choice of side regions $\{\beta_{I} < \vert I\vert -2 \colon I\in\qI\}$ to all interior semi-stability walls for a unique partition set $\qI\subset P_{\sbul}(\{1,2,3,4\})$. Equivalently, it is the interior of the convex hull of the apex~$v_{1/2}$ and the tetrahedral cell in $\partial\qB_{\sbul}$ spanned by $\{v_{I'}\colon I'\in\qI\}$.
\end{Corollary}
\begin{proof}
It follows from Proposition~\ref{prop:int-w} that interior chambers in each polytope $\qB_{\sbul}$ can be parametrized by a choice of side of all interior semi-stability walls. Any of these choices is defined and parametrized by a unique partition set under the correspondence
\[
\qI \mapsto \bigcap_{I\in\qI} \qR_{I},\qquad
\qR_{I} := \beta^{-1}_{I}\bigl((-\vert I^{\mathsf{c}}\vert,\vert I\vert -2)\bigr),
\]
given that for any partition $\{1,2,3,4\}=I\sqcup I'$, the open regions $\qR_{I}$ and $\qR_{I'}$ satisfy
\[
v_{I}\in\overline{\qR_{I}},\; v_{I'}\in\overline{\qR_{I'}},\qquad
\overline{\qR_{I}}\cup \overline{\qR_{I'}} =[0,1]^{4},\qquad
\overline{\qR_{I}}\cap \overline{\qR_{I'}} = \qH_{I,\vert I\vert -2} = \qH_{I',\vert I'\vert -2},
\]
from which the first claim follows. Similary, it follows from Proposition~\ref{prop:W-hull} that the closure of any interior chamber contains a unique tetrahedral cell in $\partial\qB_{\sbul}$, which is equal to the intersection of $\partial\qB_{\sbul}$ and the former. Any additional vertex in $\qW_{\sbul}$ would necessarily arise from the intersection of interior semi-stability walls of a given parity. Since their common intersection equals $v_{1/2}$, the second claim follows.
\end{proof}

We will stick to the following notational convention. Interior chambers will be denoted as~$\qC_{\qI}$, where $\qI$ is the partition set parametrizing each of them under the previous correspondence. On~the other hand, the exterior chamber to the boundary semi-stability wall $\qH_{I,\vert I\vert - 3}$ will be denoted as $\qC_{I}$.

The intersection of $\partial \qB_{\sbul}$ and the closure of an open chamber of the same parity is a tetrahedral cell and characterizes an interior chamber uniquely. Since tetrahedral cells of type $B$ are exterior semi-stability walls, they characterize uniquely both an interior and an exterior chamber as the common boundary of their closures. The type of an interior chamber in $\qW_{\sbul}$ (or equivalently, its associated partition set) is defined as the type of its tetrahedral cell. Reflection along an exterior semi-stability wall $\qH_{I,\vert I\vert - 3}$ exchanges the vertex $v_{I}$ and the apex $v_{1/2}$, and determines a bijection between its exterior chamber and its interior chamber of type $B$.

We will say that two interior chambers $\qC_{\qI}$ and $\qC_{\qI'}$ are \textit{neighboring} if a reflection along a~semi-stability wall $\qH_{I,\vert I\vert - 2} = \qH_{I',\vert I'\vert - 2}$ transforms one into the other.
The proof of the next corollary is immediate from Corollaries~\ref{cor:ex-chambers} and~\ref{cor:in-chambers}.

\begin{Corollary}[combinatorial wall-crossing]\label{cor:corresp}
Every interior semi-stability wall of given parity $\qH_{I,\vert I\vert -2} = \qH_{I',\vert I'\vert -2}$ is the union of $4$ tetrahedral cells in $\qW_{\sbul}$, each of which is the common boundary of one of the four possible pairs of neighboring interior chambers $\{\qC_{\qI},\qC_{\qI'}\}$ of opposite type whose partitions sets satisfy
\[
\qI\backslash\{I\} = \qI'\backslash\{I'\}.
\]
Every even $($resp.~odd$)$ exterior semi-stability wall $\qH_{I,\vert I\vert - 3}$ is the common boundary between the exterior chamber $\qC_{I}$ and the even $($resp.~odd$)$ interior chamber $\qC_{\qI(I)}$ of type $B$,
where
\[
I'\in\qI(I) \Leftrightarrow v_{I'}\in\qH_{I,\vert I\vert - 3}.
\]
\end{Corollary}

Tables \ref{tab:partition-even} and \ref{tab:partition-odd} provide explicit lists of partition sets for each parity, together with the tetrahedral cells in $\partial\qB_{\qp}$ for the corresponding open chambers they parametrize, according to Proposition~\ref{prop:W-hull}. Notice that the collections of partition sets of a given type are closed under the operation of taking complements in $P_{\sbul}(\{1,2,3,4\})$.

\begin{table}[h!]
\centering
\caption{Even partition sets, $i=1,2,3,4$, $\{j,k,l\}=\{1,2,3,4\}\backslash\{i\}$.}\vspace{1mm}\label{tab:partition-even}
\tabulinesep=1.4mm
{\small\begin{tabu}{c|c|c}
\hline
& Partition set $\qI$ & Tetrahedral cell in $\partial\qB_{0}$\\
\hline
\multirow{2}{*}{Type $A$} & $\{\varnothing , \{k,l\}, \{j,l\}, \{j,k\}\}$ & $\qB_{0}\cap\{\beta_{i} = 1\}$\\
\cline{2-3}
& $\{\{1,2,3,4\}, \{i,j\}, \{i,k\}, \{i,l\}\}$ & $\qB_{0}\cap\{\beta_{i} = 0\}$\\
\hline
\multirow{2}{*}{Type $B$} & $\{\varnothing , \{i,j\}, \{i,k\}, \{i,l\}\}$ & $\qH_{\{i\},-2}$\\
\cline{2-3}
& $\{\{1,2,3,4\}, \{k,l\}, \{j,l\}, \{j,k\}\}$ & $\qH_{\{j,k,l\},0}$\\
\hline
\end{tabu}}
\end{table}

\begin{table}[h!]
\centering
\caption{Odd partition sets, $i=1,2,3,4$, $\{j,k,l\}=\{1,2,3,4\}\backslash\{i\}$.}\vspace{1mm}\label{tab:partition-odd}
\tabulinesep=1.4mm
{\small\begin{tabu}{c|c|c}
\hline
& Partition set $\qI$ & Tetrahedral cell in $\partial\qB_{1}$ \\
\hline
\multirow{2}{*}{Type $A$} & $\{\{i\}, \{i,k,l\}, \{i,j,l\}, \{i,j,k\}\}$ & $\qB_{1}\cap\{\beta_{i} = 0\}$\\
\cline{2-3}
& $\{\{j,k,l\}, \{j\}, \{k\}, \{l\}\}$ & $\qB_{1}\cap\{\beta_{i} = 1\}$\\
\hline
\multirow{3}{*}{Type $B$} & $\{\{1\}, \{2\}, \{3\}, \{4\}\}$ & $\qH_{\varnothing ,-3}$\\
\cline{2-3}
& $\{\{i\}, \{j\}, \{i,j,l\}, \{i,j,k\}\}$ & $\qH_{\{i,j\},-1}$\\
\cline{2-3}
& $\{\{2,3,4\}, \{1,3,4\}, \{1,2,4\}, \{1,2,3\}\}$ & $\qH_{\{1,2,3,4\},1}$\\
\hline
\end{tabu}}
\end{table}

The group of symmetries of $[0,1]^{4}$ is called \textit{hexadecachoric group}, denoted as $B_{4}$. It is a~group of order 384 generated by permutation of coordinates and reflections along the hyperplanes
\mbox{$\beta_{i}=1/2$}, $i=1,2,3,4$. Its restriction to reflections along pairs of hyperplanes leads to a~distinguished index 2 subgroup, the Coxeter group $D_{4}$.

\begin{Corollary}\label{cor:reflections}
$D_{4}$ is the group of symmetries of both $\qW_{0}$ and $\qW_{1}$. Both polytopes are isomorphic under the action of the group $B_{4}/D_{4}\cong\ZZ_{2}$. Moreover, $D_{4}$ acts transitively on the sets of exterior chambers, and on the sets of interior chambers of type $A$ and $B$ respectively. The stabilizer of any open chamber is the permutation group of the vertices of its characteristic tetrahedral cell.
\end{Corollary}

\begin{Definition}
The \textit{wall-crossing graph} of $\qW_{\sbul}$ is the 1-skeleton of its dual polytope. It assigns a vertex to every open chamber in $\qW_{\sbul}$, and an edge to any pair of open chambers whose closures intersect at a tetrahedral cell.
\end{Definition}

\begin{figure}[!ht]
\centering
\includegraphics[width=3.8in]{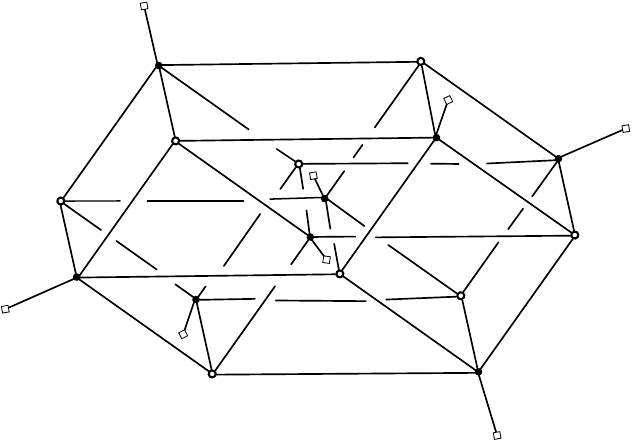}
\caption{The wall-crossing graph for both $\qW_{0}$ and $\qW_{1}$. Square vertices correspond to the 8 exterior open chambers. White round vertices correspond to the 8 interior open chambers of type $A$, while black round vertices to the 8 interior open chambers of type $B$.}\label{fig:wall-cross}
\end{figure}

{\samepage
The wall-crossing graph of $\qW_{\sbul}$ inherits a coloring from the three possible types of open chambers. The qualitative type of wall-crossing in $\qW_{\sbul}$ is then measured by the possible vertex types that are connected by a given edge in the wall-crossing graph. Up to the action of the symmetry group $D_{4}$, there are only two qualitative types of wall-crossing: $(i)$ between an exterior chamber and its neighboring chamber of type $B$, whose closures intersect in a boundary wall, and $(ii)$
between neighboring chambers of type $A$ and $B$, whose closures intersect in an interior wall (Figure~\ref{fig:wall-cross}).

}

\section{Conditional stability and basic building blocks}\label{sec:CS}

The construction of geometric models for all possible Harder--Narasimhan strata as spaces of stable orbits under $\PP(\Aut(E))$-actions, as well as the characterization of their wall-crossing under variation of parabolic weights, will be achieved by linking the classification of \textit{conditionally stable strata} (Definition~\ref{def:cond-stab}) to the combinatorics of weight polytopes. This problem depends on a choice of degree parity for a rank 2 vector bundle. We will consider the two possibilities simultaneously.

In order to classify the geometry and topology of all orbit spaces of conditionally stable strata, which in general contain non-Hausdorff loci even after restricting the type of invariant sub-bundles, their ``taming" is achieved by dissecting them into elementary Hausdorff irreducible components. The details are summarized in Definition~\ref{def:b-b-b} and Propositions~\ref{prop:CS}, \ref{prop:line-orbit} and \ref{prop:C-stable}. We begin with the following corollary.

\begin{Corollary}\label{cor:stable-open}
If $m_{1}-m_{2}\leq 2$, then for any $\pmb{\beta} := (\beta_{1},\beta_{2},\beta_{3},\beta_{4})\in(0,1)^{4}$, the set $\mathrm{QPH}^{s}_{\pmb{\beta}}(E)$ of stable elements in $\mathrm{QPH}(E)$ is nonempty, Zariski open, and contains $\mathrm{QPH}_{\CC^{*}}(E)$. If $m_{1}-m_{2}\geq 4$, then $\mathrm{QPH}^{s}_{\pmb{\beta}}(E) = \varnothing $ with respect to any $\pmb{\beta}\in[0,1)^{4}$.
\end{Corollary}

\begin{proof}
The first claim follows from Lemmas~\ref{lemma:nilpotent} and \ref{lemma:cases}, since $\mathrm{QPH}_{\CC^{*}}(E)\neq\varnothing $ if and only if $m_{1}-m_{2}\leq 2$. If $m_{1}-m_{2} \geq 4$, Lemma~\ref{lemma:cases} implies that $\mathrm{QPH}(E) = \mathrm{QPH}_{0}(E) = \mathrm{Q}(E)\sqcup \mathrm{R}(E)$, and every element in $\mathrm{QPH}(E)$ would be destabilized by $E_{1}$.
\end{proof}

\begin{Remark}
The absence of weight constraints for the existence of moduli spaces of semi-stable parabolic Higgs bundles in genus 0, following from Corollary~\ref{cor:stable-open}, is already implicit, in terms of the general non-abelian Hodge correspondence, in Simpson's solution of the Deligne--Simpson problem \cite{Sim92}. This is a sharp contrast with moduli spaces of semi-stable parabolic bundles \cite{Bau91, Bis98}, for which $\Phi\equiv 0$.
\end{Remark}

Corollary \ref{cor:stable-open} states that for even degree $d=2m$, the only Birkhoff--Grothendieck splitting types that can support stable parabolic Higgs bundles are
\[
\calO(m)\oplus \calO(m) \qquad\text{and}\qquad \calO(m+1)\oplus\calO(m-1),
\]
while for odd degree $d=2m+1$, these are
\[
\calO(m+1)\oplus \calO(m)\qquad\text{and}\qquad \calO(m+2)\oplus\calO(m-1).
\]
It follows from general principles that stability of quasi-parabolic Higgs bundles is an invariant along a given open chamber, and moreover, that the subset of points in $\mathrm{QPH}(E)$ that are stable with respect to every choice of open chamber in $\qW_{\sbul}$ is Zariski open. Nevertheless, we will verify that the intersection of the former Zariski open subset with $\mathrm{QPH}_{0}(E)$ is always empty with respect to any choice of degree parity, i.e., it either coincides with $\mathrm{QPH}_{\CC^{*}}(E)$ when $m_{1}-m_{2}\leq 2$, or is empty when $m_{1}-m_{2}\geq 3$. For any choice of open chamber $\qC\subset\qW_{\sbul}$, we will denote the stable locus associated to any $\pmb{\beta}\in\qC$ by $\mathrm{QPH}^{s}_{\qC}(E)$.

\begin{Definition}\label{def:cond-stab}
A point in $\mathrm{QPH}(E)$ will be called \textit{conditionally stable} if it is stable with respect to some choice of open chamber $\qC\subset\qW_{\sbul}$ of corresponding degree parity. The different substrata of conditionally stable quasi-parabolic Higgs bundles in $\mathrm{QPH}_{0}(E)$ will be accordingly denoted by $\mathrm{Q}^{\rm cs}(E)$, $\mathrm{R}^{\rm cs}(E)$, and $\mathrm{S}^{\rm cs}_{j}(E)$.
\end{Definition}

The subset of conditionally stable parabolic Higgs bundles in $\mathrm{QPH}(E)$ is equal to the complement of the Zariski closed subset of points that are either unstable or strictly semi-stable for any choice of weights $\pmb{\beta}\in[0,1)^{4}$. By Corollary~\ref{cor:stable-open}, the latter is a subset of $\mathrm{QPH}_{0}(E)$. The next lemma is the key to characterize conditionally stable loci systematically. The characterization is performed in Proposition~\ref{prop:CS}.

\begin{Lemma}\label{lemma:not-stable}
$(E_{\qp},\Phi)\in\mathrm{QPH}_{0}(E)\backslash\mathrm{Q}(E)$ is not stable with respect to any $\pmb{\beta}\in[0,1)^{4}$ if and only~if
\begin{gather}\label{eq:not-stable}
-\big\vert I_{L(\Phi),E_{\qp}}^{\mathsf{c}}\big\vert \geq \deg(E)-2\deg(L(\Phi)).
\end{gather}
A point $(E_{\qp},0)\in\mathrm{Q}(E)$ is not stable with respect to any $\pmb{\beta}\in[0,1)^{4}$ if and only if there is a pair of line sub-bundles $L,L'\subset E$, $L\cong\calO(m_{1})$ and $L'\cong\calO(m_{2})$, satisfying
\begin{gather}\label{eq:L-L'}
\{1,2,3,4\} = I_{L,E_{\qp}}\cup I_{L',E_{\qp}}.
\end{gather}
\end{Lemma}

\begin{proof}
Since $-\vert I^{\mathsf{c}}\vert\leq \beta_{I}$, the inequality \eqref{eq:not-stable} for some $(E_{\qp},\Phi)\in\mathrm{QPH}_{0}(E)\backslash\mathrm{Q}(E)$ implies that the inequality \eqref{eq:effective} would never be satisfied in the strict sense for any given choice of $\pmb{\beta}\in[0,1)^{4}$. The converse statement follows after taking the infimum of \eqref{eq:effective} over the values of $\beta_{I_{L(\Phi),E_{\qp}}}$.

Similarly, the condition \eqref{eq:L-L'} implies that for any $\pmb{\beta}\in[0,1)^{4}$, $(E_{\qp},0)\in\mathrm{Q}(E)$ is not stable with respect to at least one sub-bundle of the pair $\{L,L'\}$. On the other hand, if \eqref{eq:not-stable} holds for a point $(E_{\qp},0)\in\mathrm{Q}(E)$ and some $L\cong \calO(m_{1})$, then it follows from the properties of line bundle interpolation that there always exists $L'\cong \calO(m_{2})$ such that the pair $L,L'$ satisfies \eqref{eq:L-L'}. Therefore, to prove the remaining implication, it suffices to assume that $(E_{\qp},0)\in\mathrm{Q}(E)$ is not stable for any $\pmb{\beta}\in[0,1)^{4}$ while \eqref{eq:not-stable} doesn't hold for any $L\cong \calO(m_{1})$.

We will treat the possibilities $m_{1}=m_{2}$ and $m_{1}>m_{2}$ separately. If we assume that $m_{1}>m_{2}$, instability of $(E_{\qp},0)$ for all $\pmb{\beta}\in[0,1)^{4}$, and $\vert I_{E_{1},E_{\qp}}\vert < 4 -m_{1}+m_{2}$, then the instability of $(E_{\qp},0)$ for $\beta_{I_{E_{1},E_{\qp}}} < m_{2}-m_{1}$ implies the existence of a destabilizing $L'\cong\calO(m_{2})$, and we obtain
\[
I_{L',E_{\qp}} = I^{\mathsf{c}}_{E_{1},E_{\qp}}.
\]
In particular, \eqref{eq:L-L'} holds. Similarly, if we assume that $m_{1}=m_{2}=m$ and $L\cong\calO(m)$ with $\vert I_{L,E_{\qp}}\vert <4$, then the instability of $(E_{\qp},0)$ for all $\beta_{I_{L,E_{\qp}}} < 0$ implies the existence of $L'\cong \calO(m)$ such that
\[
I_{L',E_{\qp}} = I^{\mathsf{c}}_{L,E_{\qp}}.
\]
The case when $(F_{1},F_{2},F_{3},F_{4})\in\mathrm{B}_{\{1,2,3,4\},m}(E)$ holds trivially, since this is the only instance for which $I_{L,E_{\qp}}\cap I_{L',E_{\qp}} \neq\varnothing $, and then \eqref{eq:L-L'} holds if we let $L' = L$.
\end{proof}

\begin{Proposition}\label{prop:CS}
Let\footnote{The notation for $\mathrm{S}_{I,j}(E)$ is redundant when $2(\deg(L(\Phi)) +1) - \deg(E) = 4$, that is, when $E\cong\calO(m)\oplus\calO(m)$ or $E\cong\calO(m+1)\oplus\calO(m-1)$ and $j=m-1$, since then $\mathrm{S}_{\{1,2,3,4\},m-1}(E) = \mathrm{S}_{m-1}(E)$.}
\[
\mathrm{R}_{I}(E):=\mathrm{R}(E)\cap \mathrm{par}^{-1}(\mathrm{B}_{I,m_{1}}(E)),\qquad
\mathrm{S}_{I,j}(E):=\mathrm{S}_{j}(E)\cap\mathrm{par}^{-1}(\mathrm{B}_{I,j}(E)).
\]
\begin{enumerate}\itemsep=0pt
\item[$(i)$] The different types of conditionally stable substrata are classified as follows:
\begin{gather*}
\mathrm{Q}^{\rm cs}(E)=
\begin{cases}
\iota(\mathrm{U}(E)) & \text{if} \quad m_{1} - m_{2} \leq 2,
\\
\varnothing & \text{if} \quad m_{1} -m_{2} = 3,
\end{cases}
\\
\mathrm{R}^{\rm cs}(E)=
\bigsqcup_{\vert I\vert \leq 3-m_{1}+m_{2}}\mathrm{R}_{I}(E),
\\
\mathrm{S}^{\rm cs}_{j}(E)=
\begin{cases}
\bigsqcup_{\vert I\vert = 2,3}\mathrm{S}_{I,m}(E) & \text{if} \quad m_{1} = m_{2}=m,\, j=m,
\\
\mathrm{S}_{j}(E) & \text{otherwise.}
\end{cases}
\end{gather*}
\item[$(ii)$] An element in $\mathrm{Q}^{\rm cs}(E)$ is stable with respect to any choice of interior open chamber if and only if it lies in $\iota(\mathrm{U}'(E))$.
\item[$(iii)$] $\PP(\Aut(E))$ acts freely on $\mathrm{R}^{\rm cs}(E)$. If $m_{1}=m_{2}=m$ then $\PP(\Aut(E))$ acts freely on $\mathrm{S}^{\rm cs}_{m}(E)$ and $\PP(\Aut(E))\times \CC^{*}$ acts freely and transitively on $\mathrm{S}^{\rm cs}_{m-1}(E)$. If $m_{1}-m_{2}=1,2$ then $\PP(\Aut(E))$ acts freely and transitively in the factor $\mathrm{P}_{m_{2}}(E)$ of $\mathrm{S}_{m_{2}}(E)$.
 \end{enumerate}
\end{Proposition}

\begin{proof}
The first and second claims follow from Lemma~\ref{lemma:not-stable}. The third claim follows from Proposition~\ref{prop:S_{j}(E)}, the definition of the subspaces $\mathrm{B}_{I,j}(E)\subset \mathrm{QP}(E)$ in \eqref{eq:B_{I,j}(E)} together with property~\eqref{eq:I<I'}, the definition of the strata $\mathrm{B}_{I,m_{1}}(E)$ when $m_{1}>m_{2}$, and the specializations of the canonical form \eqref{eq:Phi-2} (including \eqref{eq:Phi(m-1)}) in each case. Notice that $\bigsqcup_{\vert I\vert \leq 2}\mathrm{R}_{I}(E) = \bigsqcup_{\vert I\vert =1,2}\mathrm{R}_{I}(E)$ when $m_{1}-m_{2} =1$, since $\mathrm{R}_{\varnothing }(E) = \varnothing $.
\end{proof}

\begin{Definition}\label{def:b-b-b}
The (possibly empty) \textit{basic building blocks} of a subset $I\subset\{1,2,3,4\}$, denoted by $\cK_{I}$ (resp.~$\cL_{I}$ whenever $j=m_{2}$), are the $\PP(\Aut(E))$-orbit spaces of $\mathrm{R}_{I}(E)$ (resp.~$\mathrm{S}_{I,j}(E)$). In turn, the individual orbits $\cN_{I}$ are unambiguously defined as
\[
\cN_{I}:=\iota(\mathrm{B}_{I,j}(E)\cap \mathrm{U}(E))/\PP(\Aut(E)).
\]
When $E$ is evenly-split we will denote the ``bulk" of orbits of quasi-parabolic bundles in $\mathrm{QP}(E)$ that are stable with respect to \textit{any interior open} chamber as
\[
\cN_{\qB_{\sbul}}:= \iota(\mathrm{U}'(E))/\PP(\Aut(E))\cong\CC\PP^{1}\backslash\{z_{1},0,1,\infty\}.
\]
\end{Definition}

Propositions \ref{prop:S_{j}(E)} and \ref{prop:CS} imply that the orbit space of $\mathrm{S}^{\rm cs}_{m-1}(E)$ consists of a single basic building block, which is a $\PP(\Aut(E))$-principal bundle with base a point when $E\cong \calO(m+1)\oplus\calO(m-1)$ and $\CC^{*}$ when $E\cong\calO(m)\oplus\calO(m)$. We will denote these as $\cL_{\{1,2,3,4\},0}$ and $\cL_{\{1,2,3,4\},1}$ respectively. Even though they are associated to different bundle splitting types, they glue into a punctured line (Proposition~\ref{prop:glue}). For notational simplicity, we will denote the resulting punctured line by $\cL_{\{1,2,3,4\}}$.

\begin{Proposition}\label{prop:line-orbit}\quad
\begin{enumerate}\itemsep=0pt
\item[$(i)$] There is a bijective correspondence between basic building blocks $\cN_{I}$ (resp.~$\cK_{I}$, resp.~$\cL_{I}$) and subsets $I\subset\{1,2,3,4\}$.
\item[$(ii)$] When nonempty, the action of $\PP(\Aut(E))$ on $\mathrm{S}_{I,m_{2}}(E)$ and $\mathrm{R}_{I}(E)$ is free and proper. Apart from $\cL_{\{1,2,3,4\},0} = \mathrm{pt}$ and $\cL_{\{1,2,3,4\},1}\cong\CC^{*}$, the nonempty basic building blocks are classified as follows:
\begin{gather*}
\cK_{I} \cong\begin{cases}
\CC, & \vert I\vert + m_{1}-m_{2} = 2,
\\
\CC\PP^{1}, & \vert I\vert + m_{1} - m_{2} = 3,
\end{cases}
\\
\cL_{I} \cong\begin{cases}
\CC,& \vert I\vert - m_{1} + m_{2} = 2,
\\
\CC\PP^{1}, & \vert I\vert - m_{1} + m_{2} = 3.
\end{cases}
\end{gather*}
Moreover, for any $\vert I\vert = 3 - m_{1} + m_{2}$, there is an isomorphism
\[
\cK_{I}\sqcup \bigg\{\bigsqcup_{I'\subset I,\; \vert I\backslash I'\vert = 1}\cK_{I'}\bigg\}\cong \mathrm{Bl}_{I}\big(\CC\PP^{1}\big)
\]
and for any $\vert I\vert = 3 + m_{1} - m_{2}$, there is an isomorphism
\[
\cL_{I}\sqcup\bigg\{\bigsqcup_{I'\subset I,\; \vert I\backslash I'\vert = 1}\cL_{I'}\bigg\}\cong \mathrm{Bl}_{I}\big(\CC\PP^{1}\big).
\]
\end{enumerate}
\end{Proposition}

\begin{proof}
The first claim is already implicit in Proposition~\ref{prop:U} for the basic building blocks of the form $\cN_{I}$, and the remaining possibilities are implicit in the second claim. In terms of this correspondence, the resulting orbits are individual points.
The second claim is a consequence of the $\PP(\Aut(E))$-equivariance of $\mathrm{par}$ combined with Propositions~\ref{prop:S_{j}(E)} and \ref{prop:CS}. When $I$ is such that $\vert I\vert = 3 - m_{1} + m_{2}$ (resp.~$\vert I\vert = 3 + m_{1} - m_{2}$), the divisor $(\Phi)$ associated to any element in $\mathrm{R}_{I}(E)$ (resp.~$\mathrm{S}_{I,m}(E)$) takes the form
\[
(\Phi) = z+ \sum_{j\in I^{\mathsf{c}}} z_{j},
\]
in such a way that if $z=z_{i}$ for any $i\in I$, then $F_{i} = L(\Phi)\vert_{z_{i}}$. In this case, the correspondence
\[
(E_{\qp},\Phi)\mapsto z
\]
determines the isomorphism $\cK_{I}\cong\CC\PP^{1}$ (resp.~$\cL_{I}\cong \CC\PP^{1}$). Moreover, when $I'\subset I$ is such that $\vert I'\vert = 2 - m_{1} + m_{2}$ (resp.~$\vert I'\vert = 2 + m_{1} - m_{2}$), letting $I\backslash I' = \{i\}$ we have that $z =z_{i}$ and $F_{i} \neq L(\Phi)\vert_{z_{i}}$. Then
\[
(E_{\qp},\Phi)\mapsto F_{i}
\]
gives a biholomorphism $\cK_{I'}\cong\CC$ (resp.~$\cL_{I'}\cong \CC$), and moreover, the biholomorphism
\[
\cK_{I}\sqcup\cK_{I'}\cong\mathrm{Bl}_{\{i\}}\big(\CC\PP^{1}\big)\qquad \big(\text{resp.}\ \cL_{I}\sqcup\cL_{I'}\cong\mathrm{Bl}_{\{i\}}\big(\CC\PP^{1}\big)\big),
\]
from which the claim follows.
\end{proof}

\begin{figure}[!ht]
\centering
\includegraphics[width=4.2in]{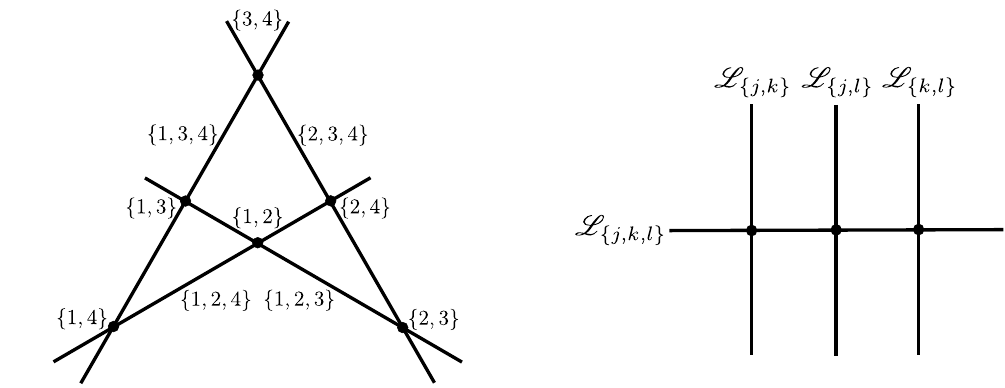}
\put(-230,-13){\makebox(0,0)[lb]{\small$a$}}
\put(-52,-13){\makebox(0,0)[lb]{\small$b$}}
\caption{$(a)$ The complete quadrilateral or $(6_{2}4_{3})$-configuration. $(b)$ Maximal Hausdorff unions of basic building blocks in the orbit space of $\mathrm{S}^{\rm cs}_{0}(E)$ when $E\cong \calO(m)\oplus\calO(m)$, parametrized by the sets $\{\{j,k\},\{j,l\},\{k,l\},\{j,k,l\}\}$, $j\neq k$, $j\neq l$, $k\neq l$.
}\label{fig:S-0-conf}
\end{figure}

\begin{figure}[!ht]
\centering
\includegraphics[width=4.2in]{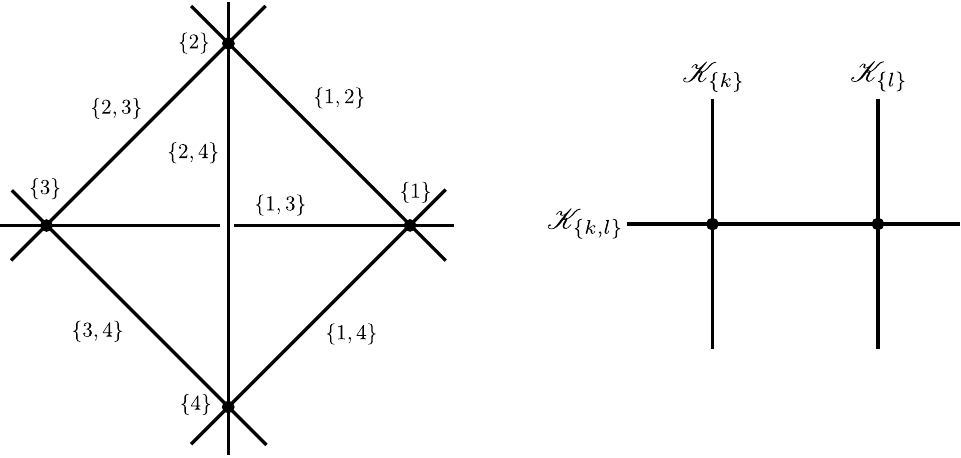}
\put(-235,-13){\makebox(0,0)[lb]{\small$a$}}
\put(-55,-13){\makebox(0,0)[lb]{\small$b$}}
\caption{$(a)$ The complete quadrangle or $(4_{3}6_{2})$-configuration. $(b)$ Maximal Hausdorff union of basic building blocks in the orbit space of $\mathrm{R}^{\rm cs}(E)$ when $E\cong \calO(m+1)\oplus\calO(m)$, parametrized by the sets $\{\{k\},\{l\},\{k,l\}\}$, $k\neq l$.
}\label{fig:R-1-conf}
\end{figure}

\begin{figure}[!ht]
\centering
\includegraphics[width=3.8in]{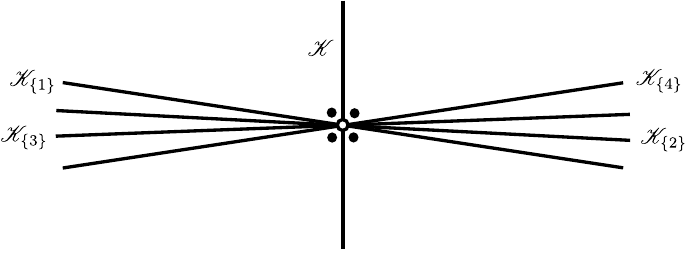}
\put(-144.5,75){\makebox(0,0)[lb]{\tiny$\varnothing$}}
\caption{The (non-Hausdorff) orbit space of $\mathrm{R}^{\rm cs}(E)$ when $E\cong \calO(m+1)\oplus\calO(m-1)$, formed by union of the punctured line $\cK_{\varnothing }$ and the quadruple line $\cup_{i=1}^{4} \cK_{\{i\}}$, with maximal Hausdorff subspaces $\cK_{\{i\}}\sqcup\cK_{\varnothing }\cong \mathrm{Bl}_{\{i\}}(\CC\PP^{1})$, $i=1,2,3,4$.}\label{fig:R-2-conf}
\end{figure}

There are three non-Hausdorff orbit spaces of conditionally stable strata in $\mathrm{QPH}_{0}(E)\backslash\mathrm{Q}(E)$, \textit{resulting from fixing a type of kernel lines}, as a consequence of the pathological orbit space topology of $\mathrm{U}(E)\backslash\mathrm{U}'(E)$ and $\mathrm{Y}(E)$. These are parametrized by the incidence properties of certain projective plane configurations.
When $m_{1} = m_{2}$, the topology of the union of basic building blocks $\cL_{I}$ is determined by the incidence properties of the complete quadrilateral $(6_{2}4_{3})$ (Figure~\ref{fig:S-0-conf}), with all basic building blocks of the form $\cL_{\{j,k,l\}}$ (resp.~all pairs of the form $\{\cL_{\{i,j\}},\cL_{\{k,l\}}\}$ for any partition $\{i,j\}\sqcup\{k,l\} = \{1,2,3,4\}$) identified under the quotient topology in the full orbit space. When $m_{1} - m_{2} = 1$, the same is true for all basic building blocks~$\cK_{I}$ and the incidence properties of the complete quadrangle $(4_{3}6_{2})$ (Figure~\ref{fig:R-1-conf}), with all basic building blocks of the form $\cK_{\{k,l\}}$ being identified under the quotient topology in the full orbit space. When $m_{1} - m_{2} = 2$, the basic building blocks $\cK_{I}$ are superimposed as an identification of four blow-ups of a line at a point (Figure~\ref{fig:R-2-conf}). A full list of conditionally stable loci grouped by $\mathrm{QPH}_{0}(E)$-stratum type, and their induced orbit space type, is compiled in
Table~\ref{tab:master-kit}.

\begin{table}[h!]
\centering
\caption{Master kit of orbit spaces for all possible conditionally stable loci, following from Pro\-po\-si\-tions~\ref{prop:CS} and \ref{prop:line-orbit}.}\vspace{1mm}\label{tab:master-kit}
\tabulinesep=1.2mm
{\small\begin{tabu}{c|c|c}
\hline
$m_{1}-m_{2}$ & Conditionally stable locus & Orbit space
\\
\hline
\multirow{3}{*}{$0$} & $\mathrm{Q}^{\rm cs}(E) = \iota(\mathrm{U}(E))$ & $\CC\PP^{1}$ w. 3 double points
\\
\cline{2-3}
& $\mathrm{S}^{\rm cs}_{m}(E) = \bigsqcup_{\vert I\vert = 2,3}\mathrm{S}_{I,m}(E)$ & $(6_{2}4_{3})$-family (Figure~\ref{fig:S-0-conf})
\\
\cline{2-3}
& $\mathrm{S}^{\rm cs}_{m-1}(E)= \mathrm{S}_{m-1}(E)$ & $\CC^{*}$
\\
\hline
\multirow{3}{*}{$1$} & $\mathrm{Q}^{\rm cs}(E)= \iota(\mathrm{U}(E))$
& $\CC\PP^{1}$ w. 4 double points
\\
\cline{2-3}
& $\mathrm{R}^{\rm cs}(E) = \bigsqcup_{\vert I\vert = 1,2}\mathrm{R}_{I}(E)$ & $(4_{3}6_{2})$-family (Figure~\ref{fig:R-1-conf})
\\
\cline{2-3}
& $\mathrm{S}^{\rm cs}_{m}(E)=\mathrm{S}_{m}(E)$ & $D_{4}$-config. (Figure~\ref{fig:D4-conf})
\\
\hline
\multirow{3}{*}{$2$} & $\mathrm{Q}^{\rm cs}(E)= \iota\bigl(\mathrm{B}_{\varnothing }(E)\backslash\mathrm{B}_{\{1,2,3,4\},0}(E)\bigr)$ & $\mathrm{point}$
\\
\cline{2-3}
& $\mathrm{R}^{\rm cs}(E) = \bigsqcup_{\vert I\vert = 0,1}\mathrm{R}_{I}(E)$ & $(1-4)$-family (Figure~\ref{fig:R-2-conf})
\\
\cline{2-3}
& $\mathrm{S}^{\rm cs}_{m-1}(E)=\mathrm{S}_{m-1}(E)$ & $\mathrm{point}$
\\
\hline
\multirow{2}{*}{$3$} & $\mathrm{Q}^{\rm cs}(E)=\varnothing $ & --
\\
\cline{2-3}
& $\mathrm{R}^{\rm cs}(E) = \mathrm{R}_{\varnothing }(E)$ & $\CC\PP^{1}$
\\
\hline
\end{tabu}}
\end{table}

For each degree parity, there is a bijective correspondence between nontrivial basic building blocks and subsets $I\subset\{1,2,3,4\}$. We conclude with a combinatorial rule to organize their families according to their biholomorphism type.

When the parities of $\vert I\vert$ and $d$ differ, basic building blocks are biholomorphic to $\CC\PP^{1}$ and in~bijective correspondence with the sets of exterior semi-stability walls. 
Concretely, the collections~are
\begin{gather}\label{eq:horizontal}
\begin{cases}
\{\cK_{I}\colon \vert I\vert = 1\} \text{ and } \{\cL_{I}\colon \vert I\vert = 3\}
& \text{if\quad $d$ is even,}
\\
\cK_{\varnothing },\, \{\cK_{I}\colon \vert I\vert = 2\} \text{ and } \cL_{\{1,2,3,4\}} & \text{if\quad $d$ is odd.}
\end{cases}
\end{gather}
In turn, a nontrivial basic building block is biholomorphic to $\CC$ if the parities of $\vert I\vert$ and $d$ coincide. It is convenient to group these in pairs corresponding to partitions $I\sqcup I' = \{1,2,3,4\}$ of even and odd type, as those pairs are in bijective correspondence with interior semi-stability walls. This will be the basis to the classification of interior wall-crossing. The resulting pairs are
\begin{gather}\label{eq:vertical}
\begin{cases}
\bigl\{\cL_{\{i,j\}},\cL_{\{k,l\}}\bigr\}\text{ and } \bigl\{\cK_{\varnothing },\cL_{\{1,2,3,4\}}\bigr\} & \text{if\quad $d$ is even},
\\
\bigl\{\cK_{\{i\}}, \cL_{\{j,k,l\}}\bigr\} & \text{if\quad $d$ is odd}.
\end{cases}
\end{gather}

We will also group $\cN_{\{1,2,3,4\}} := \iota(\mathrm{B}_{\{1,2,3,4\},m-1}(E))/\PP(\Aut(E))$ when $m_{1} = m_{2}=m$, and $\cN_{\varnothing } := \iota(\mathrm{U}(E))/\PP(\Aut(E))$ when $m_{1}-m_{2} = 2$.
With the exception of the previous pair of orbits occurring over different bundles, whose topology will be considered subsequently, all other pairs of the form $\{\cN_{I},\cN_{I^{\mathsf{c}}}\}$ are orbits in the same ambient space having intersecting closures. When~$E$ evenly-split, these are the orbits determining double points in the orbit space of $\mathrm{Q}^{\rm cs}(E)$, locally compactifying the corresponding puncture sphere model of $\cN_{\qB_{\sbul}}$.

\begin{figure}[!ht]
\centering
\includegraphics[width=4.1in]{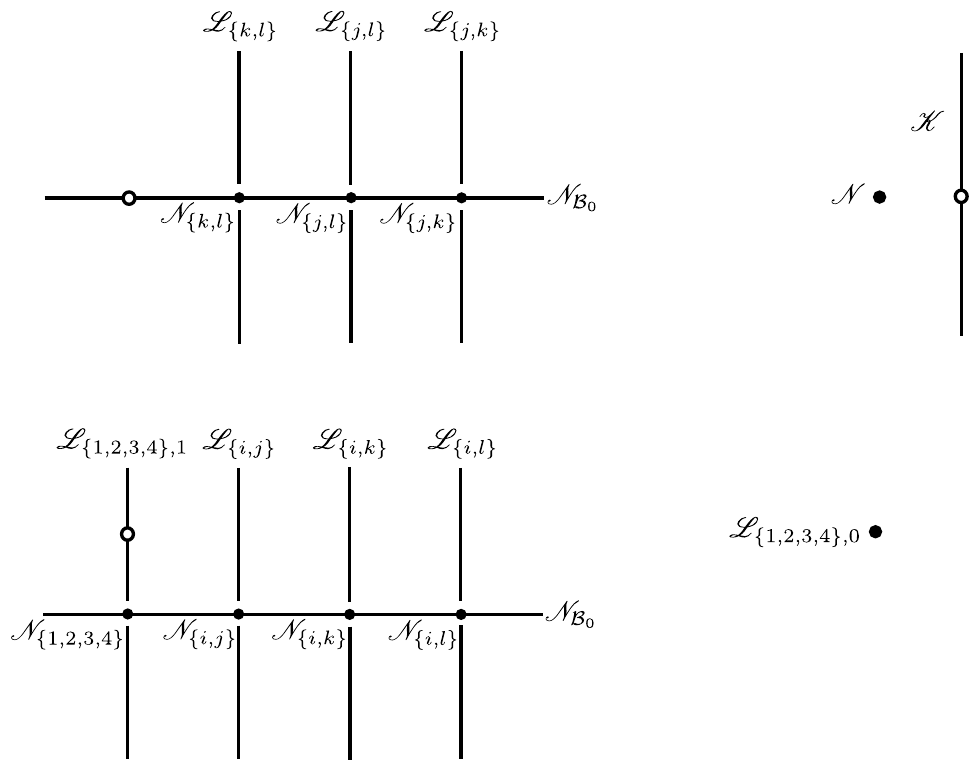}
\put(-11,190){\makebox(0,0)[lb]{\tiny$\varnothing$}}
\put(-36,168){\makebox(0,0)[lb]{\tiny$\varnothing$}}
\put(-186,115){\makebox(0,0)[lb]{\small$a$}}
\put(-186,-11){\makebox(0,0)[lb]{\small$b$}}
\caption{Nilpotent cone assembly kits for even degree interior chambers $\qC_{\qI}$ of type $A$. Left (resp. right) columns correspond to $E\cong\calO(m)\oplus\calO(m)$ (resp.~$E\cong\calO(m+1)\oplus\calO(m-1)$). $(a)$~$\qI = \{\varnothing , \{k,l\}, \{j,l\}, \{j,k\}\}$. $(b)$ $\qI = \{\{1,2,3,4\}, \{i,j\}, \{i,k\}, \{i,l\}\}$.}\label{fig:even-int-A}
\end{figure}

\begin{Definition}\label{def:a-kit}
The \textit{nilpotent cone assembly kit} of an open chamber $\qC\subset \qW_{\sbul}$ of given parity is the collection of stable components in all conditionally stable orbit spaces associated to all strata in the two admissible nilpotent loci.
\end{Definition}

\begin{Proposition}[open chamber stable orbit loci]\label{prop:C-stable}\quad
\begin{enumerate}\itemsep=0pt
\item[$(i)$] The point $\cN_{I}$ is stable with respect to an interior chamber $\qC_{\qI}$ if and only if $I\in\qI$.
\item[$(ii)$] A $I$-basic building block is stable with respect to an interior chamber $\qC_{\qI}$ if and only if $I\in \qI$, 
and stable with respect to an exterior chamber $\qC_{I'}$ if and only if $I\in\qI(I')$ and $I\neq I'$.
\end{enumerate}
\end{Proposition}

\begin{proof}
Immediate from the definition of stability and Corollaries~\ref{cor:ex-chambers}, \ref{cor:in-chambers} and \ref{cor:corresp}.
\end{proof}

\begin{figure}[!ht]
\centering
\includegraphics[width=4.1in]{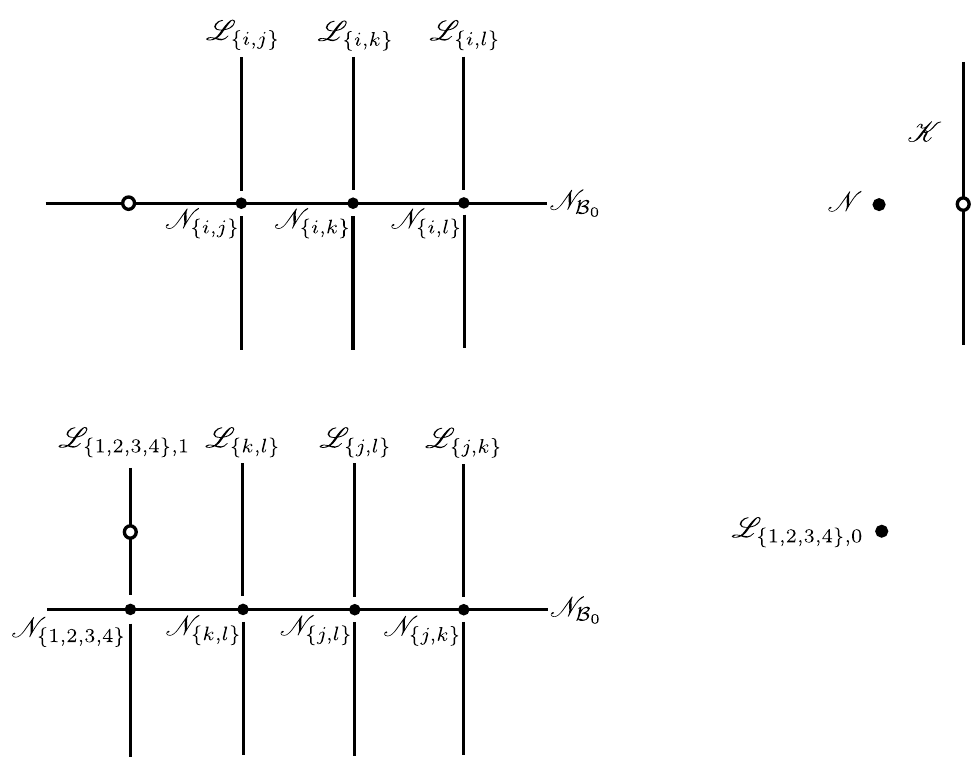}
\put(-12,186){\makebox(0,0)[lb]{\tiny$\varnothing$}}
\put(-37,165){\makebox(0,0)[lb]{\tiny$\varnothing$}}
\put(-186,115){\makebox(0,0)[lb]{\small$a$}}
\put(-186,-11){\makebox(0,0)[lb]{\small$b$}}
\caption{Nilpotent cone assembly kits for even degree interior chambers $\qC_{\qI}$ of type $B$. Left (resp. right) columns correspond to $E\cong\calO(m)\oplus\calO(m)$ (resp.~$E\cong\calO(m+1)\oplus\calO(m-1)$). $(a)$~$\qI = \{\varnothing , \{i,j\},\{i,k\},\{i,l\}\}$. $(b)$ $\qI = \{\{1,2,3,4\}, \{k,l\}, \{j,l\}, \{j,k\}\}$.}\label{fig:even-int-B}
\end{figure}

\begin{figure}[!ht]
\centering
\includegraphics[width=4.6in]{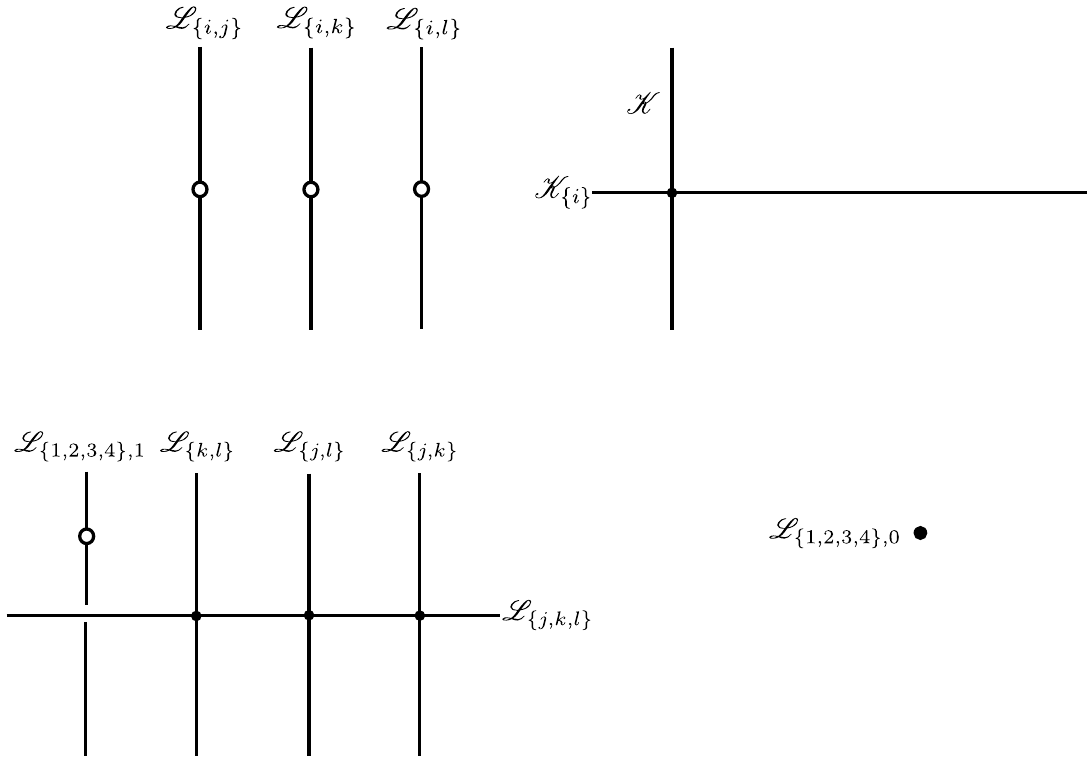}
\put(-134,194){\makebox(0,0)[lb]{\tiny$\varnothing$}}
\put(-184,115){\makebox(0,0)[lb]{\small$a$}}
\put(-184,-11){\makebox(0,0)[lb]{\small$b$}}
\caption{Nilpotent cone assembly kits for even degree exterior chambers $\qC_{I}$. Left (resp.~right) columns correspond to $E\cong\calO(m)\oplus\calO(m)$ (resp.~$E\cong\calO(m+1)\oplus\calO(m-1)$). $(a)$ $I=\{i\}$. $(b)$ $I=\{j,k,l\}$.}\label{fig:even-ext}
\end{figure}

\begin{figure}[!ht]
\centering
\includegraphics[width=2.6in]{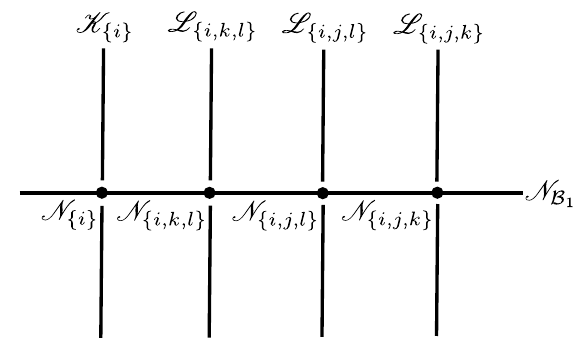}\quad
\includegraphics[width=2.6in]{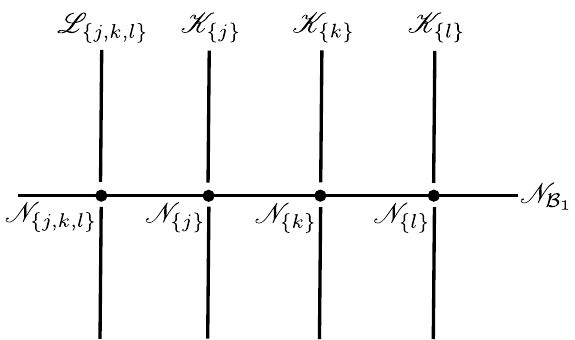}
\put(-303,-13){\makebox(0,0)[lb]{\small$a$}}
\put(-104,-13){\makebox(0,0)[lb]{\small$b$}}
\caption{Nilpotent cone assembly kits for odd degree interior chambers of type $A$, all in $E\cong\calO(m+1)\oplus\calO(m)$. $(a)$ $\qI = \{\{i\},\{i,k,l\},\{i,j,l\},\{i,j,k\}\}$. $(b)$ $\qI = \{\{j,k,l\},\{j\},\{k\},\{l\}\}$.}\label{fig:odd-Int-A}
\end{figure}

\begin{figure}[!ht]
\centering
\includegraphics[width=2.6in]{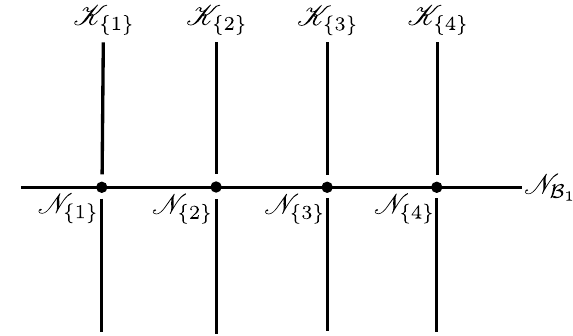}\quad
\includegraphics[width=2.6in]{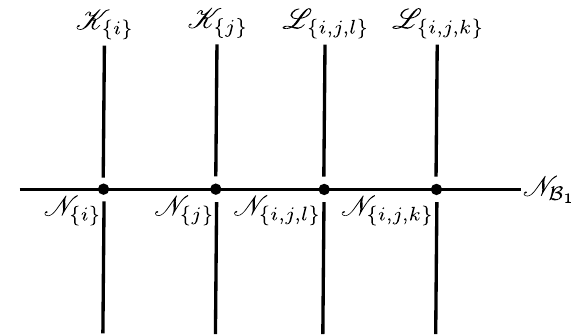}
\put(-303,-12){\makebox(0,0)[lb]{\small$a$}}
\put(-104,-12){\makebox(0,0)[lb]{\small$b$}}
\\[3ex]
\includegraphics[width=2.6in]{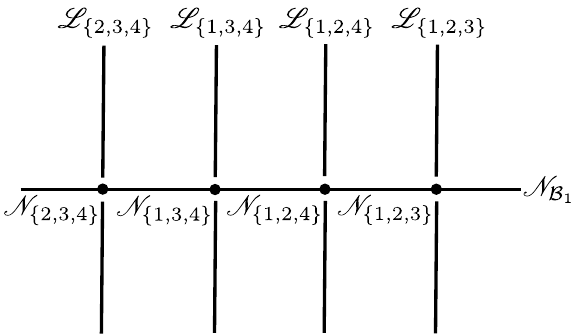}
\put(-103,-12){\makebox(0,0)[lb]{\small$c$}}
\caption{Nilpotent cone assembly kits for odd degree interior chambers of type $B$, all in $E\cong\calO(m+1)\oplus\calO(m)$. $(a)$ $\qI = \{\{1\},\{2\},\{3\},\{4\}\}$. $(b)$ $\qI = \{\{i\},\{j\},\{i,j,l\},\{i,j,k\}\}$. $(c)$ $\qI = \{\{2,3,4\},\{1,3,4\},\{1,2,4\},\{1,2,3\}\}$.}\label{fig:odd-Int-B}
\end{figure}

\begin{figure}[!ht]
\centering
\includegraphics[width=4.2in]{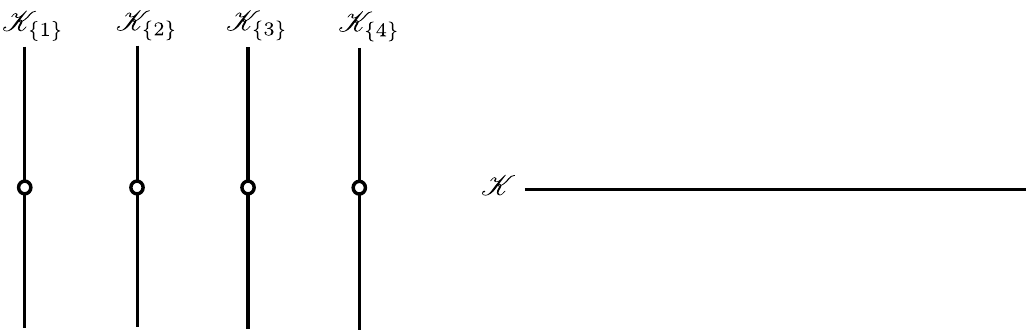}
\put(-154,37){\makebox(0,0)[lb]{\tiny$\varnothing$}}
\put(-170,-12){\makebox(0,0)[lb]{\small$a$}}
\\[3ex]
\includegraphics[width=2.56in]{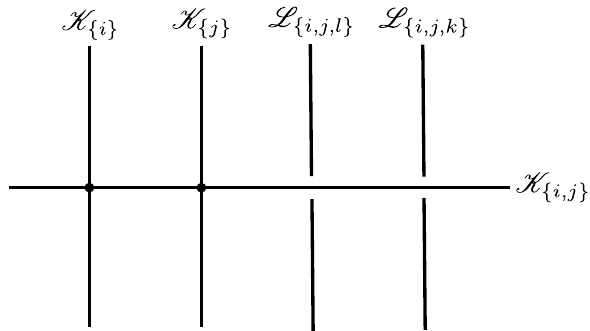}\quad
\includegraphics[width=2.6in]{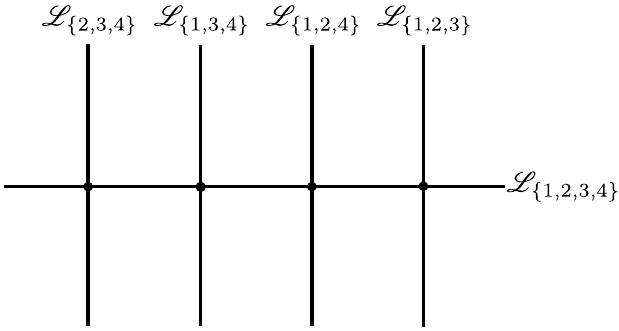}
\put(-308,-12){\makebox(0,0)[lb]{\small$b$}}
\put(-111,-12){\makebox(0,0)[lb]{\small$c$}}
\caption{Nilpotent cone assembly kits for odd degree exterior chambers $\qC_{I}$. Left (resp.~right) columns correspond to $E\cong\calO(m+1)\oplus\calO(m)$ (resp.~$E\cong\calO(m+2)\oplus\calO(m-1)$). 1) $I=\varnothing $. 2) $I=\{i,j\}$. 3) $I=\{1,2,3,4\}$.}\label{fig:odd-ext}
\end{figure}

Figures \ref{fig:even-int-A}--\ref{fig:odd-ext}, built from Definition~\ref{def:a-kit} and Propositions~\ref{prop:line-orbit} and \ref{prop:C-stable}, compile the geometry and combinatorics of all nilpotent cone assembly kits
in both parities, with central spheres represented as horizontal lines. In particular, when $d$ is odd, there is a single Harder--Narasimhan stratum for all open chambers (corresponding to $E\cong\calO(m+1)\oplus\calO(m)$), except for the exterior chamber $\qC_{\varnothing }$.

\section{Proof of main results and further properties}\label{sec:proof-m}

\subsection{Proof of Theorem~\ref{theo:main-1}}
We will combine all previous results to construct all Harder--Narasimhan strata $\cM_{\qC}(E)$ of a~given open chamber $\qC$ as $\PP(\Aut(E))$-orbit spaces. The missing information to do so is contained in~Proposition~\ref{prop:orbit-bulk}, where we relate the geometry of bulk orbit spaces and basic building blocks.\looseness=1

{\samepage\begin{Proposition}[orbit space models in $\mathrm{QPH}(E)$]\label{prop:orbit-bulk}\quad
\begin{enumerate}\itemsep=0pt
\item[$(i)$] The group $\PP(\Aut(E))$ acts freely on $\mathrm{QPH}(E)\backslash\iota(\mathrm{U}(E)^{\mathsf{c}})$.
\item[$(ii)$] For $E$ evenly-split, $\PP(\Aut(E))$ acts properly on
$\mathrm{QPH}_{\CC^{*}}(E)$, turning it into a principal $\PP(\Aut(E))$-bundle over $(\Sigma_{D}\backslash\{w_{1}\})\times \CC^{*}$ when $E\cong\calO(m)\oplus\calO(m)$ and $\Sigma_{D}\times \CC^{*}$ when $E\cong\calO(m+1)\oplus\calO(m)$.
\item[$(iii)$] If $m_{1}-m_{2}=2$, the map $\mathrm{par}\times \det$ determines the isomorphism
\[
\mathrm{QPH}_{\CC^{*}}(E)\sqcup \mathrm{S}_{0}(E)\cong\mathrm{P}_{m-1}(E)\times H^{0}\bigl(\CC\PP^{1},K_{\CC\PP^{1}}^{2}(D)\bigr),
\]
inducing the biholomorphism
\[
\bigl\{\mathrm{QPH}_{\CC^{*}}(E)\sqcup \mathrm{S}_{0}(E)\bigr\}/\PP(\Aut(E))\cong H^{0}\bigl(\CC\PP^{1},K_{\CC\PP^{1}}^{2}(D)\bigr).
\]
\end{enumerate}
\end{Proposition}}

\begin{proof}
$(i)$ Propositions \ref{prop:U} and \ref{prop:X-Y} establish the claim over $\mathrm{Q}(E)$. The triviality of $\PP(\Aut(E))$-stabilizers for points in $\mathrm{QPH}(E)\backslash\mathrm{Q}(E)$ follows from Proposition~\ref{prop:QPH-par}$(i)$. Namely, the claim is trivial over $\mathrm{par}^{-1}(\mathrm{U}'(E))$, and on $\mathrm{QPH}_{\CC^{*}}(E)\cap\mathrm{par}^{-1}(\mathrm{X}(E))$ it follows from the explicit form of the residual $\CC^{*}$-automorphism action on the canonical forms \eqref{eq:HX-0}--\eqref{eq:HX-1}. Corollary \ref{cor:strata}, Proposition~\ref{prop:S_{j}(E)}, and the canonical form \eqref{eq:Phi-1} for elements in $\mathrm{R}(E)$ when $m_{1}-m_{2}=1$ imply the claim over $\mathrm{QPH}_{0}(E)\backslash\mathrm{Q}(E)$.

$(ii)$ Both claims follow in analogy to $(i)$, after restriction to $\mathrm{QPH}_{\CC^{*}}(E)$. For any $q\in H^{0}(\CC\PP^{1},K^{2}_{\CC\PP^{1}}(D))\backslash\{0\}$, the orbits in $\det^{-1}(q)\cap\mathrm{par}^{-1}(\mathrm{X}(E))$ are fixed by $\tau$, and the limit properties of $\mathrm{X}(E)$ with respect to $\mathrm{U}(E)$ described in Proposition~\ref{prop:X-Y} imply that the orbit space of $\det^{-1}(q)$ is biholomorphic to the punctured elliptic curve $\Sigma_{D}\backslash\{w_{1}\}$ when $m_{1}=m_{2}=m$, and to $\Sigma_{D}$ when $m_{1}-m_{2}=1$.

$(iii)$ In analogy to $(ii)$, the divisor of any parabolic Higgs field $\Phi$ is necessarily trivial for both strata in $\mathrm{QPH}_{\CC^{*}}(E)\sqcup \mathrm{S}_{0}(E)$ when $m_{1} - m_{2} = 2$, and every point is uniquely determined by $\Phi$, which admits a unique decomposition of the form
\[
\Phi = \Phi_{0} + \Phi_{1},\qquad
\det(\Phi_{0})=0,\qquad
\Phi_{1}\in H^{0}\bigl(\CC\PP^{1},\End(E)\otimes K_{\CC\PP^{1}}\bigr).
\]
Moreover, the map
\[
\Phi_{1}\mapsto \det(\Phi_{0}+ \Phi_{1})
\]
is a linear isomorphism to $H^{0}(\CC\PP^{1},K^{2}_{\CC\PP^{1}}(D))$ for any fixed $\Phi_{0}$, from which the first claim follows. The second claim follows from Proposition~\ref{prop:CS}.
\end{proof}

\begin{figure}[!ht]
\centering
\includegraphics[width=3.6in]{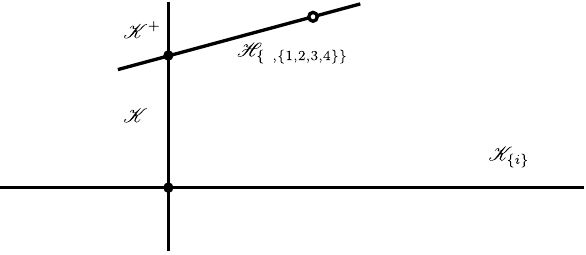}
\put(-197.5,56){\makebox(0,0)[lb]{\tiny$\varnothing$}}
\put(-144,85.5){\makebox(0,0)[lb]{\tiny$\varnothing$}}
\put(-198,93){\makebox(0,0)[lb]{\tiny$\varnothing$}}
\caption{Maximal Harder--Narasimhan strata in the splitting $E\cong\calO(m+1)\oplus\calO(m-1)$ corresponding to the exterior chambers $\qC_{\{i\}}$, consisting of the image of a Hitchin section and two nilpotent cone irreducible components.}\label{fig:HN-stratum}
\end{figure}

\begin{proof}[Proof of Theorem~\ref{theo:main-1}]
The case $m_{1}-m_{2}=3$ follows from Propositions~\ref{prop:CS}, \ref{prop:line-orbit} and \ref{prop:C-stable}, since $\mathrm{QPH}_{\CC^{*}}(E)=\varnothing $ (Lemma~\ref{lemma:cases}), and is trivial for all open chambers other than $\qC_{\varnothing }$, in which case the Harder--Narasimhan stratum is a projective line (Figures~\ref{fig:odd-Int-A}--\ref{fig:odd-ext}).

The case $m_{1}-m_{2}=2$ can be concluded from Propositions~\ref{prop:C-stable} and \ref{prop:orbit-bulk}$(iii)$ (see Figures~\ref{fig:even-int-A}--\ref{fig:even-ext}; Figure~\ref{fig:HN-stratum} illustrates the topology of the Harder--Narasimhan stratum when its nilpotent cone component is a nodal rational curve).

When $E$ is evenly-split, both claims would follow as a consequence of Proposition~\ref{prop:orbit-bulk}$(i)$ and $(ii)$ once the orbit limits of $\mathrm{QPH}_{\CC^{*}}(E)$ as $\det(\Phi)\to 0$ are established. Proposition \ref{prop:C-stable} ensures that for every possible choice of open chamber, this limit would be a rational curve, possibly noncompact and disconnected, with at most nodal points. The possibilities for both degree parities are illustrated in detail in the left column of Figures~\ref{fig:even-int-A}--\ref{fig:odd-ext}. By construction, the resulting orbit spaces are smooth surfaces.
\end{proof}

\subsection{Proof of Theorems~\ref{theo:main-2} and \ref{theo:main-3}}

In order to understand the geometry of gluing of strata, we will present a series of results on orbit space limits associated to \textit{jumping families} of even (resp.~odd) degree bundles, i.e., 1-parameter families $\cF = \{E^{t}\colon t\in\CC\}$ for which
\[
E^{0}\cong \calO(m+1)\oplus\calO(m-1)\qquad \big(\text{resp. } E^{0}\cong \calO(m+2)\oplus\calO(m-1)\big)
\]
and $E^{t}$ is evenly-split if $t\in\CC^{*}$. A holomorphic family of cocycles $\{g_{01}^{t}\colon t\in\CC\}$ with respect to the standard affine cover $\{\cU_{0}= \CC\PP^{1}\backslash\{\infty\}$, $\cU_{1}= \CC\PP^{1}\backslash\{0\}\}$ of $\CC\PP^{1}$, with coordinates $z\in\cU_{0}$, $\zeta\in\cU_{1}$ related on the intersection $\cU_{01}$ as $\zeta = 1/z$, can be given as follows
\[
g^{t}_{01}(z) = \begin{cases}
\begin{pmatrix}
z^{m+1} & 0\\
t & z^{m-1}
\end{pmatrix}\!, & d\quad \text{even},
\\[3ex]
\begin{pmatrix}
z^{m+2} & 0\\
t & z^{m-1}
\end{pmatrix}\!, & d\quad \text{odd}.
\end{cases}
\]
That the corresponding vector bundles $E^{t}$ in each degree parity are evenly-split when $t\in\CC^{*}$, is elucidated in the matrix factorizations of the form $g^{t}_{01} = g^{t}_{0}g_{01}(g^{t}_{1})^{-1}$,
\begin{gather}\label{eq:cocycle}
 \begin{cases}
\begin{pmatrix}
z^{m+1} & 0\\
t & z^{m-1}
\end{pmatrix} =
\begin{pmatrix}
1 & -z\\
0 & -t
\end{pmatrix}
\begin{pmatrix}
z^{m} & 0\\
0& z^{m}
\end{pmatrix}
\begin{pmatrix}
z^{-1} & -1\\
-t & 0
\end{pmatrix}^{-1}, & d \quad \text{even},
\\[3ex]
\begin{pmatrix}
z^{m+2} & 0\\
t & z^{m-1}
\end{pmatrix} =
\begin{pmatrix}
1 & -z^{2}\\
0 & -t
\end{pmatrix}
\begin{pmatrix}
z^{m+1} & 0\\
0& z^{m}
\end{pmatrix}
\begin{pmatrix}
z^{-1}& -1\\
-t& 0
\end{pmatrix}^{-1}, & d \text{ odd},
\end{cases}
\end{gather}
 inducing the cocycle equivalences $g^{t}_{01}\cong g_{01}$ with the standard evenly-split cocycle. The implications of the existence of jumping families on the limiting behavior of the relevant orbit spaces of quasi-parabolic structures are gathered in the next lemma.

\begin{Lemma}\label{lemma:limit}
Let $\{E^{t}\colon t\in\CC\}$ be a jumping family of vector bundles of arbitrary degree. Then
\[
\lim\limits_{t\rightarrow 0} \mathrm{Y}\big(E^{t}\big) \subset \overline{\mathrm{Y}\big(E^{0}\big)}.
\]
If $d$ is even, then
\[
\lim\limits_{t\rightarrow 0} \mathrm{B}_{\{1,2,3,4\},m-1}\big(E^{t}\big) = \mathrm{X}\big(E^{0}\big).
\]
\end{Lemma}

\begin{proof}
It follows from the definition of $\mathrm{Y}(E)$ in Proposition~\ref{prop:X-Y}, when $E$ is not evenly-split, that
\[
\overline{\mathrm{Y}(E)} =
 \begin{cases}
\displaystyle\bigsqcup_{i=1}^{4} \{\{E_{1}\vert_{z_{i}}\} \} \times \bigg\{\prod_{j\neq i}\PP(E\vert_{z_{i}}) \bigg\},& d\quad \text{even},\\
\displaystyle\mathrm{QP}(E), & d\quad \text{odd}.
\end{cases}
\]
The first claim follows from the cocycle factorizations \eqref{eq:cocycle}, after acting with the local change of trivializations $g^{t}_{0}$ and $g^{t}_{0}$ on the elements of $\mathrm{Y}(E)$, $E$ evenly-split. If $d$ is even, the second statement is verified in the same way after normalizing elements in $\mathrm{B}_{\{1,2,3,4\},1}(E^{t})$, $t\in\CC^{*}$, in the form $(z_{1},0,1,\infty)$, under the identification $\cC(E^{t})\cong \mathrm{Conf}_{4}(\mathrm{L}_{0}(E^{t}))$ and the isomorphism $\mathrm{L}_{0}(E^{t})\cong\CC\PP^{1}$.
\end{proof}

We will introduce additional notation for special components of orbit spaces in $\mathrm{QPH}_{\CC^{*}}(E)$. Recall that for $m_{1}-m_{2}\leq 2$, the connected components of the sets $\mathrm{X}(E)$ are parametrized by $I\sqcup I' = \{1,2,3,4\}$ of the corresponding degree parity. Concretely, these are
\[
 \begin{cases}
\{i,j\}\sqcup \{k,l\},& m_{1} = m_{2} = m,\\
\{i\}\sqcup\{j,k,l\}, & m_{1} - m_{2} = 1,\\
\varnothing \sqcup \{1,2,3,4\}, & m_{1} - m_{2} = 2.
\end{cases}
\]
For every partition $I\sqcup I' = \{1,2,3,4\}$, we will denote its corresponding orbit space in $\mathrm{QPH}_{\CC^{*}}(E)$ by $\mathrm{H}_{\{I,I'\}}$. In all cases, the isomorphism $\mathrm{H}_{\{I,I'\}}\cong \CC^{*}$ defined by the determinant map $\Phi\mapsto\det(\Phi)$ readily follows from Proposition~\ref{prop:orbit-bulk}. Moreover, in the case when $m_{1}-m_{2} = 2$, we have that $\mathrm{H}_{\{\varnothing ,\{1,2,3,4\}\}} = \mathrm{QPH}_{\CC^{*}}(E)/\PP(\Aut(E))$. The next corollary is a consequence of Lemma~\ref{lemma:limit} (see Figure~\ref{fig:HN-stratum}).

\begin{Corollary}\label{cor:glue-H}
Every orbit in $\mathrm{H}_{\{\varnothing ,\{1,2,3,4\}\}}$ is the limit of a holomorphic family of orbits in $\mathrm{QPH}_{\CC^{*}}(E)/\PP(\Aut(E))$, with $E\cong \calO(m)\oplus\calO(m)$.
\end{Corollary}

\begin{Proposition}[limits of $\PP(\Aut(E))$-orbit families]\label{prop:glue}\quad
\begin{enumerate}\itemsep=0pt
\item[$(i)$] If $d$ is even, there is a jumping family such that\vspace{-1ex}
\begin{gather*}
\lim\limits_{t\to 0} \cN^{t}_{\{1,2,3,4\}} = \cN^{0}_{\varnothing },
\end{gather*}
while\vspace{-1ex}
\begin{gather*}
\lim\limits_{t\to 0} \bigsqcup \cL^{t}_{\{j,k,l\}} = \bigsqcup\cK^{0}_{\{i\}},
\end{gather*}
and moreover, there is a family of orbits $\{(E^{t}_{\qp},\Phi^{t})\}\in\cL_{\{1,2,3,4\},1}$ such that\vspace{-1ex}
\begin{gather*}
\lim\limits_{t\to 0}\bigg\{\bigg(E^{t}_{\qp},\frac{1}{t}\Phi^{t}\bigg)\bigg\} = \cL^{0}_{\{1,2,3,4\},0}.
\end{gather*}
In particular, $\cL_{\{1,2,3,4\},1}\cong\CC^{*}$ is compactified by the orbits $\cL_{\{1,2,3,4\},0}$ and $\cN_{\{1,2,3,4\}}$.
\item[$(ii)$] If $d$ is odd, then there is a jumping family for which\vspace{-1ex}
\begin{gather*}
\lim\limits_{t\to 0} \cL^{t}_{\{1,2,3,4\}} = \cK^{0}_{\varnothing }.
\end{gather*}
\end{enumerate}
\end{Proposition}

\begin{proof}
Let $E = \calO(m)\oplus\calO(m)$ and consider any nonzero parabolic Higgs field of the form\vspace{-1ex}
\begin{gather*}
\Phi = \begin{pmatrix}
u_{1}v_{1} & - u_{1}^{2}\\
v_{1}^{2} & -u_{1}v_{1}
\end{pmatrix}\!,\qquad(u_{1})\neq (v_{1})
\end{gather*}
whose kernel line is an arbitrary element of $\mathrm{L}_{1}(E)$, and consider the family of quasi-parabolic Higgs bundles $\{(E^{t},\Phi^{t})\colon t\in\CC^{*}\}$ that is induced in terms of the local form\vspace{-1ex}
\begin{gather*}
\Phi^{t}\vert_{\cU_{0}} := t g^{t}_{0}(\Phi\vert_{\cU_{0}})\big(g^{t}_{0}\big)^{-1}.
\end{gather*}
Since\vspace{-1ex}
\begin{gather*}
\lim\limits_{t\to 0} \ker\big(\Phi^{t}\big) = E^{0}_{1},
\end{gather*}
it follows from Lemma~\ref{lemma:limit} that the corresponding $\PP(\Aut(E^{0}))$-orbit equals $\cL^{0}_{\{1,2,3,4\},0}$. The remaining claims follow analogously from Lemma~\ref{lemma:limit} and Proposition~\ref{prop:line-orbit}, given that\vspace{-1ex}
\begin{gather*}
\mathrm{par}\bigg( \bigsqcup_{\vert I \vert = 3 + m_{1} - m_{2}} \mathrm{S}_{I,m}(E)\bigg) = \mathrm{Y}(E)
\end{gather*}
when $E$ is evenly-split, and otherwise\vspace{-1ex}
\begin{gather*}
\mathrm{par}\bigg( \bigsqcup_{\vert I \vert = 3 - m_{1} + m_{2}} \mathrm{R}_{I}(E)\bigg) = \mathrm{Y}(E).
\tag*{\qed}
\end{gather*}
\renewcommand{\qed}{}
\end{proof}

It follows from Propositions~\ref{prop:orbit-bulk} and \ref{prop:glue} that for each choice of degree parity, the non-Hausdorff phenomena on conditionally stable loci concentrates on the basic building block data \eqref{eq:horizontal} and~\eqref{eq:vertical}. This is stated in the following corollary.

\begin{Definition}
The \textit{central sphere} $\cS_{\qI}$ for an interior open chamber $\qC_{\qI}$ is the associated moduli space of stable parabolic bundles $\cN_{\qC_{\qI}}$. The \textit{central sphere} $\cS_{I}$ for an exterior open chamber $\qC_{I}$ is its corresponding stable basic building block \eqref{eq:horizontal}.
\end{Definition}

\begin{Corollary}\label{cor:n-H-blocks}
The union of all exterior central spheres $\cS_{I}$ of a given degree parity
determines a line of $8$-tuple points. The pair \eqref{eq:vertical} of a given partition $I\sqcup I' =\{1,2,3,4\}$ determines a punctured line of double points compactified by $\{\cN_{I},\cN_{I'}\}$. The union of all conditionally stable loci of given degree parity conforms a $D_{4}$-configuration whose nodes are $10$-tuple points, its central punctured sphere is formed by $9$-tuple points, and its $4$ complementary components are formed by double points.
\end{Corollary}

{\sloppy\begin{proof}[Proof of Theorem~\ref{theo:main-2}]
It follows from Proposition~\ref{prop:orbit-bulk} together with Corollary~\ref{cor:glue-H} that when~$d$ is even, $\mathrm{H}_{\{\varnothing ,\{1,2,3,4\}\}}$ glues into the bulk $\mathrm{QPH}_{\CC^{*}}(E)$ for $m_{1}=m_{2}=m$ as a branch locus component under $\tau$ over the ramification point $w_{1}$, in a compatible way with respect to the map $\det$. This determines an isomorphism for both degree parities between $\det^{-1}(H^{0}(\CC\PP^{1},K^{2}_{\CC\PP^{1}}(D))\allowbreak\backslash\{0\})$ and $\Sigma_{D}\times\CC^{*}$ which is independent of the choice of open chamber. Similarly, it follows from Proposition~\ref{prop:glue} that for all open chambers, the different components in the nilpotent cone assembly kit (whose explicit classification is compiled in Figures~\ref{fig:even-int-A}--\ref{fig:odd-ext}) glue into a $D_{4}$-configuration (Figure~\ref{fig:D4-conf}), on which $\det$ extends to 0 as a consequence of Proposition~\ref{prop:orbit-bulk}. The resulting explicit geometry of the map $\det$ implies the isomorphism $\cM_{\qC}\cong \cM_{\mathrm{toy}}$ for all degrees and all choices of open chambers.
\end{proof}}

{\sloppy\begin{proof}[Proof of Theorem~\ref{theo:main-3}]
Since wall-crossing is concentrated along the locus $\det^{-1}(0)$ in each~$\cM_{\qC}$, all that remains to be done after having established Theorems~\ref{theo:main-1} and~\ref{theo:main-2} is to understand the resulting transformations in $\det^{-1}(0)$ when either an interior or exterior semi-stability wall is crossed. The claim for both possibilities follows if we combine Corollaries~\ref{cor:corresp} and \ref{cor:reflections} with Corollary~\ref{cor:n-H-blocks}, as the latter ensures that all nilpotent cone components that can be pairwise exchanged according to the combinatorial rules dictated by the former share the same topology in their respective geometric model. 
\end{proof}}

\begin{Remark}\label{rem:Blaavand}
In \cite[Section 6.7]{Bla15}, Blaavand considers parabolic Higgs bundles of parabolic degree zero on vector bundles of degree $-1$, subject to the constraints $\alpha_{11} = \alpha_{21} = \alpha_{31} = 0$ and $\alpha_{12}+\alpha_{22}+\alpha_{32}+\alpha_{41}+\alpha_{42} =1$. Those constraints confine the parabolic weights $\pmb{\beta}$ to lie in the open chamber $\qC_{\{1,2,3,4\}}$. Consequently, the moduli spaces that he considers correspond to the geometric model $\cM_{\qC_{\{1,2,3,4\}}}$.
\end{Remark}

\begin{Remark}
The isomorphisms for the different open chamber geometric models are also implicit in Corollary~\ref{cor:reflections}.
In fact, all the constructed geometric models are obviously isomorphic as complex surfaces, independently of the choice of degree parity. One standard strategy to construct isomorphisms between moduli spaces associated to parabolic structures is known in the literature under the name of \textit{elementary transformations} \cite{LS15,LSS13}. We have chosen not to invoque that notion in this work, as our method of construction of geometric models is not compatible with it: the groups and associated orbits are fundamentally different for each choice of degree parity.
\end{Remark}

\subsection[The Hitchin fibration, \protect{$C\textasciicircum{}{*}$}-action, and Hitchin sections]
{The Hitchin fibration, $\boldsymbol{\CC^{*}}$-action, and Hitchin sections}

Each geometric model $\cM_{\qC}$ can be endowed with a product of transversal projections
\[
\pi_{1}\times \pi_{2}\colon\ \cM_{\qC}\rightarrow \CC\PP^{1}\times H^{0}\bigl(\CC\PP^{1},K^{2}_{\CC\PP^{1}}(D)\bigr),
\]
where $\pi_{1}$ maps onto the central sphere of $\cM_{\qC}$, while $\pi_{2}$ is the standard \textit{Hitchin elliptic fibration}
\[
 \pi_{2}([(E_{*},\Phi)]) = \det(\Phi).
\]
By Proposition~\ref{prop:orbit-bulk}, the nonzero fibers $\pi_{2}^{-1}(q)$, $q\neq 0$ are isomorphic to the elliptic curve $\Sigma_{D}$ of the pair $(\CC\PP^{1},D)$, with an involution induced from \eqref{eq:def-tau} that we will keep denoting by $\tau$. It~follows from Proposition~\ref{prop:orbit-bulk} that the resulting branch loci are parametrized by the connected components of the union of the sets $\mathrm{X}(E)$ of a given degree parity, and consequently, correspond to the four components $\mathrm{H}_{\{I,I'\}}$ parametrized by partitions $I\sqcup I' = \{1,2,3,4\}$ of the given parity.

The zero fiber $\pi_{2}^{-1}(0)$ is called the \textit{nilpotent cone} of the geometric model $\cM_{\qC}$. In general the Zariski open sets $\cM_{\qC}\backslash \pi^{-1}_{2}(0)$ are independent of the choice of open chamber $\qC$.

The definition of $\pi_{1}$ requires careful attention to details and will be presented in several stages. In essence, the restriction $\pi_{1}\times\pi_{2}\vert_{\cM_{\qC}\backslash \pi^{-1}_{2}(0)}$ is defined as the induced branched cover
\begin{gather}\label{eq:branched-cover}
\cM_{\qC}\backslash \pi^{-1}_{2}(0)\xrightarrow{2:1} \CC\PP^{1}\times H^{0}\bigl(\CC\PP^{1},K^{2}_{\CC\PP^{1}}(D)\bigr)\backslash\{0\},
\end{gather}
under identification of the first factor in the image with the central sphere, and the extension to the nilpotent cone is then defined as its collapse to the former. We will make this precise in terms of the global properties of the $\CC^{*}$-action induced by \eqref{eq:def-C*} on any geometric model $\cM_{\qC}$, described in Propositions~\ref{prop:limits} and \ref{prop:Hitchin-s}.

Yet another consequence of Propositions~\ref{prop:QPH-par} and \ref{prop:CS} -- in analogy to the proof of Proposition~\ref{prop:line-orbit} -- is that for any basic building block in a pair \eqref{eq:vertical}, there is exactly one $\PP(\Aut(E))$-orbit contained in $\mathrm{par}^{-1}(\mathrm{X}(E))$. In the special case when $m_{1}-m_{2} = 2$ and $I= \{1,2,3,4\}$, this orbit coincides with $\cL_{\{1,2,3,4\},0}$. Depending on the basic building block they belong to in a given pair, we will denote them as
\[
\cK^{+}_{I}\qquad \big(\text{resp.} \ \cL^{+}_{I}\big).
\]
For any partition $I\sqcup I' = \{1,2,3,4\}$, Corollary~\ref{cor:n-H-blocks}$(ii)$ implies that the pair $\{\cK^{+}_{I},\cL^{+}_{I'}\}$ (or $\{\cL^{+}_{I},\cL^{+}_{I'}\}$ if $\vert I\vert =\vert I'\vert =2$) conforms a double point.

\begin{Proposition}\label{prop:limits}
Any central sphere is pointwise-fixed by the $\CC^{*}$-action.
For any complementary nilpotent cone component in a pair \eqref{eq:vertical}, $\cK^{+}_{I}$ $($resp.~$\cL^{+}_{I})$ is fixed by the $\CC^{*}$-action, while $\cK_{I}\backslash\cK_{I}^{+}$ $($resp.~$\cL_{I}\backslash\cL_{I}^{+})$ is a $\CC^{*}$-principal homogeneous space. Moreover, in the latter case, 
any orbit $\{(E_{\qp},\Phi)\}\in \cK_{I}\backslash\cK_{I}^{+}$ $($resp.~$\cL_{I}\backslash\cL_{I}^{+})$ satisfies
\[
\lim\limits_{c\to \infty}c\cdot\{(E_{\qp},\Phi)\} = \cK^{+}_{I}\qquad \big(\text{resp.~}\cL^{+}_{I}\big),
\]
while for any partition $I\sqcup I' = \{1,2,3,4\}$ of given parity and any orbit $\{(E_{\qp},\Phi)\}\in \mathrm{H}_{\{I,I'\}}$,
\[
\lim\limits_{c\to 0}c\cdot\{(E_{\qp},\Phi)\} = \big\{\cK^{+}_{I},\cL^{+}_{I'}\big\}\qquad \big(\text{or}\ \big\{\cL^{+}_{I},\cL^{+}_{I'}\big\}\ \text{if $\vert I\vert =\vert I'\vert =2$}\big).
\]
\end{Proposition}

\begin{proof}
The first statement is clear when the central sphere is a moduli space of stable parabolic bundles. Otherwise, 
a conditionally stable $\PP(\Aut(E))$-orbit is fixed under the induced $\CC^{*}$-action if and only if for any choice of $c\in \CC^{*}$, the equation
\begin{gather}\label{eq:C*}
(F_{1},F_{2},F_{3},F_{4},c\Phi) = (g\cdot F_{1},g\cdot F_{2},g\cdot F_{3},g\cdot F_{4},\Ad(g)(\Phi))
\end{gather}
can be solved for some $g\in\PP(\Aut(E))$.
We will independently consider the two families of basic building blocks according to their stratum type in $\mathrm{QPH}_{0}(E)$. Under a choice of bundle isomorphism $E\cong \calO(m_{1})\oplus\calO(m_{2})$, every orbit in a basic building block $\cL_{I}$ contains a representative $(F_{1},F_{2},F_{3},F_{4},\Phi)$ whose parabolic Higgs field takes the form
\[
\Phi = \begin{pmatrix}
0 & 0\\
w_{2-m_{1}+m_{2}} & 0
\end{pmatrix}\!,\qquad m_{1}-m_{2} = 0,1,2.
\]
In all cases, the set of these canonical forms is a principal homogeneous space for a subgroup in $\PP(\Aut(E))$ isomorphic to $\CC^{*}$. It follows from Proposition~\ref{prop:line-orbit} and
the definition of $\mathrm{Y}(E)$ for every admissible splitting type, that when an orbit belongs to the central sphere $\cS_{I}$ of an exterior chamber $\qC_{I}$, the equation \eqref{eq:C*} can be solved on its subset of canonical form elements, and the $\CC^{*}$-action leaves the orbit invariant. In the case of a nilpotent cone's complementary component, the same holds for the special orbit $\cL^{+}_{I}$ by the definition of $\mathrm{X}(E)$, while the action of the subgroup stabilizing canonical forms is trivial on quasi-parabolic structures associated to $\cL_{I}\backslash\cL^{+}_{I}$. Hence it coincides with the standard $\CC^{*}$-action, and determines a principal homogeneous space structure on it. The final claim, on the orbit limits under the $\CC^{*}$-action, follows once again from the limiting properties of $\mathrm{X}(E)$ and $\mathrm{Y}(E)$ with respect to the corresponding $\PP(\Aut(E))$ action on the loci $\mathrm{U}(E)$.

On the other hand, in the case of conditionally stable basic building blocks $\cK_{I}$, the parabolic Higgs field of all elements in every given orbit is always in the normal form
\[
\Phi = \begin{pmatrix}
0 & -v_{2+ m_{1} - m_{2}}\\
0 & 0
\end{pmatrix}\!,\qquad m_{1} - m_{2} = 1, 2, 3.
\]
Such parabolic Higgs field's normal forms are $\mathrm{R}(\Aut(E))$-invariant by definition, and the residual group of bundle automorphisms acting on the underlying space $\{E_{\qp}\}$ of quasi-parabolic bundles associated to the orbit is isomorphic to $\CC^{*}$. The same argument as before also verifies the claim in this case.
\end{proof}

\begin{Proposition}\label{prop:Hitchin-s}
For every interior open chamber $\qC_{\qI}$ of given parity, the elliptic fibration $\pi_{2}$ possesses four natural holomorphic sections, extending the connected components of the branch locus of \eqref{eq:branched-cover} to $\pi^{-1}_{2}(0)$, parametrized by the subsets $I\in \qI$,
\[
\sigma_{I}\colon\ H^{0}\bigl(\CC\PP^{1},K^{2}_{\CC\PP^{1}}(D)\bigr)\rightarrow \cM_{\qC_{\qI}}
\]
such that
\[
\sigma_{I}\bigl(H^{0}\bigl(\CC\PP^{1},K^{2}_{\CC\PP^{1}}(D)\bigr)\backslash \{0\}\bigr) = \mathrm{H}_{\{I,I'\}},\qquad \sigma_{I}(0) = \cK^{+}_{I'}\ \big(\text{or}\ \cL^{+}_{I}\big).
\]
For every exterior chamber $\qC_{I''}$, an analogous construction holds in terms of its neighboring interior chamber $\qC_{\qI(I'')}$ of type $B$.
\end{Proposition}

\begin{proof}
Recall that \eqref{eq:HX-0} and~\eqref{eq:HX-1} establish correspondences between nonzero elements
\[
q\in H^{0}\bigl(\CC\PP^{1},K^{2}_{\CC\PP^{1}}(D)\bigr)\backslash\{0\}
\]
and orbits of parabolic Higgs bundles with parabolic Higgs field
\[
\Phi = \begin{pmatrix}
0 & -v_{2+m_{1}-m_{1}}\\
w_{2-m_{1}+m_{2}} & 0
\end{pmatrix}
\]
in terms of factorizations $q = v_{2+m_{1}-m_{1}}w_{2-m_{1}+m_{2}}$, parametrized by partitions $I\sqcup I' = \{1,2,3,4\}$ induced by fixing the divisors $(v_{2+m_{1}-m_{1}})$ and $(w_{2-m_{1}+m_{2}})$. For each degree parity, there are exactly four partitions, and for any such partition, we can assume without loss of generality that \mbox{$I\in \qI$}. By Corollary~\ref{cor:n-H-blocks} and Proposition~\ref{prop:limits}, the induced isomorphism $H^{0}(\CC\PP^{1},K^{2}_{\CC\PP^{1}}(D))\backslash\{0\}\allowbreak\cong\mathrm{H}_{\{I,I'\}}$ can be extended uniquely to a section $\sigma_{I}$ as $\sigma_{I}(0):= \cK^{+}_{I}$ (or $\sigma_{I}(0):= \cL^{+}_{I}$ if that were the case for $I$). Since exterior wall-crossing leaves all basic building blocks in the pairs~\eqref{eq:vertical} invariant, in the case of an exterior chamber $\qC_{I''}$, the corresponding holomorphic sections can be defined to be the same as the corresponding sections of $\cM_{\qC_{\qI(I'')}}$.
\end{proof}

We will refer to the previous holomorphic sections as the \textit{Hitchin sections} of the geometric model $\cM_{\qC}$. When $d$ is even, three of the Hitchin sections are defined in the Zariski open Harder--Narasimhan stratum of $\cM_{\qC}$, while the fourth one belongs to the complementary stratum. This fourth section corresponds to the holomorphic section constructed by Gothen--Oliveira \cite{GO19}. In~turn, when $d$ is odd, all Hitchin sections always belong to the Zariski open Harder--Narasimhan stratum of $\cM_{\qC}$.

\begin{Corollary}\label{cor:T*CP1}
Under the conventions of Proposition~$\ref{prop:Hitchin-s}$, for every open chamber $\qC\subset \qW_{\sbul}$ there is a biholomorphism
\begin{gather}\label{eq:comp-H-sec}
\cM_{\qC}\big\backslash\bigsqcup \sigma_{I}\bigl(H^{0}\bigl(\CC\PP^{1},K^{2}_{\CC\PP^{1}}(D)\bigr)\bigr) \cong T^{*}\CC\PP^{1}
\end{gather}
mapping the $\CC^{*}$-action into the standard cotangent bundle $\CC^{*}$-action.
The limit $c\to 0$ collapses the Zariski open set
$\cM_{\qC}\backslash\bigsqcup \sigma_{I'}(H^{0}(\CC\PP^{1},K^{2}_{\CC\PP^{1}}(D)))$
to the central sphere, and every component $\mathrm{H}_{\{I,I'\}}$ to the point $\cK^{+}_{I}$ $($or $\cL^{+}_{I})$.
\end{Corollary}

\begin{Remark}
For each interior chamber $\qC_{\qI}$, the points $\{\cN_{I}\colon I\in\qI\}\subset \cS_{\qI}$ correspond to isomorphism classes of stable parabolic bundles supporting nilpotent Higgs fields. Altogether, they are called the \textit{wobbly locus} in the literature. By Corollary~\ref{cor:n-H-blocks}, the induced loci on all central spheres recovers the pair $(\CC\PP^{1},D)$ for any open chamber.
\end{Remark}

Summarizing, we conclude that the projection $\pi_{1}$ for an interior chamber $\qC_{\qI}$ can be defined as an extension of the forgetful map $\pi_{1}([(E_{*},\Phi)]) = [E_{*}]$
over $T^{*}\cS_{\qI}\subset\cM_{\qC_{\qI}}$, under the implicit isomorphism $\cS_{\qI}\cong \CC\PP^{1}$. In virtue of Corollary~\ref{cor:n-H-blocks}, $\pi_{1}$ is also defined for each exterior chamber~$\qC_{I}$ as a projection onto the corresponding exterior central sphere $\cS_{I}\cong\CC\PP^{1}$.
In both cases, the restriction $\pi_{1}\vert_{\pi^{-1}_{2}(0)}$ is the map collapsing the nilpotent cone to its central sphere $\cS_{\qC}$ along the wobbly locus.

\begin{figure}[!ht]
\centering
\includegraphics[width=4.7in]{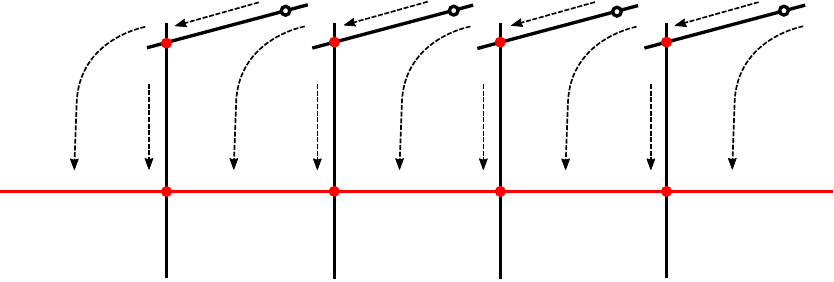}
\caption{The topology of the union of the nilpotent cone and the four Hitchin sections, indicating in red the four 0-dimensional and one 1-dimensional components of fixed points of the global $\CC^{*}$-action in the geometric model, to which $\cM_{\qC}$ collapses in the limit $c\to 0$.}\label{fig:C*-action}
\end{figure}

In view of the restrictions \eqref{eq:branched-cover} and \eqref{eq:comp-H-sec} and Corollary~\ref{cor:T*CP1}, every geometric model $\cM_{\qC}$ is the simultaneous extension of each of the two complex surfaces
\[
T^{*}\CC\PP^{1}\qquad \text{and}\qquad \Sigma_{D}\times\CC^{*}
\]
towards the Hitchin section loci and the nilpotent cone, respectively, in such a way that the intersection of any elliptic fiber of $\pi_{2}$ and the Hitchin section loci corresponds to the set of Weierstrass points of the former. In the case when $\qC$ is an interior chamber and $\cN_{\qC}\neq \varnothing $, $\pi_{1}$ corresponds to the dominant abelianization morphism constructed in \cite{BNR89} for the only possible spectral curve $\Sigma_{D}$, whose very stable locus equals $\cN_{\qB_{\sbul}}$.
The $\CC^{*}$-action on $\cM_{\qC}$ corresponds to the standard $\CC^{*}$-action on both surfaces. The characterization of the limits of this action follows from Propositions~\ref{prop:limits} and \ref{prop:Hitchin-s} (Figure~\ref{fig:C*-action}). There is always a bijective correspondence between $\CC^{*}$-fixed points not in the central sphere and points in the wobbly locus.

\subsection{Weight degeneration to semi-stability walls}

We will finish this section with a brief discussion of an interpretation of the notion of \textit{parabolic $S$-equivalence} for any choice of parabolic weights along the intersection of any given semi-stability wall with the interior of $[0,1]^{4}$. Strictly speaking, it is necessary to exclude the points in $\partial [0,1]^{4}$, as they lead to a degeneration of the moduli problem by definition, to a strictly smaller number of marked points in $\CC\PP^{1}$. Recently, Godinho--Mandini \cite{GM21} have studied the geometry of the special degeneration $\pmb{\beta}=0$ of the Zariski open set of holomorphically trivial parabolic Higgs bundles (in the case when $d=0$), in terms of a correspondence with the so-called \textit{null hyperpolygons}.\looseness=1

{\samepage\begin{Proposition}\label{prop:ss}\qquad
\begin{enumerate}\itemsep=0pt
\item[$(i)$] For any $\pmb{\beta}\in\mathring{\qH}_{I, \vert I\vert - 3}$, all points in the central spheres $\cS_{I}$ and $\cS_{\qI(I)} = \cN_{\qC_{\qI(I)}}$ are strictly semi-stable.
\item[$(ii)$] For any $\pmb{\beta}\in\mathring{\qH}_{I,\vert I\vert -2} = \mathring{\qH}_{I',\vert I'\vert -2}$ and a corresponding pair $\{\qC_{\qI},\qC_{\qI'}\}$ of neighboring interior chambers, all points in the associated pair \eqref{eq:vertical}, as well as the pair $\{\cN_{I},\cN_{I'}\}$, are strictly semi-stable. Additionally, the orbits of $\iota(\mathrm{B}_{I,j}(E)\cap\mathrm{B}_{I',j}(E))\subset \iota(\mathrm{X}(E))$ are strictly semi-stable.
 \end{enumerate}
\end{Proposition}}

\begin{proof}
The result is a straightforward consequence from Lemma~\ref{lemma:not-stable}, Proposition~\ref{prop:C-stable}, and the definition of strict semi-stability for parabolic Higgs bundles.
\end{proof}

\begin{Corollary}[parabolic $S$-equivalence on geometric models]
For any $\pmb{\beta}$ in the interior of a semi-stability wall, the pointwise identification of strictly semi-stable components as in Proposition~$\ref{prop:ss}$ determines a surface biholomorphic to $\cM_{\mathrm{toy}}$.
\end{Corollary}

\appendix

\section{Residue models}\label{sec:Res-mod}

We will describe another approach to explicitly construct moduli spaces of stable parabolic Higgs bundles based on residue evaluations. We will refer to the resulting orbit spaces $\cM^{\mathrm{ev}}_{\qC}(E)$ as \textit{residue models} for Harder--Narasimhan strata.

The desingularization by blow-up of the cone of nilpotent endomorphisms $\mathfrak{n}(V)$ of a 2-dimensional vector space $V$ is constructed as the smooth complex hypersurface
\[
\widetilde{\mathfrak{n}(V)}:= \bigl\{(\phi,F)\in \mathfrak{n}(V)\times \PP(V)\colon F\subset \ker(\phi)\bigr\}.
\]
The restriction to $\widetilde{\mathfrak{n}(V)}$ of the projection onto the second factor is again a projection which we will denote by
$\mathrm{pr}$, and endows $\widetilde{\mathfrak{n}(V)}$ with the structure of a holomorphic line bundle of degree~$-2$, identifying its exceptional divisor $\cE\cong \PP(V)$ with the zero section of the latter. The fiber~$\mathfrak{n}(F)$ over a given line $F\in \PP(V)$ is identified with the Lie algebra of nilpotent endomorphisms of~$V$ whose kernel contains~$F$. $\widetilde{\mathfrak{n}(V)}$ is the fiber product of $\mathfrak{n}(V)\backslash\{0\}$ (viewed as a principal $\CC^{*}$-bundle) under the standard $\CC^{*}$-representation $\rho$ on $\CC$,
\[
\widetilde{\mathfrak{n}(V)} = \bigl(\mathfrak{n}(V)\backslash \{0\}\bigr)\times_{\rho}\CC.
\]
Our goal is to parametrize the spaces $\mathrm{QPH}(E)$ and their stratifications in terms of the geometry of a natural evaluation map associated to $D$. Consider the maps
\[
\mathrm{ev}_{z_{i}}\colon\quad \mathrm{QPH}(E)\rightarrow \widetilde{\mathfrak{n}(E\vert_{z_{i}})},\qquad (E_{\qp},\Phi) \mapsto (\Res_{z_{i}}\Phi, F_{i}), \qquad i=1,2,3,4.
\]
We define the total evaluation to be the map $\mathrm{Ev} = \prod_{i=1}^{4}\mathrm{ev}_{z_{i}}$. Letting
$\mathrm{Pr}$ be the projection map associated to the product of projections $\prod_{i=1}^{4}\mathrm{pr}_{z_{i}}$, we obtain the following commutative diagram
\begin{gather}\label{eq:comm-diag}
\begin{split}&
\xymatrix{
\mathrm{QPH}(E) \ar[rr]^{\mathrm{Ev}} \ar[rd]_{\mathrm{par}} && \displaystyle\prod_{i=1}^{4}\widetilde{\mathfrak{n}(E\vert_{z_{i}})}\ar[ld]^{\mathrm{Pr}}\\
 & \mathrm{QP}(E) &
}\end{split}
\end{gather}
which is equivariant with respect to the induced $\PP(\Aut(E))$-actions. In other words, the total evaluation map is derived from the residue map on parabolic Higgs fields
\[
\Res_{E}\colon\quad \mathrm{Int}(E)\rightarrow \prod_{i=1}^{4}\mathfrak{n}(E\vert_{z_{i}}),\qquad \Phi\mapsto (\Res_{z_{1}}\Phi, \Res_{z_{2}}\Phi, \Res_{z_{3}}\Phi, \Res_{z_{4}}\Phi),
\]
by taking simultaneous blow-ups on its domain and codomain along the loci defined by the identities \eqref{eq:blow-up-loci}.

\begin{Remark}
The surface $\widetilde{\mathfrak{n}(V)}$ can be also described in terms of explicit coordinate charts.\footnote{A choice of isomorphism $V\cong \CC^{2}$ determines an identification $\mathfrak{n}(V)\cong \ZZ_{2}\backslash\CC^{2}$ under the explicit parametrization
\[
(Z_{0},Z_{1}) \mapsto \begin{pmatrix}
Z_{0}Z_{1} & -Z_{0}^{2}\\
Z_{1}^{2} & -Z_{0}Z_{1}
\end{pmatrix}\!,
\]
which is invariant under the $\ZZ_{2}$-action given by $(Z_{0},Z_{1})\mapsto (-Z_{0},-Z_{1})$, and associates the kernel line $[Z_{0} : Z_{1}]\in\PP(\CC^{2})$ whenever $(Z_{0},Z_{1})\neq(0,0)$. Thus, the introduction of the standard affine charts $\cU_{0}$, $\cU_{1}$ on $\PP(\CC^{2})$ induces two affine charts $\cW_{0},\cW_{1}\cong \CC^{2}$ for the blown-up surface $\widetilde{\ZZ_{2}\backslash\CC^{2}}$, defined on their domain intersections in the complement of the exceptional divisor by the maps
\[
(Z_{0}/Z_{1},Z_{1}^{2})\ \reflectbox{\ensuremath{\mapsto}} \
Z_{1}^{2}
\begin{pmatrix}
Z_{0}/Z_{1} & -(Z_{0}/Z_{1})^{2}\\
1 & -Z_{0}/Z_{1}
\end{pmatrix}
=Z_{0}^{2}\begin{pmatrix}
Z_{1}/Z_{0} & -1\\
(Z_{1}/Z_{0})^{2} & -Z_{1}/Z_{0}
\end{pmatrix}\mapsto \big(Z_{1}/Z_{0},Z_{0}^{2}\big).
\]
}
These coordinates determine another way to parametrize $\mathrm{QPH}(E)$ from local trivializations on~$E$, from which a projective model for Harder--Narasimhan strata can be constructed repeating the ideas of this section (cf. \cite{KS19}). We do not pursue that possibility here, since a coordinate-invariant approach is self-sufficient for our general purposes.
\end{Remark}

\begin{figure}[!ht]
\centering
\includegraphics[width=3.0in]{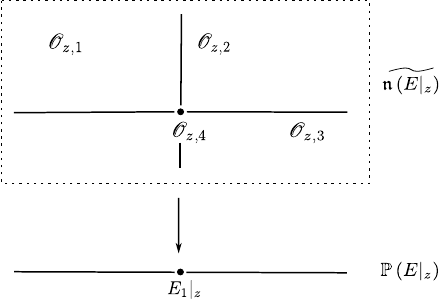}
\caption{Orbit stratification of $\widetilde{\mathfrak{n}(E\vert_{z})}$ when $m_{1} > m_{2}$ $\forall z\in\CC\PP^{1}$.}\label{fig:orbit-strat}
\end{figure}

For every $z\in\CC\PP^{1}$, there is a natural action of the group $\Aut(E)\vert_{z}$ on the space $\mathfrak{n}(E\vert_{z})$ by conjugation. Since the action of central elements is trivial, there is an induced action of $\PP(\Aut(E\vert_{z}))$ on $\widetilde{\mathfrak{n}(E\vert_{z})}$ and $\PP(E\vert_{z})$, making the projection $\mathrm{pr}_{z}$ equivariant. When $m_{1}=m_{2}$ the induced action of $\PP(\Aut(E)|_{z})$ on $\PP(E|_{z})$ is transitive with stabilizer
\[
\PP(\Aut(E|_{z}))_{F} = \PP(\mathrm{P}(F))
\]
for any $F\in \PP(E|_{z})$. Moreover, since the action of $\PP(\Aut(E)\vert_{z})$ on $\mathfrak{n}(E\vert_{z})\backslash\{0\}$ is also transitive, and such that for any $\phi\in \mathfrak{n}(E\vert_{z})\backslash\{0\}$, its stabilizer is equal to the subgroup
\[
\PP(\Aut(E|_{z}))_{\phi} = \PP(\mathrm{R}(F)),
\]
this action partitions the surface $\widetilde{\mathfrak{n}(E\vert_{z})}$ into a total of two orbits, corresponding to the exceptional divisor $\cE_{z}$ and its complement, which we will denote by $\cO_{z}$, that is
\[
\widetilde{\mathfrak{n}(E\vert_{z})} = \cO_{z} \sqcup \cE_{z},\qquad \cO_{z}\cong \mathfrak{n}(E\vert_{z})\backslash\{0\}.
\]
When $m_{1}>m_{2}$, the $\mathrm{pr}_{z}$-inverse image of each of the two orbits $\cV_{z}(E)\sqcup\{E_{1}\vert_{z}\} =\PP(E\vert_{z})$ is in turn partitioned into the two $\PP(\Aut(E)|_{z})$-orbits which result from the intersection with $\cE_{z}$ and its complement $\cO_{z}$.
This way, the action of $\PP(\Aut(E)|_{z})$ on $\widetilde{\mathfrak{n}(E\vert_{z})}$ partitions it into a total of four orbits (Figure~\ref{fig:orbit-strat}), that we will denote by $\cO_{z,1}$, $\cO_{z,2}$, $\cO_{z,3}$ and $\cO_{z,4}$,
where
\[
\cO_{z,1}:=\mathrm{pr}_{z}^{-1}(\cV_{z}(E))\backslash\cE_{z} = \big(\cE_{z}\cup\mathfrak{n}(E_{1}\vert_{z})\big)^{\mathsf{c}}
\]
is Zariski open and a principal homogeneous space for $\PP(\Aut(E)|_{z})$. The remaining orbits
\[
\cO_{z,2}:=\mathfrak{n}(E_{1}\vert_{z})\backslash \cO_{z,4},\qquad \cO_{z,3}:=\cE_{z}\backslash\cO_{z,4},\qquad \cO_{z,4}:=\cE_{z}\cap\mathfrak{n}(E_{1}\vert_{z})
\]
are respectively isomorphic to $\CC^{*}$, $\CC$, and a single point, since the stabilizer of any point $(\phi,E_{1}\vert_{z})\in\cO_{z,2}$ is equal to $\PP(\mathrm{R}(E_{1}\vert_{z}))$.

Given that $\mathrm{Ev}(\mathrm{Q}(E))= \prod_{i=1}^{4}\cE_{z_{i}}$, the inclusion $\iota\colon \mathrm{QP}(E)\rightarrow \mathrm{QPH}(E)$ satisfies $\mathrm{Pr}\circ \mathrm{Ev}\circ\iota = \mathrm{id}$. Since the sum of an arbitrary $\End (E)$-valued holomorphic differential and any parabolic Higgs field is again a parabolic Higgs field, there is an additional affine action of the (possibly trivial) vector space $H^{0}(\CC\PP^{1},\End(E)\otimes K_{\CC\PP^{1}})$ on $\mathrm{QPH}(E)$. The orbit of any $(E_{\qp},\Phi)\in\mathrm{QPH}(E)$ is the level set $\mathrm{Ev}^{-1}(\mathrm{Ev}(E_{\qp},\Phi))$, and $\mathrm{Ev}$ is precisely the invariant map for this action, endowing $\mathrm{QPH}(E)$ with the structure of an (possibly trivial) affine bundle over its image.

\begin{Proposition}\label{prop:Ev}\quad
\begin{enumerate}\itemsep=0pt
\item[$(i)$] If $E$ is evenly-split, then $\mathrm{Ev}$ is a $\PP(\Aut (E))$-equivariant holomorphic embedding, whose image is the smooth 5-dimensional variety of residue constraints.
\item[$(ii)$] When $m_{1}-m_{2} = 2$, $\mathrm{QPH}_{\CC^{*}}(E)\sqcup\mathrm{S}_{m-1}(E)$ and $\mathrm{Q}(E)\sqcup\mathrm{R}(E)$ are $H^{0}(\CC\PP^{1},\End(E)\otimes K_{\CC\PP^{1}})$-invariant, $\mathrm{Ev}(\mathrm{QPH}_{\CC^{*}}(E)\sqcup\mathrm{S}_{m-1}(E))\subset \prod_{i=1}^{4}\cO_{z_{i},1}$ and $\mathrm{Ev}(\mathrm{Q}(E)\sqcup\mathrm{R}(E)) = \prod_{i=1}^{4}\cO^{\mathsf{c}}_{z_{i},1}$.
\item[$(iii)$] When $m_{1}-m_{2} = 3$, $\mathrm{Ev}(\mathrm{R}_{\varnothing }(E)) = \prod_{i=1}^{4}\cO_{z_{i},2}$.
\end{enumerate}
\end{Proposition}

\begin{proof}
$(i)$ $\mathrm{Ev}$ is injective if and only if $E$ is evenly-split. The definition of $\mathrm{Ev}$ implies that, in general, $d\mathrm{Ev}$ has constant rank at any point $(E_{\qp},\Phi)\in\mathrm{QPH}(E)$. Consequently, $\mathrm{Ev}$ is a~holomorphic embedding if $E$ is evenly-split.

($ii$) Since nonzero elements in $H^{0}(\CC\PP^{1},\End(E)\otimes K_{\CC\PP^{1}})$ preserve $E_{1}$, and these exist if and only if $m_{1}-m_{2}\geq 2$, while $\mathrm{QPH}(E) = \mathrm{Q}(E)\sqcup\mathrm{R}(E)$ when $m_{1}-m_{2}\geq 3$, $H^{0}(\CC\PP^{1},\End(E)\otimes K_{\CC\PP^{1}})$-invariance is trivial unless $m_{1}-m_{2}=2$, and then it follows from \eqref{eq:Phi-1}, since $\mathrm{QPH}_{\CC^{*}}(E)\sqcup \mathrm{S}_{0}(E)$ is characterized by $w_{0}\neq 0$, which is unchanged under addition of $\End(E)$-valued holomorphic differentials, and is independent of the isomorphism $E\cong\calO(m+1)\oplus\calO(m-1)$. The second claim is verified from \eqref{eq:Phi-1} in analogy to the proof of Proposition~\ref{prop:QPH-par}$(ii)$ and $(iii)$.

$(iii)$ When $m_{1}-m_{2}=3$ by definition we have that $\mathrm{Ev}(\mathrm{R}_{\varnothing }(E)) = \prod_{i=1}^{4}\cO_{z_{i},2}$.
\end{proof}

\begin{Corollary}
Assume that $\mathrm{QPH}^{s}_{\qC}(E)\neq\varnothing $.
\begin{enumerate}\itemsep=0pt
\item[$(i)$] If $E$ evenly-split, there is an isomorphism $\cM_{\qC}(E)\cong \cM^{\mathrm{ev}}_{\qC}(E)$, where $\cM^{\mathrm{ev}}_{\qC}(E)$ is the $\PP(\Aut(E))$-orbit space of
\[
\mathrm{Ev}(\mathrm{QPH}^{s}_{\qC}(E))\subset \prod_{i=1}^{4}\widetilde{\mathfrak{n}(E\vert_{z_{i}})}.
\]
\item[$(ii)$] When $m_{1}-m_{2}= 2$, $\cM^{\mathrm{ev}}_{\qC}(E):= \mathrm{Ev}(\mathrm{QPH}^{s}_{\qC}(E))/\mathrm{R}(\Aut(E))\cong \cM_{\qC}(E)$.
\item[$(iii)$] When $m_{1}-m_{2}= 3$, $\mathrm{R}(\Aut(E))$ acts freely and transitively on $\mathrm{Ev}(\mathrm{QPH}^{s}_{\qC}(E))$ and
\[
\cM_{\qC}(E)\cong \PP\bigl(H^{0}\bigl(\CC\PP^{1},\End(E)\otimes K_{\CC\PP^{1}}\bigr)\bigr).
\]
\end{enumerate}
\end{Corollary}

\begin{proof}
$(i)$ The claim follows from the equivariance of the commutative diagram \eqref{eq:comm-diag}, Proposition~\ref{prop:Ev}$(i)$, and Theorem~\ref{theo:main-1}.

$(ii)$--$(iii)$ When $m_{1} - m_{2} \geq 2$, the affine $H^{0}(\CC\PP^{1},\End(E)\otimes K_{\CC\PP^{1}})$-action on $\mathrm{QPH}(E)$ is nontrivial and commutes with the action of the subgroup $\mathrm{R}(\Aut(E))$. $\mathrm{R}(\Aut(E))$ acts trivially on $H^{0}(\CC\PP^{1},\End(E)\otimes K_{\CC\PP^{1}})$, and the induced action of $\PP(\Aut(E))/\mathrm{R}(\Aut(E))$ on it equals the square of its standard $\CC^{*}$-action. When $m_{1}-m_{2}=2$, it follows from the classification of nilpotent cone assembly kits and Proposition~\ref{prop:Ev}$(ii)$ that the action of $\PP(\Aut(E))$ on the fibers of $\mathrm{Ev}\vert_{\mathrm{QPH}^{s}_{\qC}(E)}$ is transitive, and the result follows. When $m_{1}-m_{2}=3$ there is a single case when $\mathrm{QPH}^{s}_{\qC}(E)\neq\varnothing $, namely $\qC =\qC_{\varnothing }$, for which $\mathrm{QPH}^{s}_{\qC}(E) = \mathrm{R}_{\varnothing }(E)$. Every element in $\mathrm{R}_{\varnothing }(E)$ can be uniquely expressed as the sum of an element in $\iota(B_{\varnothing ,m+2}(E))$ and a nonzero element in $H^{0}(\CC\PP^{1},\End(E)\otimes K_{\CC\PP^{1}})$. The result follows from Proposition~\ref{prop:Ev}$(iii)$, since $\mathrm{R}(\Aut(E))$ acts freely and transitively on $\prod_{i=1}^{4}\cO_{z_{i},2}$.
\end{proof}

\subsection*{Acknowledgements}

I would like to thank Hartmut Wei\ss, whose encouragement and support where crucial in promp\-ting the appearance of the present manuscript, Steven Rayan for providing many insightful remarks, and the anonymous referee for the careful revision of the manuscript and the constructive criticism provided. This work was supported by the DFG SPP 2026 priority programme ``Geometry at infinity''.

\pdfbookmark[1]{References}{ref}
\LastPageEnding


\begin{thebibliography}{99}
\footnotesize\itemsep=0pt

\bibitem{Bau91}
Bauer S., Parabolic bundles, elliptic surfaces and {${\rm
 SU}(2)$}-representation spaces of genus zero {F}uchsian groups, \href{https://doi.org/10.1007/BF01459257}{\textit{Math.
 Ann.}} \textbf{290} (1991), 509--526.

\bibitem{BNR89}
Beauville A., Narasimhan M.S., Ramanan S., Spectral curves and the generalised
 theta divisor, \href{https://doi.org/10.1515/crll.1989.398.169}{\textit{J.~Reine Angew. Math.}} \textbf{398} (1989), 169--179.

\bibitem{BDHK18}
B\'erczi G., Doran B., Hawes T., Kirwan F., Geometric invariant theory for
 graded unipotent groups and applications, \href{https://doi.org/10.1112/topo.12075}{\textit{J.~Topol.}} \textbf{11}
 (2018), 826--855, \href{https://arxiv.org/abs/1601.00340}{arXiv:1601.00340}.

\bibitem{BJK18}
B\'erczi G., Jackson J., Kirwan F., Variation of non-reductive geometric
 invariant theory, in Surveys in Differential Geometry 2017. {C}elebrating the
 50th Anniversary of the {J}ournal of {D}ifferential {G}eometry, \textit{Surv.
 Differ. Geom.}, Vol.~22, \href{https://doi.org/10.4310/SDG.2017.v22.n1.a2}{Int. Press}, Somerville, MA, 2018, 49--69,
 \href{https://arxiv.org/abs/1712.02576}{arXiv:1712.02576}.

\bibitem{Bis98}
Biswas I., A criterion for the existence of a parabolic stable bundle of rank
 two over the projective line, \href{https://doi.org/10.1142/S0129167X98000233}{\textit{Internat.~J. Math.}} \textbf{9} (1998),
 523--533.

\bibitem{Bla15}
Blaavand J.L., The {D}irac--{H}iggs bundle, Ph.D.~Thesis, {U}niversity of
 Oxford, 2015.

\bibitem{FMSW22}
Fredrickson L., Mazzeo R., Swoboda J., Weiss H., Asymptotic geometry of the
 moduli space of parabolic $\mathrm{SL}(2,\mathbb{C})$-{H}iggs bundles,
 \href{https://doi.org/10.1112/jlms.12581}{\textit{Proc. London Math. Soc.}}, {t}o appear, \href{https://arxiv.org/abs/2001.03682}{arXiv:2001.03682}.

\bibitem{GMN13}
Gaiotto D., Moore G.W., Neitzke A., Wall-crossing, {H}itchin systems, and the
 {WKB} approximation, \href{https://doi.org/10.1016/j.aim.2012.09.027}{\textit{Adv. Math.}} \textbf{234} (2013), 239--403,
 \href{https://arxiv.org/abs/0907.3987}{arXiv:0907.3987}.

\bibitem{GM13}
Godinho L., Mandini A., Hyperpolygon spaces and moduli spaces of parabolic
 {H}iggs bundles, \href{https://doi.org/10.1016/j.aim.2013.04.026}{\textit{Adv. Math.}} \textbf{244} (2013), 465--532,
 \href{https://arxiv.org/abs/1101.3241}{arXiv:1101.3241}.

\bibitem{GM21}
Godinho L., Mandini A., Quasi-parabolic {H}iggs bundles and null hyperpolygon
 spaces, \href{https://doi.org/10.1090/tran/8450}{\textit{Trans. Amer. Math. Soc.}} \textbf{374} (2021), 7411--7447,
 \href{https://arxiv.org/abs/1907.01937}{arXiv:1907.01937}.

\bibitem{GO19}
Gothen P.B., Oliveira A.G., Topological mirror symmetry for parabolic {H}iggs
 bundles, \href{https://doi.org/10.1016/j.geomphys.2018.08.020}{\textit{J.~Geom. Phys.}} \textbf{137} (2019), 7--34,
 \href{https://arxiv.org/abs/1707.08536}{arXiv:1707.08536}.

\bibitem{Ham19}
Hamilton E., Stratifications and quasi-projective coarse moduli spaces for the
 stack of {H}iggs bundles, \href{https://arxiv.org/abs/1911.13194}{arXiv:1911.13194}.

\bibitem{Hau98}
Hausel T., Compactification of moduli of {H}iggs bundles, \href{https://doi.org/10.1515/crll.1998.096}{\textit{J.~Reine
 Angew. Math.}} \textbf{503} (1998), 169--192, \href{https://arxiv.org/abs/math.AG/9804083}{arXiv:math.AG/9804083}.

\bibitem{HH16}
Heller L., Heller S., Abelianization of {F}uchsian systems on a 4-punctured
 sphere and applications, \href{https://doi.org/10.4310/JSG.2016.v14.n4.a4}{\textit{J.~Symplectic Geom.}} \textbf{14} (2016),
 1059--1088, \href{https://arxiv.org/abs/1404.7707}{arXiv:1404.7707}.

\bibitem{HL19}
Heu V., Loray F., Flat rank two vector bundles on genus two curves,
 \href{https://doi.org/10.1090/memo/1247}{\textit{Mem. Amer. Math. Soc.}} \textbf{259} (2019), v+103~pages,
 \href{https://arxiv.org/abs/1401.2449}{arXiv:1401.2449}.

\bibitem{Hit87a}
Hitchin N.J., The self-duality equations on a {R}iemann surface, \href{https://doi.org/10.1112/plms/s3-55.1.59}{\textit{Proc.
 London Math. Soc.}} \textbf{55} (1987), 59--126.

\bibitem{KW18}
Kim S., Wilkin G., Analytic convergence of harmonic metrics for parabolic
 {H}iggs bundles, \href{https://doi.org/10.1016/j.geomphys.2018.01.023}{\textit{J.~Geom. Phys.}} \textbf{127} (2018), 55--67,
 \href{https://arxiv.org/abs/1705.08065}{arXiv:1705.08065}.

\bibitem{KS19}
Komyo A., Saito M.H., Explicit description of jumping phenomena on moduli
 spaces of parabolic connections and {H}ilbert schemes of points on surfaces,
 \href{https://doi.org/10.1215/21562261-2019-0016}{\textit{Kyoto~J. Math.}} \textbf{59} (2019), 515--552, \href{https://arxiv.org/abs/1611.00971}{arXiv:1611.00971}.

\bibitem{LS15}
Loray F., Saito M.H., Lagrangian fibrations in duality on moduli spaces of rank
 2 logarithmic connections over the projective line, \href{https://doi.org/10.1093/imrn/rnt232}{\textit{Int. Math. Res.
 Not.}} \textbf{2015} (2015), 995--1043, \href{https://arxiv.org/abs/1302.4113}{arXiv:1302.4113}.

\bibitem{LSS13}
Loray F., Saito M.H., Simpson C.T., Foliations on the moduli space of rank two
 connections on the projective line minus four points, in Geometric and
 Differential {G}alois Theories, \textit{S\'emin. Congr.}, Vol.~27, Soc. Math.
 France, Paris, 2013, 117--170, \href{https://arxiv.org/abs/1012.3612}{arXiv:1012.3612}.

\bibitem{MS80}
Mehta V.B., Seshadri C.S., Moduli of vector bundles on curves with parabolic
 structures, \href{https://doi.org/10.1007/BF01420526}{\textit{Math. Ann.}} \textbf{248} (1980), 205--239.

\bibitem{Men18}
Meneses C., Remarks on groups of bundle automorphisms over the {R}iemann
 sphere, \href{https://doi.org/10.1007/s10711-017-0309-y}{\textit{Geom. Dedicata}} \textbf{196} (2018), 63--90,
 \href{https://arxiv.org/abs/1607.03865}{arXiv:1607.03865}.

\bibitem{Mena}
Meneses C., Geometric models for moduli of rank~2 parabolic {H}iggs bundles in
 genus~0 and applications, {i}n preparation.

\bibitem{MT21}
Meneses C., Takhtajan L.A., Logarithmic connections, {WZNW} action, and moduli
 of parabolic bundles on the sphere, \href{https://doi.org/10.1007/s00220-021-04183-y}{\textit{Comm. Math. Phys.}} \textbf{387}
 (2021), 649--680, \href{https://arxiv.org/abs/1407.6752}{arXiv:1407.6752}.

\bibitem{Mir89}
Miranda R., The basic theory of elliptic surfaces, \textit{Dottorato di Ricerca in
 Matematica}, ETS Editrice, Pisa, 1989.

\bibitem{Muk03}
Mukai S., An introduction to invariants and moduli, \textit{Cambridge Studies
 in Advanced Mathematics}, Vol.~81, \href{https://doi.org/10.1017/CBO9781316257074}{Cambridge University Press}, Cambridge,
 2003.

\bibitem{Nak96}
Nakajima H., Hyper-{K}\"ahler structures on moduli spaces of parabolic {H}iggs
 bundles on {R}iemann surfaces, in Moduli of Vector Bundles ({S}anda, 1994;
 {K}yoto, 1994), \textit{Lecture Notes in Pure and Appl. Math.}, Vol.~179,
 Dekker, New York, 1996, 199--208.

\bibitem{Ray17}
Rayan S., The quiver at the bottom of the twisted nilpotent cone on
 {$\mathbb{P}^1$}, \href{https://doi.org/10.1007/s40879-016-0120-6}{\textit{Eur.~J. Math.}} \textbf{3} (2017), 1--21,
 \href{https://arxiv.org/abs/1609.08226}{arXiv:1609.08226}.

\bibitem{RS21}
Rayan S., Schaposnik L.P., Moduli spaces of generalized hyperpolygons,
 \href{https://doi.org/10.1093/qmath/haaa036}{\textit{Q.~J.~Math.}} \textbf{72} (2021), 137--161, \href{https://arxiv.org/abs/2001.06911}{arXiv:2001.06911}.

\bibitem{Sim90}
Simpson C.T., Harmonic bundles on noncompact curves, \href{https://doi.org/10.2307/1990935}{\textit{J.~Amer. Math.
 Soc.}} \textbf{3} (1990), 713--770.

\bibitem{Sim92}
Simpson C.T., Products of matrices, in Differential Geometry, Global Analysis,
 and Topology ({H}alifax, {NS}, 1990), \textit{CMS Conf. Proc.}, Vol.~12,
 Amer. Math. Soc., Providence, RI, 1991, 157--185.

\bibitem{Thad02}
Thaddeus M., Variation of moduli of parabolic {H}iggs bundles, \href{https://doi.org/10.1515/crll.2002.051}{\textit{J.~Reine
 Angew. Math.}} \textbf{547} (2002), 1--14, \href{https://arxiv.org/abs/math.AG/0003222}{arXiv:math.AG/0003222}.

\end{thebibliography}
\end{document}